\documentclass[12pt]{amsart}
\usepackage{amssymb,amscd}
\usepackage{verbatim}

\usepackage{graphicx,color}

\usepackage{graphicx}

\usepackage{here}   

\newtheorem*{Whitney towers}{Theorem~\ref{Whitney towers}}
\newtheorem*{h-towers}{Theorems ~\ref{half} \& \ref{$(n)$-solvable}}

\newtheorem*{surgery curves}{Theorem~\ref{surgery curves}}
\newtheorem*{cg=0}{Theorem~\ref{vanish}}

\newtheorem{thm}{Theorem}[section]

\newtheorem{pr}[thm]{Proposition} 
\newtheorem{lem}[thm]{Lemma}

\newtheorem{cla}[thm]{Claim}

\theoremstyle{definition}
\newtheorem{defn}[thm]{Definition}
\newtheorem{note}[thm]{Note}
\newtheorem{que}[thm]{Question}

\numberwithin{equation}{section}
\numberwithin{figure}{section}

\newcommand{\x}{\times}
\newcommand{\np}{\newpage}
\newcommand{\Z}{\mathbb{Z}}
\newcommand{\N}{\mathbb{N}}

\newcommand{\R}{\mathbb{R}}

\def\yen{{\setbox0=\hbox{Y}Y\kern-.97\wd0\vbox{hrule height.lex width.98%
\wd0\kern.33ex\hrule height.lex width.98\wd0\kern.45ex}}}

\def\np{\newpage}

\begin{document}
\pagestyle{plain}

\title{
A spinning construction for virtual 1-knots and 2-knots,   
and the fiberwise and welded equivalence of virtual 1-knots 
}
\author{Louis H. Kauffman,  Eiji Ogasa, and Jonathan Schneider}


\date{}
\begin{abstract}   
Spun-knots  (respectively, spinning tori) in $S^4$ 
made from classical 1-knots compose an  important class of 
spherical 2-knots (respectively, embedded tori) contained in $S^4$. 
Virtual 1-knots are generalizations of classical 1-knots.  
We generalize these constructions to the virtual 1-knot case  
by using what we call, in this paper, the spinning construction of a submanifold.
The construction proceeds as follows: 
It has been known that 
there is a consistent way to make an embedded circle $C$ contained in \\
(a closed oriented surface $F$)$\x$(a closed interval $[0,1]$) from any virtual 1-knot $K$. 
%
%
%
Embed $F$ in $S^4$ by an embedding map $f$.  
Let $F$ also denote $f(F).$
We can regard the tubular neighborhood of $F$ in $S^4$ 
as $F\x D^2$. 
Let $[0,1]$ be a radius of $D^2$. 
We can regard $F\x D^2$ 
as the result of rotating 
 $F\x [0,1]$ around $F\x \{0\}$.    
Assume $C\cap(F\x\{0\})=\phi$. 
Rotate $C$ together 
when we rotate $F\x [0,1]$ around $F\x \{0\}$. 
Thus we obtain  an embedded torus $Q\subset S^4$.  
We prove the following: 
The embedding type $Q$ in $S^4$ depends only on $K$, and does not depend on $f$. 
Furthermore, 
the submanifolds,   $Q$ and  the embedded torus made from $K,$ defined by Satoh's method, 
of $S^4$ are isotopic.   

We generalize this construction in the virtual 1-knot case, and 
we also succeed to make a consistent 
construction of one-dimensional-higher tubes from 
any virtual 2-dimensional knot.
Note that Satoh's method says nothing about the virtual 2-knot case. 
Rourke's interpretation of Satoh's method is that 
one puts `fiber-circles' on each point of each virtual 1-knot diagram. 
If there is no virtual branch point in a virtual 2-knot diagram, 
our way gives such fiber-circles to each point of the  virtual 2-knot diagram. 
Furthermore we prove the following:
If a virtual 2-knot diagram $\alpha$ has a virtual branch point,  
$\alpha$ cannot be covered by such fiber-circles.

We obtain a new equivalence relation, 
 the $\mathcal E$-equivalence relation 
of the set of virtual 2-knot diagrams, 
that is much connected with 
the welded equivalence relation and our spinning construction. 
We prove that 
there are virtual 2-knot diagrams, $J$ and $K$, 
that  are virtually nonequivalent 
but  are $\mathcal E$-equivalent. 

Although Rourke claimed that 
two virtual 1-knot diagrams $\alpha$ and $\beta$ are 
fiberwise equivalent if and only if  
$\alpha$ and $\beta$ are welded equivalent,  
we state that this claim is wrong. 
We prove that 
two virtual 1-knot diagrams $\alpha$ and $\beta$ are 
fiberwise equivalent if and only if  
$\alpha$ and $\beta$ are rotational welded equivalent (the definiton of rotational welded equivalence is given in the body of the paper). 
\end{abstract}   

\maketitle

\np
\tableofcontents

\section{Introduction}\label{jobun}
\subsection{
Spinning tori
}\label{i1}\hskip20mm\\%
Spun-knots (respectively, spinning tori) in $S^4$ 
made from classical 1-knots compose an  important class of 
spherical 2-knots (respectively, embedded  tori) contained in $S^4$. 
See \cite{Zeeman} for the definition of spun-knots. 
We review the construction of them below.\\
 
Let $\R^4=\{(x,y,z,w)|x,y,z,w\in\R\}$. 
Regard 
$\R^4$ as the result of rotating 
$H=\{(x,,y,z,w)|x\geqq0,w=0\}$ 
around 
$A=\{(x,,y,z,w)|x=0,w=0\}$ 
as the axis. 
Take a 1-knot $K$ in $H$ 
so that $K\cap A$ is  an arc (respectively, the empty set).  
Rotate $K-{\rm Int}(K\cap A)$ 
around $A$ together 
when we rotate  $H$ around $A$. 
The resultant submanifold of $\R^4$
is the spun knot  (respectively, the spinning tori) of $K$. 
We can easily also regard them as submanifolds of $S^4$. 
We can define a spun link if $K$ is a link although we discuss the knot case mainly in this paper. 
Our discussion can be easily generalized to the link case. 

One of our themes in this paper is to generalize the spun knots of classical knots to the virtual knot case.
We begin by explaining why virtual knots are important.

\bigbreak
\subsection{
History of relations between virtual knots and QFT,  
and a reason why virtual knots are important 
}\label{i2}\hskip20mm\\%
Virtual 1-links are defined in \cite{Kauffman1,Kauffman, Kauffmani} as generalizations 
of classical 1-links.  One motivation for virtual 1-links is as follows. 
 Jones \cite{Jones}  defined the Jones polynomial for classical 1-links in $S^3$. 
The following had been well-known before the Jones polynomial was found:
The Alexander polynomial for classical 1-links in $S^3$
 is defined in terms of the topology of the complement of the link and can be generalized to 
give invariants of closed oriented 3-manifolds and of links within the 3-manifold.\\

Jones \cite[page 360, \S10]{Jones} tried to define a 3-manifold invariant associated with the Jones polynomial, 
and succeeded in some cases. 
Of course, when the Jones polynomial was found, the following question was regarded as a very natural one:  

\smallbreak 
{\bf Question J.}  Can we generalize the definition of the Jones polynomial for classical 1-links in $S^3$ to that in any 3-manifold?  
\smallbreak

Note the result may not be a polynomial but a function of $t$. 
\\

Witten \cite{W} wrote a quantum field theoretic path integral for any 1-link $L$ in any compact oriented 3-manifold $M$. 
His path integral included the Jones polynomial for 1-links in $S^3$, 
its generalizations and new (at the time) invariants of 3-manifolds.
This was a breakthrough for the philosophy of physics in that 
one of the most natural geometrically intrinsic interpretations of a mathematical object 
was done by using a path integral, and had not been done by any other way. 

\smallbreak
\noindent
{\bf Note.} 
Here, `geometrically intrinsic interpretation' means the point of view that would define a link invariant in terms of the embedding of the link in the ambient 3-dimensional manifolds  
just as one can do naturally and easily in the case of the Alexander polynomial of 1-knots. 
Jones \cite{Jones} defined the Jones polynomial by using representations of braid groups to an operator algebra (the Temperley-Lieb algebra). Representations, braid groups, operator algebras are mathematically explicit objects so some people may feel that that  is enough to consider the meaning of the Jones polynomial. 
\bigbreak

If $M=S^3$, we can say at the physics level that the Witten path integral represents the Jones polynomial 
for 1-links in $S^3$. 
Reshetikhin, Turaev, Lickorish and others \cite{KM, Lickorish, Lickorishl, RT} etc.  
generalized the result in \cite[page 360, \S10]{Jones} and   
created rigorous definitions for invariants of 3-manifolds that parallel Witten's ideas, without using the functional integral.  
They succeeded to define new invariants of closed oriented 3-manifolds and 
invariants of links embedded in 3-manifolds
that we today call quantum invariants.  
(Note, here, we distinguish the above invariants of links embedded in 3-manifolds
with the Jones polynomial for them as below.) 
In both Witten's version and the Reshetikhin-Turaev versions the invariants of 3-manifolds are obtained by representing
the 3-manifold as surgery on a framed link and summing over invariants corresponding to appropriate representations decorating the surgery link. The same technique
applies when one includes an extra link component that is not part of the surgery data. 
In this way, one obtains quantum invariants of links in 3-manifolds.
Another technique, formalized by Crane \cite{Crane} and by Kohno \cite{Kohno} uses a Heegard decomposition of the 3-manifold and algebraic structure of the conformal field theory 
for the surface of the Heegard decompositon. These methods produce invariants for 3-manifolds and, in principle, invariants for links in 3-manifolds, but  are much more indirect than the original physical idea of Witten that would integrate directly over the many possible evaluations of the Wilson loop for the knot or link in the 3-sphere, or the original combinatorial 
skein techniques that produce the invariant of a link from its diagrammatic combinatorics.
%
%
%
See \cite{Kauffmanp}. \\

The Witten path integral is written also in the case where $L\neq\phi$ and $M\neq S^3$. 
It corresponded to 
Question J , which had been considered before the Witten path integral appeared.    
\bigbreak

In \cite{KauffmanJ} Kauffman found a definition of the Jones polynomial as a state summation over combinatorial states of the link diagram and found a diagrammatic interpretation of the Temperley-Lieb algebra that put the original definition of Jones in a wider context of generalized partition functions and statistical mecnanics on graphs and knot and link diagrams.
In \cite{Kauffman1,Kauffman, Kauffmani} Kauffman generalized the Jones polynomial
in the case where $M$ is (a closed oriented surface)$\x[-1,1]$. 
In fact, 
\cite{Kauffman1,Kauffman, Kauffmani} defined virtual 1-links 
as 
another way of describing 
1-links in  (a closed oriented surface)$\x[-1,1]$: 
the set of virtual 1-links 
is the same as 
that of 1-links in  (a closed oriented surface)$\x[-1,1],$ taken up to handle stabilization.
See Theorem \ref{vk}.  We make the point here that the virtual knot theory is a context for links in the fundamental 3-manifolds of the form $F \times I$ where $F$ is a closed 
surface. The state summation approach to the Jones polynomial generalizes to invariants of links in such thickened surfaces. This provides a significant and direct arena for 
examining such structures without the functional integral. It also provides challenges for corresponding approaches that use the functional integral methods. It remains a serious 
challenge to produce ways to work with the functional integrals that avoid difficulties in analysis.
\\

Path integrals represent the superposition principle dramatically. 
This is a marvelous idea of Feynman. 
The Witten path integral also represents a geometric idea of the Jones polynomial and quantum invariants physics-philosophically very well. 
Witten found a Lagrangian via the Chern-Simons 3-form and Wilson line with a tremendous insight, and he calculated  the path integral of the Lagrangian rigorously at physics level,  
and showed that the result of the calculation is the Jones polynomial for links in $S^3$, and the quantum invariants of any closed oriented 3-manifold with or without embedded circles.  It is a wonderful work of Witten. 
However recall the following facts:
The Witten path integral for any 1-link in any closed 3-manifold has not been calculated 
in mathematical level nor in physics level in any way that can be regarded as direct. This means that Question J is open in the general case.
%
That is, 
nobody has succeeded to generalize the Jones polynomial in a direct way, and mathematically rigorously to the case where $M$ is not $S^3$, (respectively, $B^3$, $\R^3$), nor 
(a closed oriented surface $F$)$\x[-1,1]$.  (Note the last manifold is not closed. 
Note that the discussion in the $S^3$ case is the same as that in  the $B^3$ (respectively, $\R^3$) case. Virtual knot theory can also discuss the case where $F$ is compact and non-closed,  but then we need to fix the embedding type of $F$ in $F\x[-1,1]$.)\\

Recall the following fact: Even if we make a (seemingly) meaningful Lagrangian,   
the path integral associated with the Lagrangian cannot always be calculated.
An example is 
the Witten path integral associated with 
the general case of Question J. 
Another one is the following. Today they do not know how to calculate the path integral 
if we replace Chern-Simons-3-form on 3-manifolds with 
Cern-Simons-$(2p+1)$-form on $(2p+1)$-manifolds, 
where $p$ is any integer$\geq2$, 
in the Witten path integral.  
Indeed nowadays they only calculate path integrals only when they can calculate them.
If the path integral of the Lagrangian  is not calculated explicitly, neither mathematicians nor physicists 
regard the theory of the Lagrangian as a meaningful one.  
Furthermore, even if we calculate path integrals,  the result of the calculation is sometimes what we do not expect. 
See an example of \cite{LeeYang} explained in 
\cite[the last part of section 5.1]{Ryder}.\\

The heuristics of the Witten path integral have not been fully mined. See \cite{KauffmanPath} for a survey of the results of some of these heuristics in relation to the Jones
polynomial and Vassiliev invariants. It is possible that good heuristics will emerge for understanding invariants of links in 3-manifolds. But at the present time it is worth examining the 
cases we do understand for working with generalizations of the Jones polynomial for links in thickened surfaces. 
We had begun considering Question J 
before the Witten path integral 
appeared in this discussion.
Question J is also natural and important 
even if we do not consider path integrals. 
\\


\noindent
{\bf Note.} 
(1) 
We can observe some historical correspondences. 
Feynman discovered path integrals 
by using an analogy with (quantum) statistical mechanics, 
and he interpreted quantum theory by using path integrals.  
Operator algebras,  path integrals, (quantum) statistical mechanics are closely related. 
The Jones polynomial is discovered by using operator algebras (\cite{Jones}), next 
is interpreted via (quantum) statistical mechanics (\cite{KauffmanJ}), then by using path integrals (\cite{W}).
Operator algebras,  path integrals, and (quantum) statistical mechanics are related again with topology in 
the background.

\smallbreak\noindent
(2)
The Jones polynomial of 1-links in (a closed oriented surface)$\x$(the interval) 
is discovered in \cite{Kauffman1,Kauffman, Kauffmani},  
by  using the analogy with state sums in (quantum) statistical mechanics in \cite{KauffmanJ}.

\smallbreak\noindent
(3) \cite{KauffmanSaleur} found a relation between the Alexander-Conway polynomial between 
1-dimensional classical knots and quantum field theory. The relation gives a different aspect from the Homflypt polynomial and the Witten path integral. 
\cite{Ogasapath} found a relation between the degree of the Alexander polynomial of high dimensional knots and the Witten index of a supersymmetric quantum system. 
It is also an outstanding open question whether we can define an analog to the Jones polynomial for high dimensional knots.
\bigbreak

Virtual 1-links have many other important properties than the above one. 
See  \cite{Kauffman1,Kauffman, Kauffmani}.  
Thus it is very natural to consider whether any property of classical 1-knots is possessed by virtual 1-knots, as below.

\bigbreak
\subsection{
Main results
}\label{i3}\hskip20mm\\%
We generalize the construction of spun-knots (respectively, spinning tori) of classical 1-knots  
to the virtual 1-knot case  as follows. 
Recall that, in \cite{Kauffman1, Kauffman, Kauffmani} there is given  a consistent way 
to make an embedded circle $C$ contained in 
(a closed oriented surface $F$)$\x$(a closed interval $[0,1]$) from any virtual 1-knot $K$ diagram 
(see Theorem \ref{vk}).  
Note the following. 
When we construct spun knots (spinning tori), 
we regard $\R^4$ itself as the total space of the normal bundle of $A$ in $\R^4$.
Recall that $A$ is defined in \S\ref{i1}. 
Embed $F$ in $\R^4\subset S^4$ by an embedding map $f$.  
Let $F$ stand for $f(F).$.
Note that the tubular neighborhood of $F$ in $S^4$ is diffeomorphic to $F\x D^2$. 
Let $[0,1]$ be a radius of $D^2$. 
We can regard $F\x D^2$ 
as the result of rotating 
 $F\x [0,1]$ around $F\x \{0\}$.    
Assume $C\cap(F\x\{0\})=\phi$. 
Rotate $C$ together 
when we rotate $F\x [0,1]$ around $F\x \{0\}$. 
Thus we obtain  an embedded torus $Q\subset S^4$.  \\

We prove the following  (Theorems \ref{honto} and \ref{mainkore}): 
The embedding type $Q$ in $S^4$ depends only on $K$, and does not depend on $f$. 
Furthermore the submanifolds, 
$Q$ and the embedded torus made from $K$ defined by Satoh in \cite{Satoh}, 
of $S^4$ are isotopic.\\ 

This construction of $Q$ is an example of what we call the spinning construction of submanifolds 
in Definition \ref{spinningsubmanifold}. 
This paper does not discuss the case where $C\cap(F\x\{0\})\neq\phi$. 
\\

There are classical 1-knots, virtual 1-knots, and classical 2-knots 
so it is natural to consider virtual 2-knots. We define virtual 2-knots in Definition \ref{JV}.  
It is very natural to consider whether 
any property of `classical 1-, and 2-knots and virtual 1-knots' 
is possessed by virtual 2-knots. 
It is natural to ask whether we can define one-dimensional-higher tubes
for  virtual 2-knots 
(Question \ref{North Carolina}) 
since
we succeed in the virtual 1-knot case 
as explained above.
Note that Satoh's mehtod in \cite{Satoh} does not treat the virtual 2-knot case.  \\

In the virtual 1-knot case, 
in \cite{Rourke},  
Rourke  interpreted Satoh's method  as follows:  
Let $\alpha$ be any virtual 1-knot diagram.  
Put `fiber-circles' on each point of $\alpha$ and obtain a one-dimensional-higher tube.  
(We review this construction in Theorem \ref{Montana} and Definition \ref{Nebraska}).   
If we try to generalize Rourke's way to the virtual 2-knot case, 
we encounter the following situation.\\
%

Let $\alpha$ be any virtual 1-knot diagram.  There are two cases:

\smallbreak\noindent(1) 
The case where $\alpha$ has no virtual branch point.  
(We define virtual branch point in Definitions \ref{oyster} and \ref{JV}.)

\smallbreak\noindent(2) 
The case where $\alpha$ has a virtual branch point.  

\smallbreak
In the  case (1) , we can make a tube by Rourke's method. 
See \cite[section 3.7.1]{J}, Note \ref{kaiga}, and Definition \ref{suiri}.   
In the case (2), however, Schneider \cite{J} found 
it difficult to define a tube near any virtual branch point. \\

Thus we consider the following two problems. \\

Can we put fiber-circles  over each point of any virtual 2-knot 
in a consistent way as described above, 
and make a one-dimensional-higher tube
(Question \ref{North Dakota})? \\

Is there a one-dimensional-higher tube construction which is defined for all virtual 2-knots, and which agrees with the method in the case (1)  written above   
when there are no virtual branch points  
(Question \ref{South Dakota})? \\

In Theorem \ref{vv} we give an affirmative answer to Question \ref{South Dakota}. 
Our solution is a generalization of our method in the virtual 1-knot case 
used in \S\S\ref{E}-\ref{Proof}. 
We also use 
a spinning construction of submanifolds 
explained in Definition \ref{spinningsubmanifold}. \\

In Theorem \ref{Rmuri} 
we give a negative answer to Question \ref{North Dakota}. \\
\\

We obtain a new equivalence relation, 
 the $\mathcal E$-equivalence relation of 
the set of virtual 
1- and 2-knot diagrams (Definition \ref{zoo}). 
It is done by using the above spinning construction. 
 The $\mathcal E$-equivalence relation is closley connected with 
the welded equivalence relation and our spinning construction. 
Welded 1-links are defined in \cite{Rourke} associated with virtual 1-links. 
Welded 1-links are related to tubes very much as we discuss in this paper.  
We introduce welded 2-knots in Definition \ref{JW}. \\

We prove that there are virtual 2-knot diagrams, $J$ and $K$, 
that are virtually nonequivalent 
but are $\mathcal E$-equivalent 
(Theorem \ref{Maine}).      
Welded 1-,and 2-knots are recipients  of the tube construction or the above spinning construction.  
The above spinning construction is related to the fiberwise equivalence explained below. 
We will explain their connection in this paper and this is a theme of this research.\\

Although Rourke claimed in \cite[Theorem 4.1]{Rourke}  that 
two virtual 1-knot diagrams $\alpha$ and $\beta$ are fiberwise equivalent if and only if  
$\alpha$ and $\beta$ are welded equivalent,  
we state that this claim is wrong. 
(See 
\cite{Rourke} and Definitions \ref{Nevada} 
of this paper for the definition of the fiberwise equivalence, 
and 
 \cite{Rourke, Satoh} for that of the welded equivalence.) 
The reason for the failure of Rourke's claim is given in   
Theorems \ref{smooth} and  \ref{fwrw}, and Claim \ref{panda}. 
We prove in Theorems \ref{smooth} and \ref{fwrw} that {\it virtual 1-knot diagrams, $\alpha$ and $\beta$, 
are    fiberwise equivalent if and only if  they are rotational welded equivalent.} The reader can recall that in virtual 1-knot theory there are Reidemeister-type moves for virtual 
crossings. Rotational equivalence for virtual knots is obtained by making the virtual curl (analog of the first Reidemeister move) forbidden. Rotational equivalence for welded knots 
also forbids the virtual curl move in the context of the rules for welded knots.
(See \cite{Kauffman, Kauffmanrw, J}  for rotational welded equivalence.)  
Our  result is proved by using the property of virtual 2-knots found in Theorem \ref{Rmuri}.
Virtual 2-knots themselves are important,  and furthermore they are also important for research in virtual 1-knots. 
Our main results are 
Theorems \ref{honto}, \ref{mainkore}, 
\ref{vv},  \ref{Rmuri}, 
\ref{Maine}, 
 \ref{smooth},  
and \ref{Montgomery}. \\

\bigbreak
\section{$\mathcal K(K)$ for a virtual 1-knot $K$}\label{K}
\noindent 
We work in the smooth category unless we indicate otherwise. 
In a part of \S\ref{New Mexico}
we will use the PL category in order to prove our results in the smooth category. 
See Note \ref{haruwa}.
We review some facts on virtual 1-knots in this section before we state two of our main results, 
Theorems \ref{honto} and \ref{mainkore}, 
in the following section. \\

\begin{figure}
     \includegraphics[width=140mm]{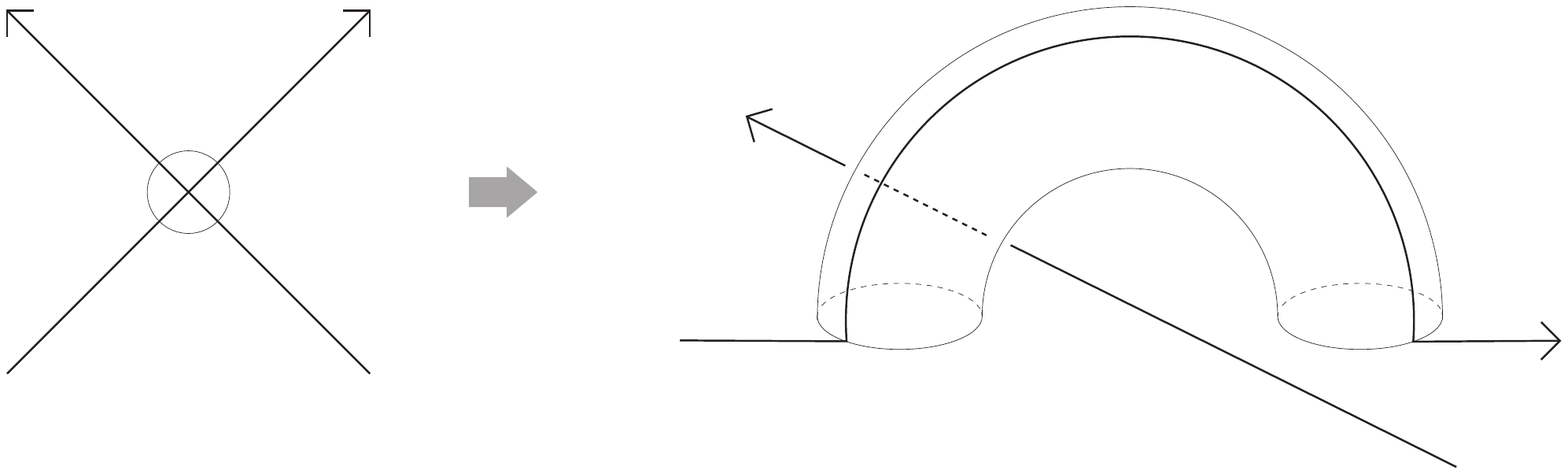}  
\vskip-50mm
\caption{{\bf A virtual crossing point and a surgery by a 1-handle}\label{Alabama}}   
\end{figure}

Let $\alpha$ be a virtual 1-knot diagram. 
In this paper we use Greek lowercase letters for virtual diagrams and 
Roman capital letters for virtual knots. 
See \cite{Kauffman1, Kauffman, Kauffmani} 
for the definition and properties of virtual 1-knot diagrams and 
those of virtual 1-knots. 
For $\alpha$ there are a nonnegative integer $g$ and 
an embedded circle contained in $\Sigma_g\x[0,1]$ as follows, 
where $\Sigma_g$ is a closed oriented surface with genus $g$. 
Take $\alpha$ in $\R^2$. (Recall that we can make the infinity point $\{*\}$ 
and $\R^2$ into $S^2$.)  
Carry out a surgery on $\R^2$ 
by using a 3-dimensional 1-handle near a virtual crossing point as shown in 
Figure \ref{Alabama}   
and obtain $T^2-\{*\}$. 
Note that the virtual 1-knot $K$ is oriented and that the arrows in 
Figure \ref{Alabama}   
 denote the orientation.
Segments are changed as shown in the right figure of Figure \ref{Alabama}.
Do this procedure near all virtual crossing points.  
Suppose that $\alpha$ has $g$ copies of virtual crossing point ($g\in\N\cup\{0\}$). 
Here, $\N$ denotes the set of natural numbers. 
Note that $a$ is a natural number if and only if $a$ is a positive integer.   
What we obtain is $\Sigma_g-\{*\}$.  We call it $\Sigma_g^\bullet$.  
(In \S\ref{Proof}, for a closed oriented surface $F$, 
we define $F^\circ$ to be $F-$(an open 2-disc). 
So, here, we use $^\bullet$  not  $^\circ$.)
Thus we obtain an immersed circle in $\Sigma_g^\bullet$ from $\alpha$. 
Call it $\mathcal I(\alpha)$. 
Note that it is an immersion in ordinary sense (that is, it has no `virtual crossing point').  
Regard  $\Sigma_g$ as an abstract manifold. 
Make $\Sigma_g^\bullet\x[0,1]$.  
There is a naturally embedded circle $\mathcal L(\alpha)$ contained in $\Sigma_g^\bullet\x[0,1]$  
whose projection by the projection $\Sigma_g^\bullet\x[0,1]\to\Sigma_g^\bullet\x\{0\}$ is 
$\mathcal I(\alpha)$. 
Suppose that $\mathcal L(\alpha)\cap(\Sigma_g^\bullet\x\{0\})=\phi$. 
Let ${\mathcal{K}}(\alpha)$ be 
 an embedded circle in $\Sigma_g\x[0,1]$ which we obtain naturally from $\mathcal L(\alpha)$.  
$\Sigma_g$ is called a {\it representing surface}. 
 $\Sigma_g^\bullet=\Sigma_g-\{*\}$ 
is also sometimes called a {\it representing surface}.  
(The closure of ) any neighborhood of the immersed circle
in $\Sigma_g$ is also 
called a {\it representing surface}.\\

\begin{thm}\label{vk}
{\rm (\cite{Kauffman1, Kauffman, Kauffmani}.)} 
Let $\alpha$ and $\beta$ be virtual 1-knot diagrams. 
$\alpha$ and $\beta$ represent the same virtual 1-knot 
if and only if 
${\mathcal{K}}(\alpha)$ is obtained from ${\mathcal{K}}(\beta)$ by 
a sequence of the following operations.

\smallbreak\noindent
$(1)$ A surgery on the surface by a 3-dimensional 1-handle, where 

\hskip2mm$($The attached part of the handle$)\cap($the projection of the embedded circle$)=\phi$.   

\smallbreak\noindent
$(2)$ A surgery on the surface by a 3-dimensional 2-handle, where 

\hskip2mm$($The attached part of the handle$)\cap($the projection of the embedded circle$)=\phi$. 

\smallbreak\noindent
$(3)$ An orientation preserving diffeomorphism map of the surface.  




\end{thm}

Hence the following definition makes sense. 
Let $K$ be a virtual 1-knot. 
Let $\alpha$ be a virtual 1-knot diagram of $K$. 
 Define $\mathcal K(K)$ to be  $\mathcal K(\alpha)$. 


\bigbreak
\section{$\mathcal E(K)$ for a virtual 1-knot $K$}\label{E}
\noindent 
We generalize spun knots and spinning tori,  
and introduce a new class of submanifolds. \\

Let $n$ be a positive integer. 
Two submanifolds $J$ and  $K$ $\subset S^n$ are {\it $($ambient$)$ isotopic}  
if there is a smooth orientation preserving family of diffeomorphisms $\eta_t$ of $S^n$, $0\leqq t\leqq1$, with $\eta_0$ the identity and $\eta_1(J)=K$.   \\

\begin{figure}
\bigbreak
     \includegraphics[width=120mm]{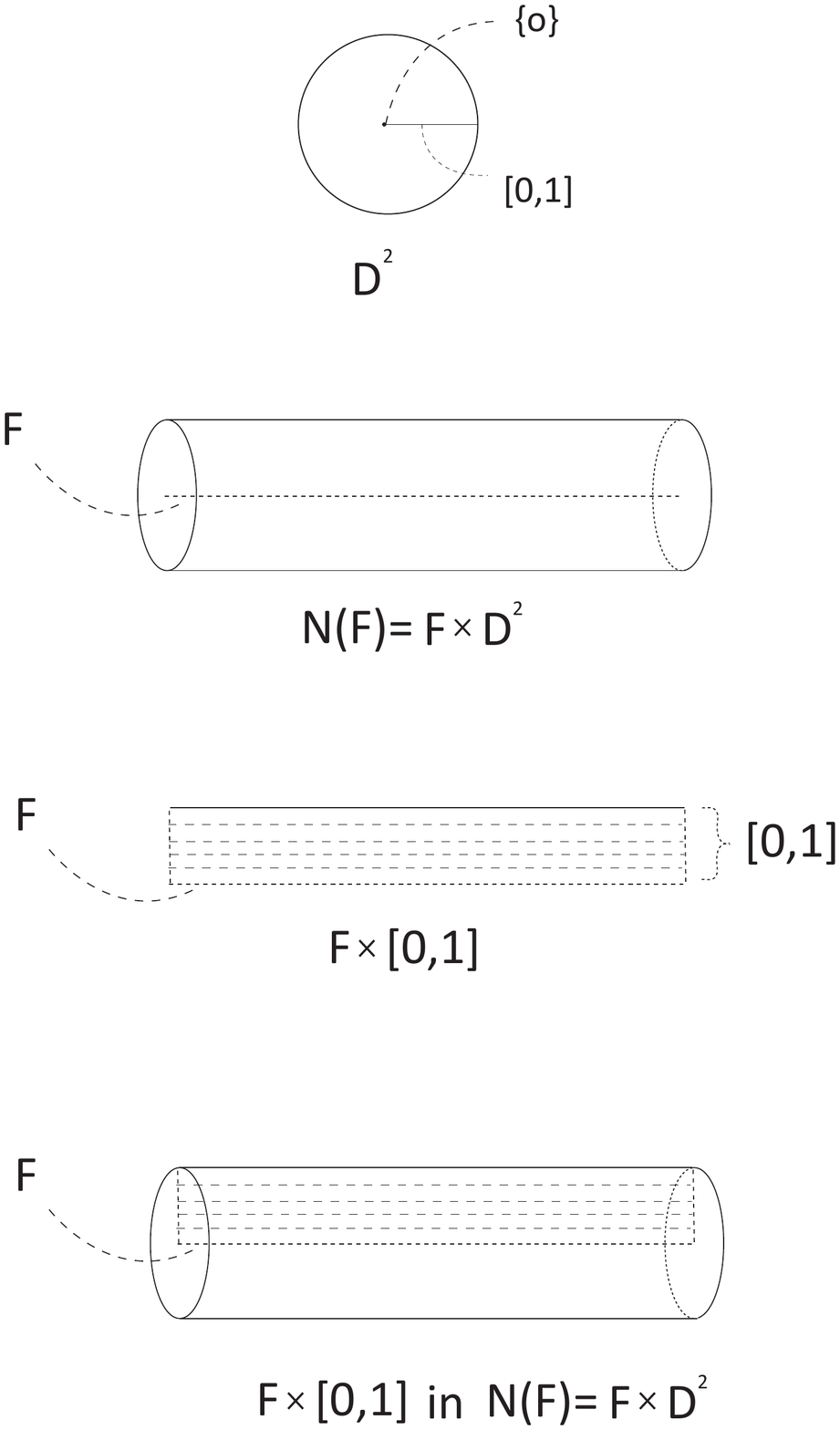}  
\vskip-10mm
\caption{{\bf The tubular neighborhood which is a product $D^2$-bundle}\label{tube}}   
\bigbreak  
\end{figure}

\begin{defn}\label{spinningsubmanifold}
Let $F$ be a codimension two submanifold contained in a manifold $X$. 
Suppose that the tubular neighborhood $N(F)$ of $F$ in $X$ 
is the product bundle. That is, we can regard $N(F)$ as $F\x D^2$. 
See Figure \ref{tube}. 
We can regard the closed 2-disc $D^2$ as the result of 
rotating a radius $[0,1]$ around the center $\{o\}$ as the axis.  
We can regard $N(F)$ as the result of 
rotating $F\x[0,1]$ around $F=F\x\{0\}$ 
as the axis.  
Suppose that a submanifold  $P$ contained in $X$ is embedded in $F\x[0,1]$. 
Let $P'$ be a submanifold $P\cap(F\x\{0\})$ of $F\x\{0\}$.  
When we rotate $F\x[0,1]$ around $F$ and make $F\x D^2$, 
rotate $P$ together, and call the resultant submanifold $Q$. 
This submanifold $Q$ contained in $X$ is called 
the {\it spinning submanifold} made from 
$P$ by the rotation in  $F\x D^2$ under the condition that 
$P\cap(F\x\{0\})$ is the submanifold $P'$.  
%
%
This way of construction of $Q$ is called a {\it spinning consruction} of submanifolds. 
If $P$ is a subset not a submanifold, we can define $Q$ as well. 
\end{defn}

Spun knots and spinning tori are spinning submanifolds. 
\cite[Proof of Claim in page 3114]{Ogasa98SL} and 
\cite[Lemma 5.3]{Ogasa98n} used spinning construction. 
By the uniqueness of the tubular neighborhood, we have the following. 

\begin{cla}\label{atarimae1}
Let $\check f$ $($respectively, $\check g)$ 
$:F\x D^2\hookrightarrow X$ be an embedding map. 
We can regard $\check f(\Sigma_g\x D^2)$ 
$($respectively, $\check g(\Sigma_g\x D^2))$     
as the tubular neighborhood of 
$\check f(\Sigma_g\x\{o\})$ $($respectively, $\check g(\Sigma_g\x\{o\}))$. 
Let $\check f|_{\Sigma_g\x\{o\}}$ 
be isotopic to  $\check g|_{\Sigma_g\x\{o\}}$. 
Then submanifolds,  
$\check f(\Sigma_g\x\{o\})$ and $\check g(\Sigma_g\x\{o\})$, 
of $X$  
are isotopic. 
\end{cla}\bigbreak

Let $\alpha$ be a virtual 1-knot diagram.    
Take $\Sigma_g\x[0,1]$ and  $\mathcal K(\alpha)$ as in \S\ref{K}, 
that is, $K(\alpha)$ is a 1-knot in $\Sigma_g\x[0,1]$, where $\Sigma_g$ representing $\alpha$. 
Assume $\mathcal K(\alpha)\cap(\Sigma_g\x\{0\})=\phi$.
Suppose $\mathcal K(\alpha)\cap(\Sigma_g\x\{0\})=\phi$.  
Make $\Sigma_g\x D^2$, where we regard $[0,1]$ as a radius of $D^2$.  
Let $\check f:\Sigma_g\x D^2\hookrightarrow S^4$ be an embedding map. 
Let $\mathcal E_{\check f}(\alpha)$  be 
the spinning submanifold made from $\mathcal K(\alpha)$ 
by the rotation in $\check f(\Sigma_g\x D^2)$.  
Note $\mathcal E_{\check f}(\alpha)\subset S^4$.
Let $f$ be $\check f|_{\Sigma_g\x\{o\}}$. 
By Claim \ref{atarimae1} it makes sense that 
 we call  $\mathcal E_{\check f}(\alpha)$,  
$\mathcal E_f(\alpha)$. \\

Suppose that $\alpha$ represents a virtual 1-knot $K$. 
Theorem \ref{honto} is one of our main results. \\

\begin{thm}\label{honto}
For an arbitrary virtual 1-knot $K$, 
the submanifold type $\mathcal E_f(\alpha)$ 
of $S^4$ does not depend on the choice of a set $(\alpha, f)$. 
\end{thm} \bigbreak

By Theorem \ref{honto}  we can define $\mathcal E(K)$ for any virtual 1-knot $K$. \\

Let $\mathcal S(\alpha)$ be an embedded $S^1\x S^1$ contained in $S^4$ 
for a virtual 1-knot diagram $\alpha$, defined by Satoh in \cite{Satoh}.    
It was proved  
there that if $\alpha$ and $\beta$ represent the same virtual 1-knot, 
the submanifolds, 
 $\mathcal S(\alpha)$ and $\mathcal S(\beta)$, 
of $S^4$ are isotopic. 
So we can define $\mathcal S(K)$ for any virtual 1-knot $K$. \\

We will prove the following in \S\ref{Proof}. 

Theorem \ref{mainkore} is one of our main results. 

\begin{thm}\label{mainkore}
Let $K$ be a virtual 1-knot. Then 
the submanifolds, 
$\mathcal E(K)$ and $\mathcal S(K)$, 
of  $S^4$ are isotopic. 
\end{thm}

\begin{note}\label{yuenchi}  
If $K$ in Theorem \ref{mainkore}  is a classical knot, 
$\mathcal E(K)$ and $\mathcal S(K)$ is the spinning torus of $K$. 
\end{note}  

\bigbreak\noindent{\bf Note.}
\cite[section 10.2]{B1}, 
\cite[section 3.1.1]{B2} and 
\cite{Dylan} proved only a special case of Theorem \ref{mainkore}, 
which is only Theorem \ref{mainmaenotame} of this paper. 
%
We prove the general case. 
Our result is stronger than the result in 
\cite[section 10.2]{B1}, 
\cite[section 3.1.1]{B2} and 
\cite{Dylan}.

\bigbreak
\section{Proof of Theorems \ref{honto}    and \ref{mainkore}}\label{Proof}  
\noindent 
We first prove a special case.

\begin{thm}\label{mainmaenotame}
Take a virtual 1-knot diagram $\alpha$ in \S\ref{K}.  
Let $\check\iota:\Sigma_g\x D^2\to S^4$ be 
an embedding map whose image of  
$\Sigma_g^\bullet$ 
by $\check\iota$  
 is $\Sigma_g^\bullet$ 
 in \S\ref{K}.  
Let $\iota$ be $\check\iota\vert_{\Sigma_g}$. 
Then the submanifolds, 
$\mathcal E_{\iota}(\alpha)$ and $\mathcal S(\alpha)$, 
of $S^4$ are isotopic. 
\end{thm}

\noindent{\bf Proof of Theorem \ref{mainmaenotame}.}
Let $\R^4=\{(x,y,u,v)|x,y,u,v\in\R\}$,  
$\R^2_b=\{(x,y)|x,y\in\R\}$, 
and $\R^2_F=\{(u,v)|u,v\in\R\}$. 
Note 
$\R^4=\R^2_b\x\R^2_F$. 
Regard $\R^3$ in \S2 as \\
$\R^2_b\x\{(u,v)| u\in\R, v=0\}$. 
Take the tubular neighborhood of $\Sigma_g^\bullet$ in $\R^3$.  
It is diffeomorphic to $\Sigma_g^\bullet\x[-1,1]$. 
We can suppose that 
 $\Sigma_g^\bullet, \Sigma_g^\bullet\x[0,1]
\subset\R^2_b\x\{(u,v)| u\geqq0, v=0\}$.  \\

\noindent{\bf Note.} 
  Let $\Sigma_g\subset S^4$. 
Suppose that $\{*\}\in S^4$ is included in $\Sigma_g$. 
Then \newline
$S^4-\Sigma_g=(S^4-\{*\})-(\Sigma_g-\{*\})=R^4-\Sigma_g^\bullet$. \\

Take the tubular neighborhood $N(\Sigma_g^\bullet)$ of $\Sigma_g^\bullet$ in $\R^4$.  
Note that 
$N(\Sigma_g^\bullet)$ is diffeomorphic to $\Sigma_g^\bullet\x D^2$. 
We can regard $N(\Sigma_g^\bullet)$ 
as the result of 
rotating $\Sigma_g^\bullet\x[0,1]$ 
around $\Sigma_g^\bullet$ as the axis 
(diffeomorphically not isometrically). 
Suppose that $\mathcal L(\alpha)\cap(\Sigma_g^\bullet\x\{0\})=\phi$.  
Make the spinning submanifold 
 $\mathcal E_\iota(\alpha)$ 
from $\mathcal L(\alpha)$. 
\\

We can suppose that 
each fiber $D^2$ of 
$N(\Sigma_g^\bullet)$ 
is parallel to $\{(x,y)|x=0, y=0\}\x\R^2_F$ 
by using an isotopy of  an embedding map of the tubular neighborhood.  \\

\begin{figure}
\includegraphics[width=130mm]{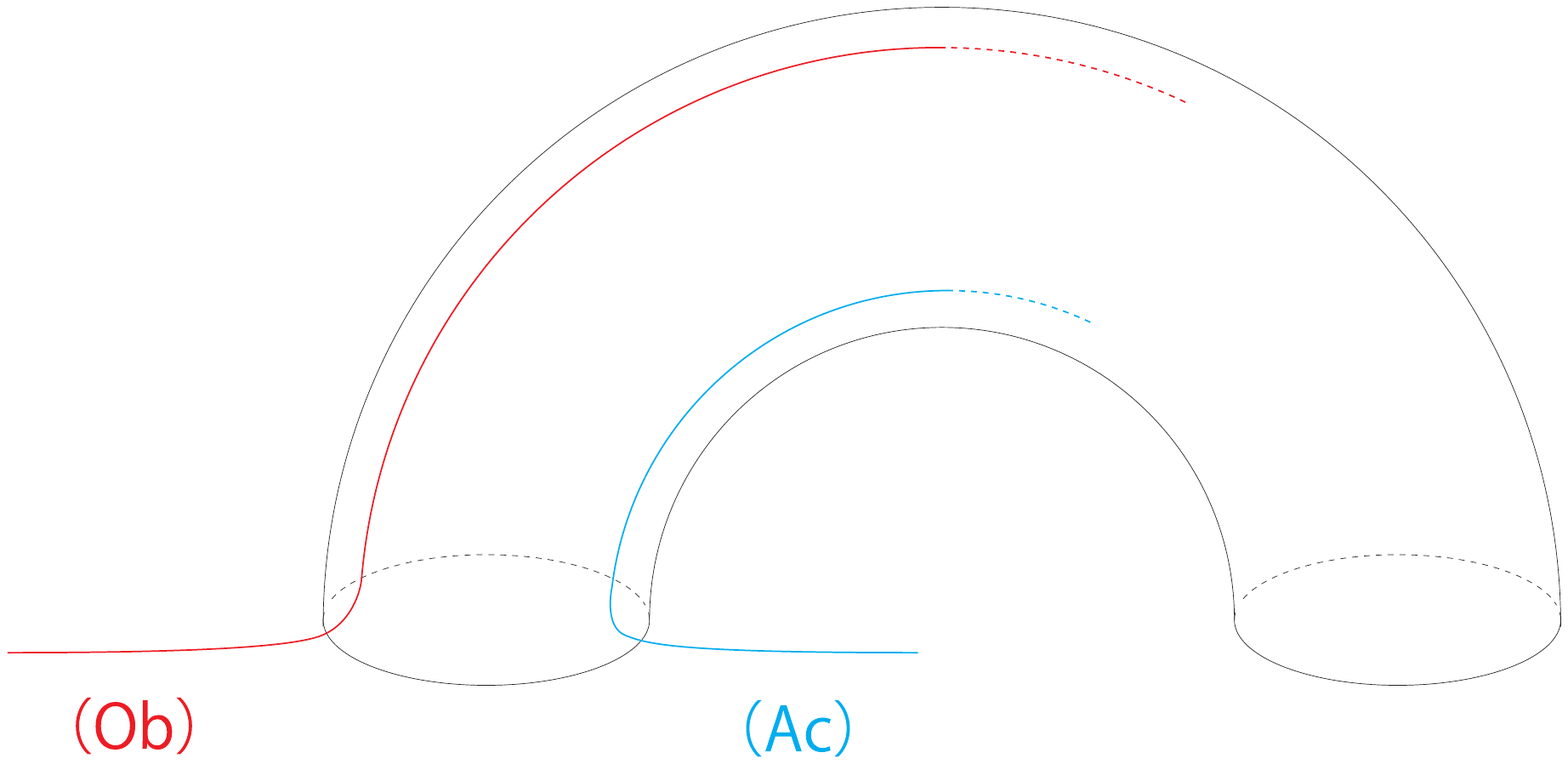}  
\vskip-50mm
\caption{{\bf (Ac) and (Ob). 
}\label{Alaska}}   
\end{figure}

We can suppose that $\mathcal I(\alpha)$ intersects each fiber $D^2$ transversely.  
{\it Reason.} 
Note a 1-handle  drawn in the right-side of Figure \ref{Alabama}.  
If $\mathcal I(\alpha)$ near the 1-handle is put like (Ac) in Figure \ref{Alaska}, 
$\mathcal I(\alpha)$ does not intersect each fiber $D^2$ transversely.  
However we can do the following operation. 
By using an isotopy of a part of $\mathcal I$ 
we change 
the part of  $\mathcal I(\alpha)$  
from (Ac) to  (Ob)  in Figure \ref{Alaska}. 
After this operation, $\mathcal I(\alpha)$ intersects each fiber $D^2$ transversely.  \\

\begin{note}\label{kabuto}
We will explain a property of (Ac), in Note \ref{kuwagata}. 
It is important. We will use it in Alternative proof of Claim \ref{shichi} of \S\ref{v2}.  
\end{note} \bigbreak

Note that, even if a part of  $\mathcal I(\alpha)$ is (Ac), 
we can make a spinning submanifold $\mathcal E_\iota(\alpha)$. 
However, if there is not (Ac), we have an advantage as below. \\

\begin{figure}
\includegraphics[width=150mm]{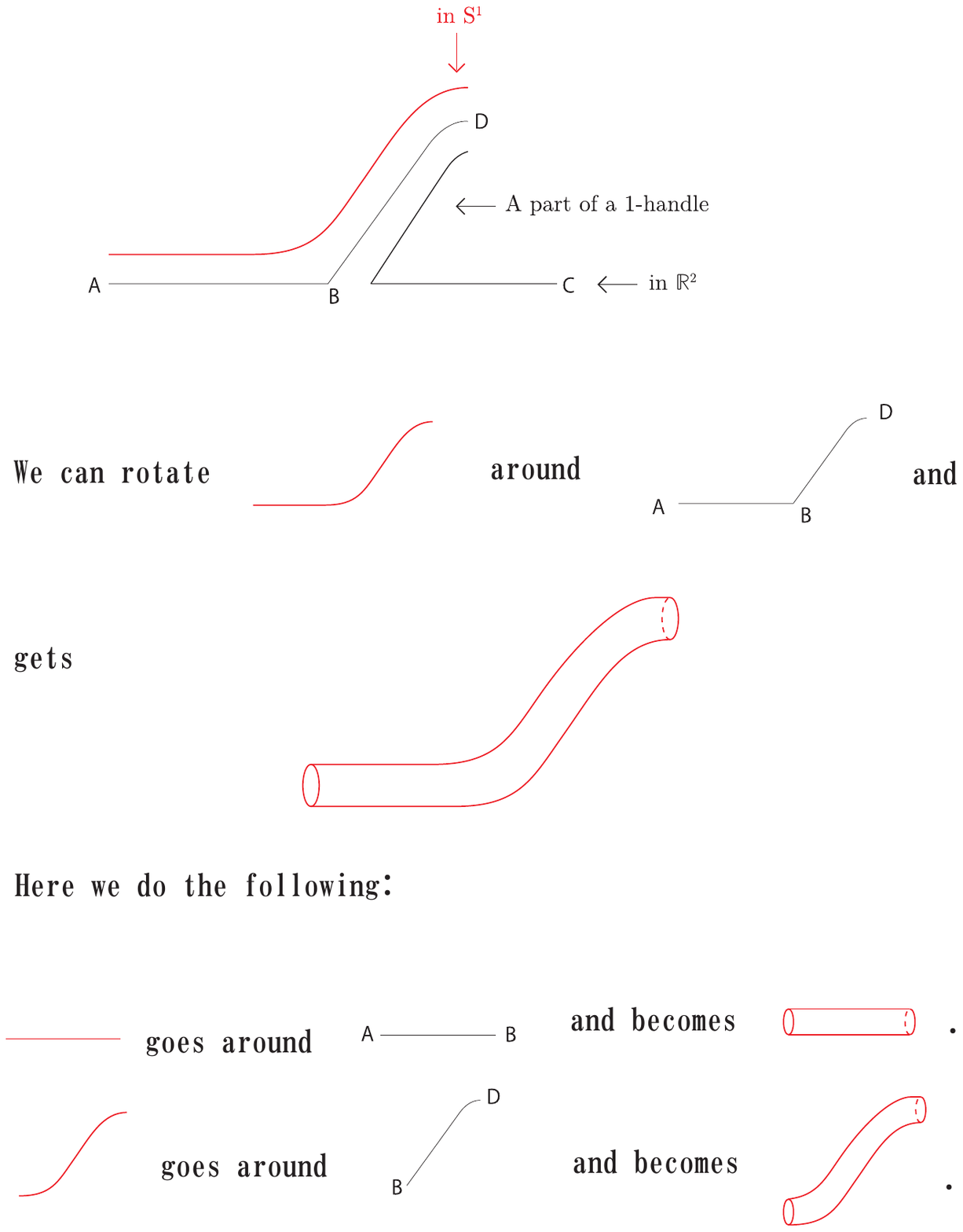}  
\vskip-40mm
\caption{{\bf Rotation around a part near (Ob). 
The reason why (Ob) is useful for us. 
}\label{Arizona}}   
\end{figure}

\begin{figure}
\vskip-7mm
     \includegraphics[width=150mm]{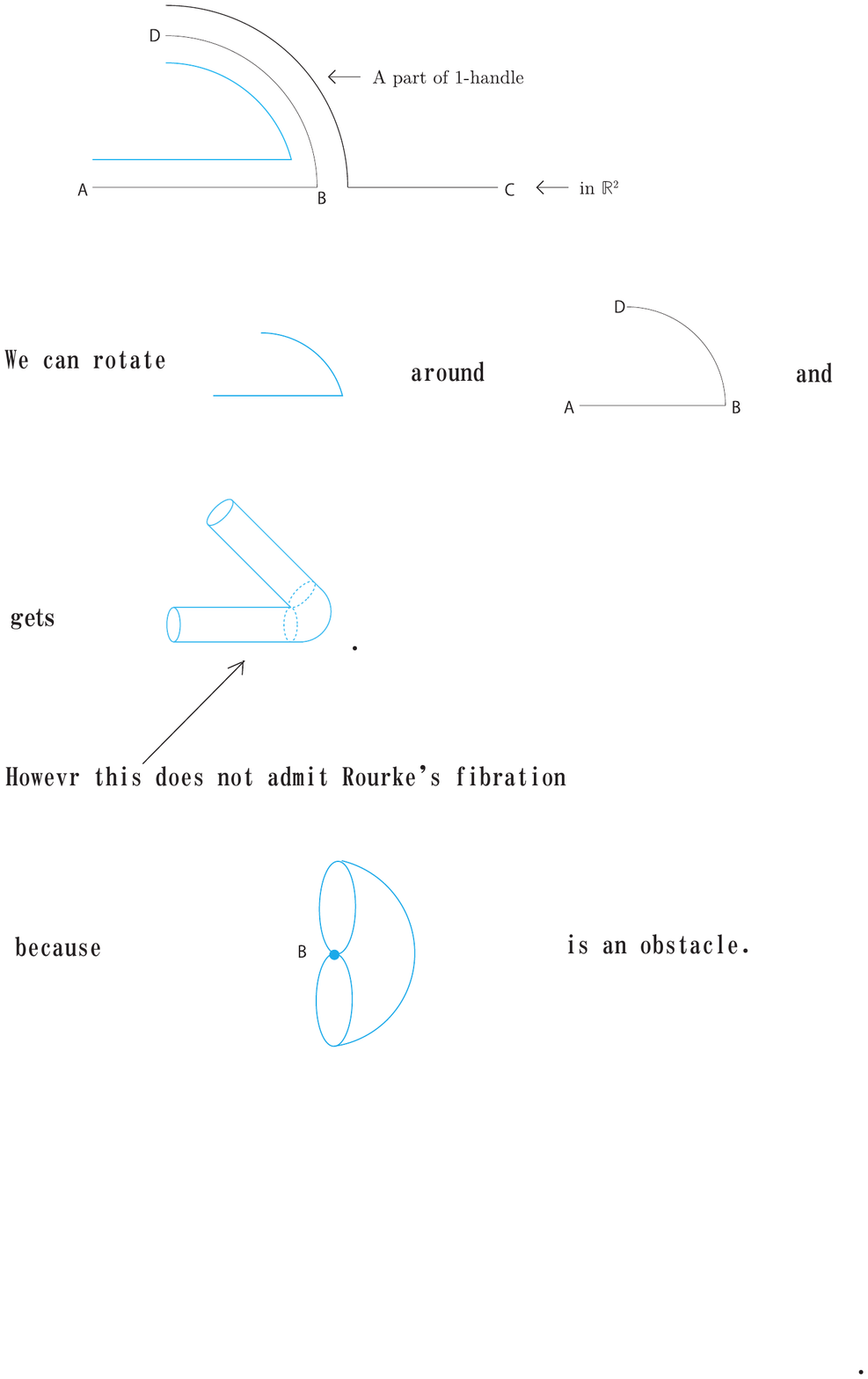}  
\vskip-47mm
\caption{{\bf Rotation around a part near (Ac). 
The reason why (Ac) is not useful for us. 
This property of (Ac) is used in in Alternative proof of Claim \ref{shichi} of \S\ref{v2}. 
}\label{Arkansas}}   
\bigbreak  
\end{figure}

If there is not (Ac) in $\mathcal I(\alpha)$, we have the following. \\
 
Take a point $q\in\R^2_b\x\{(u,v)| u=0, v=0\}$. 
Note $\alpha\subset\R^2_b\x\{(u,v)| u=0, v=0\}$. 
By the above construction of $\mathcal E_\iota(\alpha)$, we have the following. \\

\smallbreak\noindent (i)  
If $q\cap \alpha=\phi$, 
 $(\{q\}\x\R^2_F)\cap\mathcal E_\iota(\alpha)=\phi.$ 

\smallbreak\noindent (ii)
If $q$ is a normal point of $\alpha$, 
 $(\{q\}\x\R^2_F)\cap\mathcal E_\iota(\alpha)$ is a single circle in $\{q\}\x\R^2$.

\smallbreak\noindent (iii)
If $q$ is a real crossing point of $\alpha$, 
 $(\{q\}\x\R^2_F)\cap\mathcal E_\iota(\alpha)$ 
is two circles in $\{q\}\x\R^2$ 
such that one of the two is in the inside of the other.  
The inner (respectively, outer) circle corresponds 
to the lower (respectively, upper) point of the singular point.

\smallbreak\noindent (iv)
If $q$ is a virtual crossing point of $\alpha$, 
 $(\{q\}\x\R^2_F)\cap\mathcal E_\iota(\alpha)$ 
is two circles in $\{q\}\x\R^2$ 
such that each of the two is in the outside of the other each other.   


It is Rourke's description of $\mathcal S(\alpha)$ 
in  Theorem \ref{Montana} which is cited below from \cite{Rourke}.  
(However \cite{Rourke} does not write a proof.  
So \cite{J} wrote a proof.)\\


\begin{thm}\label{Montana}
{\bf (\cite{Rourke}.)} 
Let $\alpha$ be a virtual 1-knot diagram. 
Take an embedding map 
$\varphi: S^1_b\x S^1_f\hookrightarrow \R^2_b\x\R^2_f$ with the following properties. 

\smallbreak\noindent$(1)$  
Let $\pi:\R^2_b\x\R^2_f\to\R^2_b$ be the natural projection. 
$\pi\circ\varphi(S^1_b\x S^1_f)\subset\R^2_b$ defines 
$\alpha$ without the notations of virtual crossings. 

\smallbreak\noindent$(2)$  
For points in $\R^2_b$, we have the following:

\begin{enumerate} 
\item
If $q\notin\alpha$, 
we have $\pi^{-1}(q)=\phi.$ 

\smallbreak
\item
If $q$ is a normal point of $\alpha$, 
we have that $\pi^{-1}(q)$ is a circle in $\{q\}\x\R^2$.

\smallbreak
\item
If $q$ is a real crossing point of $\alpha$, 
we have that $\pi^{-1}(q)$ is  two circles in $\{q\}\x\R^2$ 
such that one of the two is in the inside of the other. 
The inner $($respectively, outer$)$ circle corresponds 
to the lower $($respectively, upper$)$ point of the singular point.

\smallbreak
\item
If $q$ is a virtual crossing point of $\alpha$, 
we have that $\pi^{-1}(q)$ is  two circles in $\{q\}\x\R^2$ 
such that each of the two is in the outside of the other each other.   
\end{enumerate}

\smallbreak
Then the submanifolds, 
$\mathcal S(\alpha)$ and $\varphi(S^1_b\x S^1_f)$, 
of $S^4$ are isotopic 

\end{thm}

This completes the proof of Theorem \ref{mainmaenotame}. 
\qed\\

\begin{defn}\label{Nebraska}  
In Theorem \ref{Montana}, 
each circle
$f(S^1_b\x S^1_f)\cap$(each fiber $\R^2_f$) 
is called a {\it fiber-circle}. 
We say that $f(S^1_b\x S^1_f)$ admits {\it Rourke's fibration}.
\end{defn}\bigbreak

\begin{note}\label{kuwagata}  
As we preannounced in Note \ref{kabuto}, we state a comment on (Ac).  
If the projection on a surface includes (Ac), 
$\mathcal E(\alpha)$ does not admit Rourke's fibration.  
The reason is  explained in Figures \ref{Arizona} and \ref{Arkansas}.
We will use this property, 
which is raised by the difference between (Ac) and (Ob),     
in Alternative proof of Claim \ref{shichi} of \S\ref{v2}. 
\end{note}\bigbreak


We next prove the general case. \\

\noindent{\bf Proof of Theorems \ref{honto} and \ref{mainkore}.}  
We prove Theorem \ref{oh} below. 
The key idea of the proof is Claim \ref{wow}. 
Let $\Sigma$ be a closed oriented surface.  
Let $G_1$ and $G_2$ be submanifolds of $S^4$ 
which are orientation preserving diffeomorphic to $\Sigma$. 
It is known that 
there is a case that 
the submanifolds, $G_1$ and $G_2$, of  $S^4$ are non-isotopic. 
 Let   $\Sigma^\circ$ denote $\Sigma-(\text{an open 2-disc})$. 
Let $$G^\circ_i\\=G_i-(\text{an open 2-disc})$$ be a submanifold of $S^4$ 
which are orientation preserving diffeomorphic to $\Sigma^\circ$ $(i=1,2)$.  \\

\begin{cla}\label{wow} 
The submanifolds, $G^\circ_1$ and $G^\circ_2$, of $S^4$ are isotopic. 
\end{cla}

\noindent{\bf Proof of Claim \ref{wow}.} 
$\Sigma^\circ$ has a handle decomposition 
which consists of one 0-handle, 1-handles, and no 2-handle.  
\qed\\

Let $i\in\{1,2\}$. 
We can regard the tubular neighborhood of $G_i$ in $S^4$ as  $G_i\x D^2$. 
Embed $S^1$ in  $G_i\x [0,1]$, where we regard $[0,1]$ as a radius of $D^2$,  
and call the image $J_i$. Assume that  $J_i\cap(G_i\x\{0\})=\phi$. 
 Suppose that there is 
a bundle map $\check\sigma:G_1\x D^2\to G_2\x D^2$ 
such that $\check\sigma$ covers 
an orientation preserving diffeomorphism map $\sigma:G_1\to G_2$ and 
such that $\check\sigma(J_1)=J_2$. \\

Define a submanifold $E_i$ contained in $S^4$ 
to be the spinning submanifold made from $J_i$ 
by the rotation  in $G_i\x D^2$. 

\begin{thm}\label{oh}
The submanifolds, $E_1$ and $E_2$, of $S^4$ are isotopic. 
\end{thm}

\noindent{\bf Proof of Theorem \ref{oh}.} 
We can suppose that 
$J_i\subset G^\circ_i\x [0,1]$. 
By the existence of $\sigma$,  
there is a bundle map $\check\tau:G^\circ_1\x D^2\to G^\circ_2\x D^2$ 
such that $\check\tau$ covers an orientation preserving diffeomorphism map 
$\tau: G^\circ_1\to  G^\circ_2$ and 
such that $\check\tau(J_1)=J_2$. \\

Note the following: 
Let $f:\Sigma^\circ\to S^4$ be an embedding map. 
We can regard $\tau$ as a diffeomorphism map 
$\Sigma^\circ\to \Sigma^\circ$. 
By Claim \ref{wow} , 
the submanifolds, 
$f(\Sigma^\circ)$ and $f(\tau(\Sigma^\circ))$,  
of $S^4$ are isotopic. 
Therefore the submanifolds, $E_1$ and $E_2$, of $S^4$ are isotopic, 
\qed\\

Theorems \ref{vk} and \ref{oh}  imply Theorems \ref{honto} and \ref{mainkore}      \qed\\

 We can extend all discussions in \S\S\ref{K}-\ref{Proof} and the following \S\ref{rt} to the virtual 1-link case easily.
 When we define $\mathcal E(\alpha)$ in \S\ref{E}, 
we assume 
$\mathcal L(\alpha)\cap\Sigma_g^\bullet=\phi$. 
Suppose that 
$\mathcal L(\alpha)\cap(\Sigma_g^\bullet)$ is an arc instead. 
Then we obtain a spherical 2-knot in $\R^4$ 
as the spinning submanifold. 
The class of such spherical 2-knots is also a generalization of 2-dimensional spun-knots of 1-knots, 
and is also worth studying. 
As we state in \S\ref{jobun}, we do not discuss this class in this paper.

\bigbreak
\section{Immersed solid tori} \label{rt}
\begin{figure}
     \includegraphics[width=150mm]{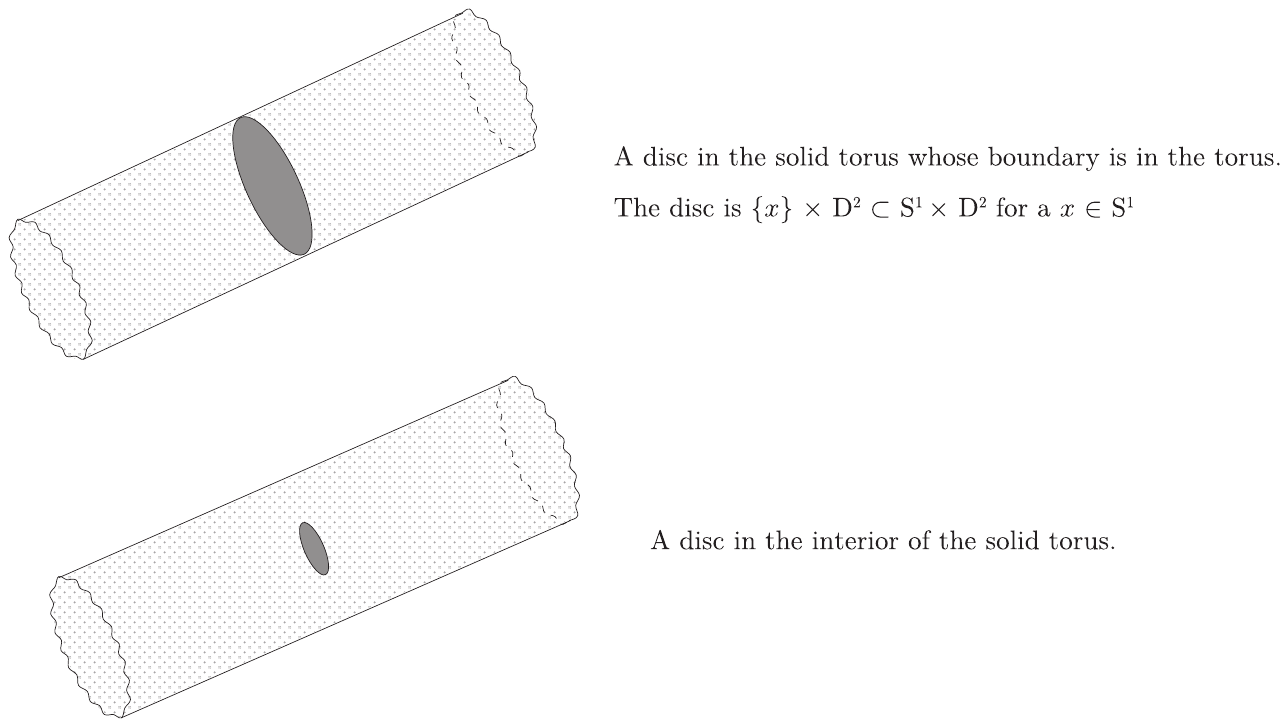}  
\vskip-40mm
\caption{{\bf  $\zeta^{-1}($each closed 2-disc$)$}\label{Colorado}}   
\end{figure}

\noindent By the definition of $\mathcal S(\hskip2mm)$ in \cite{Satoh}, we have (i)$\Rightarrow$(ii). 

\smallbreak\noindent
(i) An embedded torus $Y$ contained in $S^4$ is isotopic 
to  $\mathcal S(\alpha)$ for a virtual 1-knot diagram $\alpha$.

\smallbreak\noindent
(ii) There is an immersion map $\zeta:S^1\x D^2\looparrowright S^4$ 
with the following properties: \newline
$\zeta(S^1\x\partial D^2)$ is $Y$. 
The singular point set of $\zeta$ consists of double points and is a disjoint union of closed 2-discs,  
and $\zeta^{-1}($each closed 2-disc$)$ is as shown in Figure \ref{Colorado}. 
 
\smallbreak   
By using the construction of $\mathcal E(\alpha)$,  
we can also describe the immersed solid torus \\ in (ii) 
as follows: 
By using the projection  `$\mathcal L(\alpha)\to \mathcal I(\alpha)$'  in \S\ref{K}, 
we can make an immersed annulus in 
$\Sigma_g\x[0,1]$ naturally.
Note that (the immersed annulus)$\cap\Sigma_g\x\{0\}\neq\phi$.  
Make a subset from this immersed annulus 
by a spinning construction around $\Sigma_g$, defined in Definition \ref{spinningsubmanifold}. 
Then the result is an immersed solid torus in (ii). 
\smallbreak

We prove the converse of the above claim,  
that is, the following.

\begin{thm}\label{grt} 
{\rm(ii)}$\Rightarrow${\rm(i)}. 
\end{thm}

We prove this theorem as an application of our results  in \S\ref{Proof} 
although it may be also proved in another way. 
\\

\noindent{\bf Proof of Theorem \ref{grt}.}
Let $q\in\partial D^2$. 
Let $C$ be $\zeta(S^1\x\{q\})$. 
In the following paragraphs, 
for $Y$, we will make an embedded oriented surface $F$ contained in $S^4$ 
so that we will put $C$ in the tubular neighborhood $N(F)$ of $F$ in $S^4$. 
We will make $C\cap F=\phi$. 
We will make $Y$ so that it will be the spinning submanifold of $C$ around $F$.  
Let $\{o\}$ be the center of $D^2$. 
We will let $F$ include $\zeta(S^1\x\{o\})$. \\

Let $\xi:S^1\x D^2\x I\looparrowright S^4$ be an immersion map, where $I=[-1,1]$,  
to satisfy that  
$\xi\vert_{S^1\x D^2\x \{0\}}=\zeta$ and that 
\newline\hskip3cm$\xi(\{x\}\x \{o\}\x I) \quad\bot\quad \xi(\{x\}\x D^2\x \{0\})$ \newline
for each $x$ if we give appropriate metrics to $S^4$ and $S^1\x D^2\x I$. 
Then we can suppose the following:  

\smallbreak\noindent(1) 
$P=\xi(S^1\x \{o\}\x I)$ is a boundary-connected-sum of $n$ copies of the annulus ($n\in\N$). 

\begin{figure}
\begin{center}
\includegraphics[width=70mm]{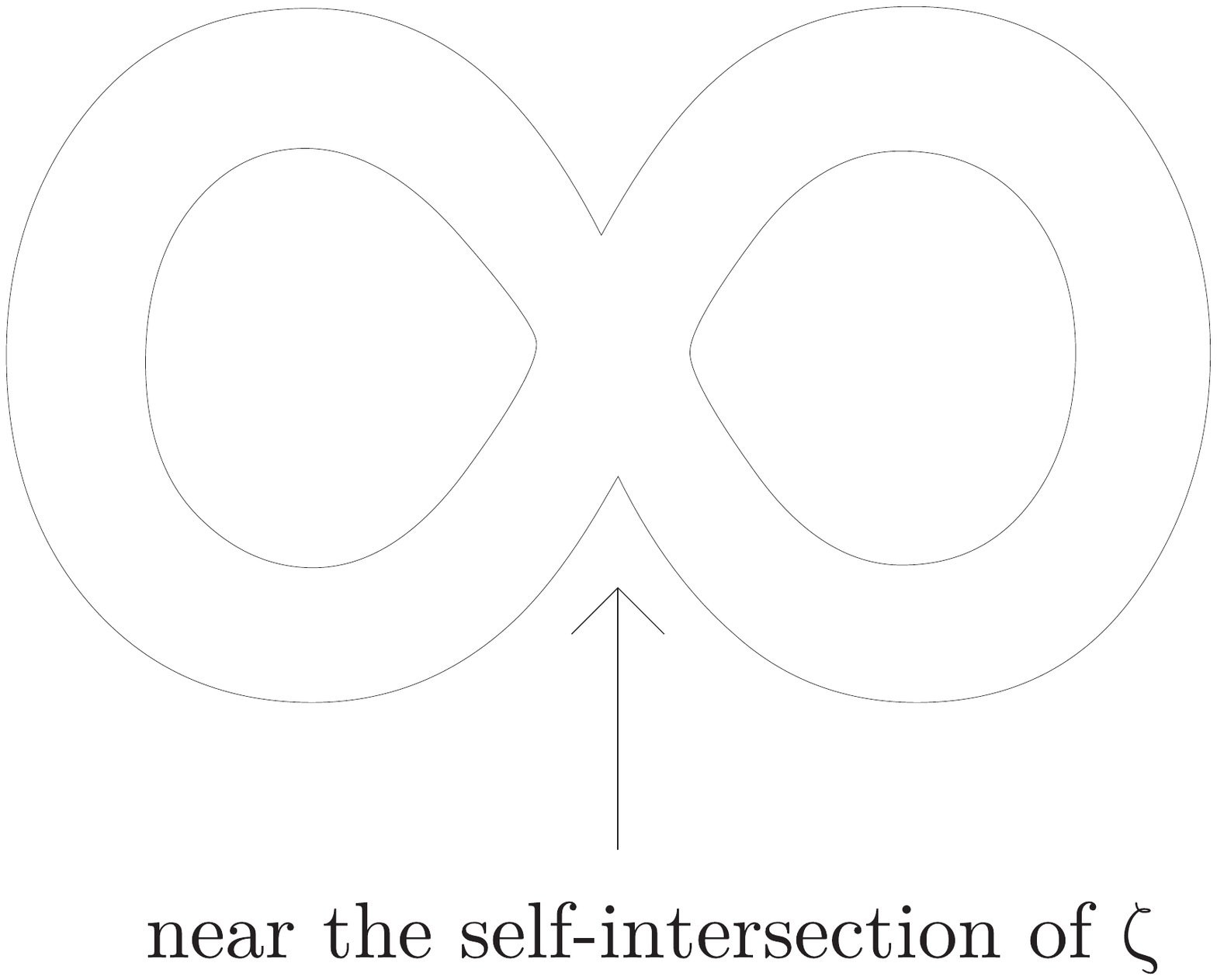}  
\end{center}
\vskip-20mm
\caption{{\bf  An example of $P$}\label{Connecticut}}   
\end{figure}

\quad See Figure \ref{Connecticut} 
for an example of $P$.  

\smallbreak\noindent(2)  
$Q=\xi(S^1\x D^2\x I)$ is a boundary-connected-sum of $m$-copies of $S^1\x B^3$
($m\in\N$). 

\smallbreak\noindent(3)  
$\partial P\subset\partial Q$. 

\smallbreak\noindent(4)  
$Q$ is the tubular neighborhood of $P$ in $S^4$. 
$Q$ is diffeomorphic to $P\x D^2$. \\

By the Mayer-Vietoris sequence,  
$H_1(S^4, Q;\Z)\cong 
H_1(S^4- {\rm Int} Q, \partial Q;\Z)\cong0$. 
Hence there is an embedded oriented compact surface-with-boundary 
$G$ contained in $S^4-{\rm Int} Q$ such that 
$\partial G=\partial P$ and that 
$G\cup P$ is an oriented closed surface $F$.   
({\it Reason}: Consider a simplicial decomposition of $S^4-{\rm Int} Q$.)
We can regard $Y$ as the spinning submanifold made from  $\zeta(S^1\x\{q\})$ around $F$. 
Hence we can regard $\zeta(S^1\x\{q\})$ as 
$\mathcal K(\beta)$ for a virtual 1-knot diagram  $\beta$ 
in a fashion which is explained in \S\S\ref{E}-\ref{Proof},  
and can regard $Y$ as $\mathcal E(\beta)$.\\

This completes the proof of Theorem \ref{grt}. \qed

\bigbreak
\section{The virtual 2-knot case}\label{v2}
\noindent

\noindent
Virtual 2-knot theory is defined analogously to Virtual 1-knot theory, using
generic surfaces in 3-space as knot diagrams and using Roseman moves for knot
equivalence, and allowing the double-point arcs to have classical or virtual crossing
data.
See \cite{Takeda, J}. 
Virtual 2-knot diagrams (respectively, virtual 2-knots) in \cite{J} and this paper 
are the same as 
virtual surface-knots (respectively, virtual surface-knot  diagrams) in \cite{Takeda}.

\begin{defn}\label{oyster} 
Let $F$ be a closed surface. 
A smooth map $f : F\to\R^3$ is considered {\it quasi-generic} 
if it fails to be one-to-one only at transverse crossings of orders 2 and 3
as shown in Figures \ref{sashimid} and \ref{sashimit}, 
\begin{figure}
\begin{center}
\includegraphics[width=40mm]{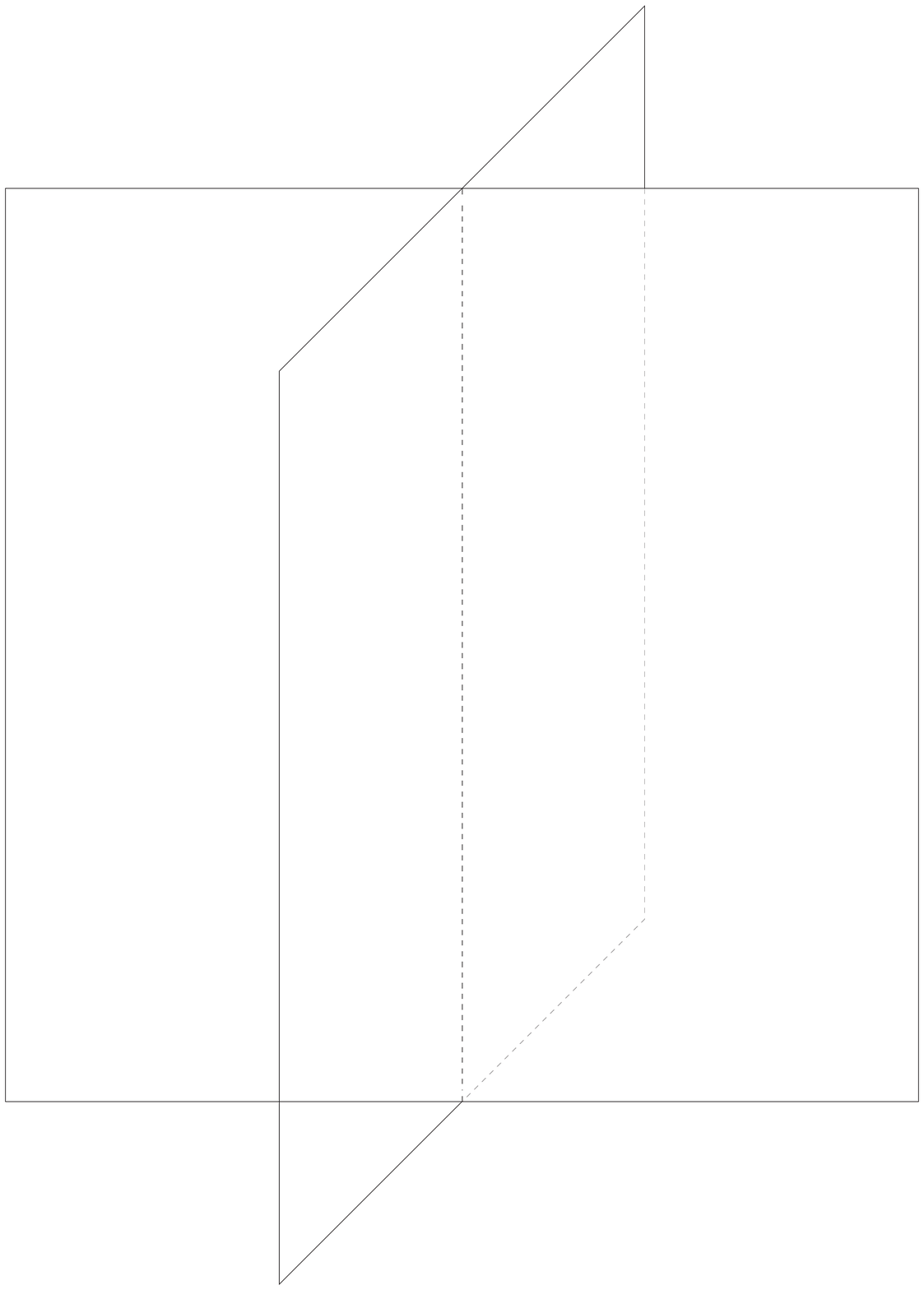}  
\end{center}
\caption{{\bf Transversal double points
}\label{sashimid}}   
\end{figure}
\begin{figure}
 \includegraphics[width=70mm]{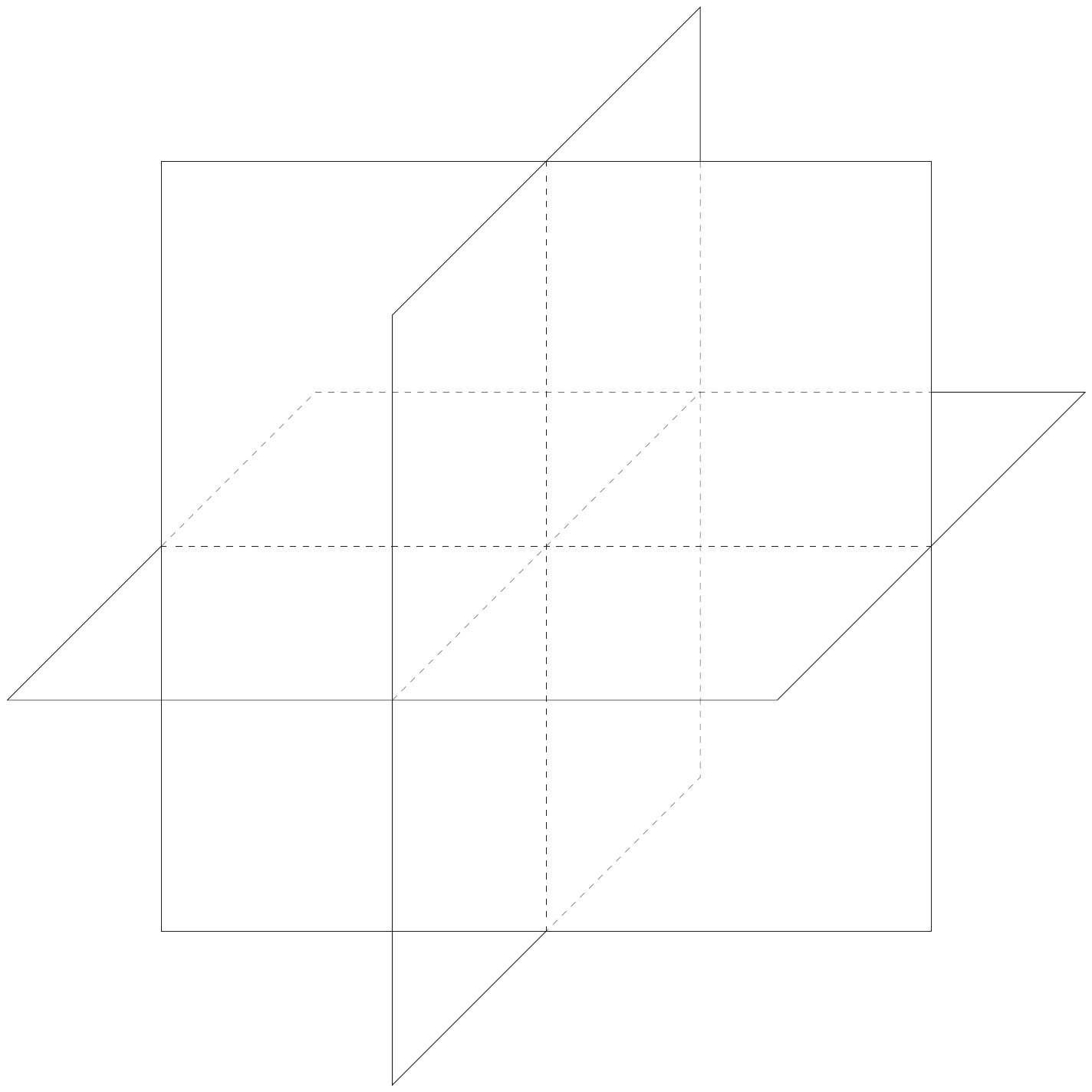}  
\caption{{\bf Transversal triple points}\label{sashimit}}   
\end{figure}
and it fails to be regular only at isolated {\it branch points}  
where, locally, the image of a disk looks like the cone over a loop, with no other parts of the surface touching the vertex. 
See  Figure \ref{gbra}. 
Branch points include the cone over any closed, regular, transversely self-intersecting
curve.  
In particular, the cone over a figure-$\infty$ curve is called 
a {\it Whitney branch point}. See Figure \ref{sashimiv}.\\

\begin{figure}
 \includegraphics[width=70mm]{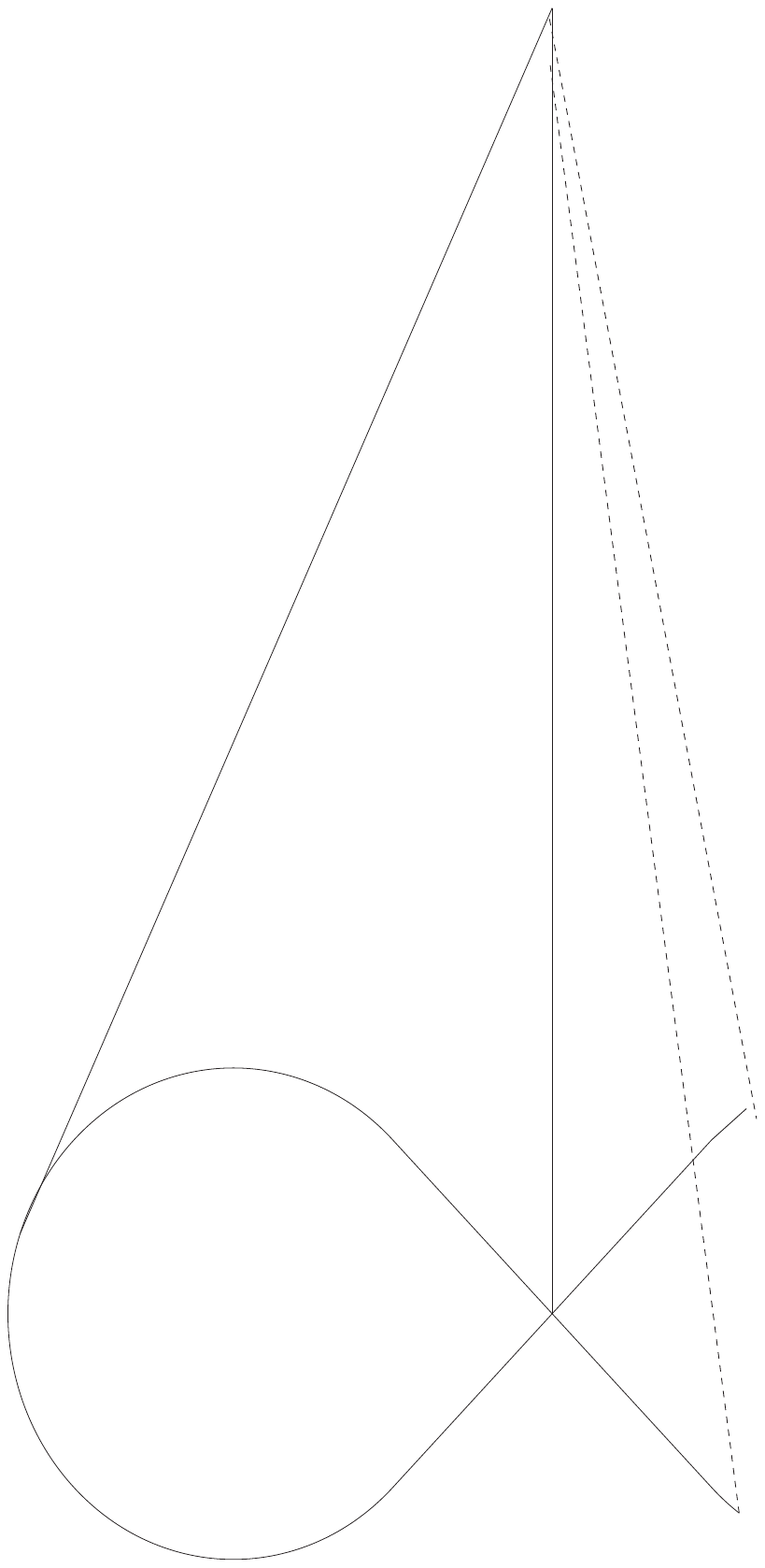}  
\caption{{\bf A general branch point}\label{gbra}}   
\end{figure}

\begin{figure}
 \includegraphics[width=70mm]{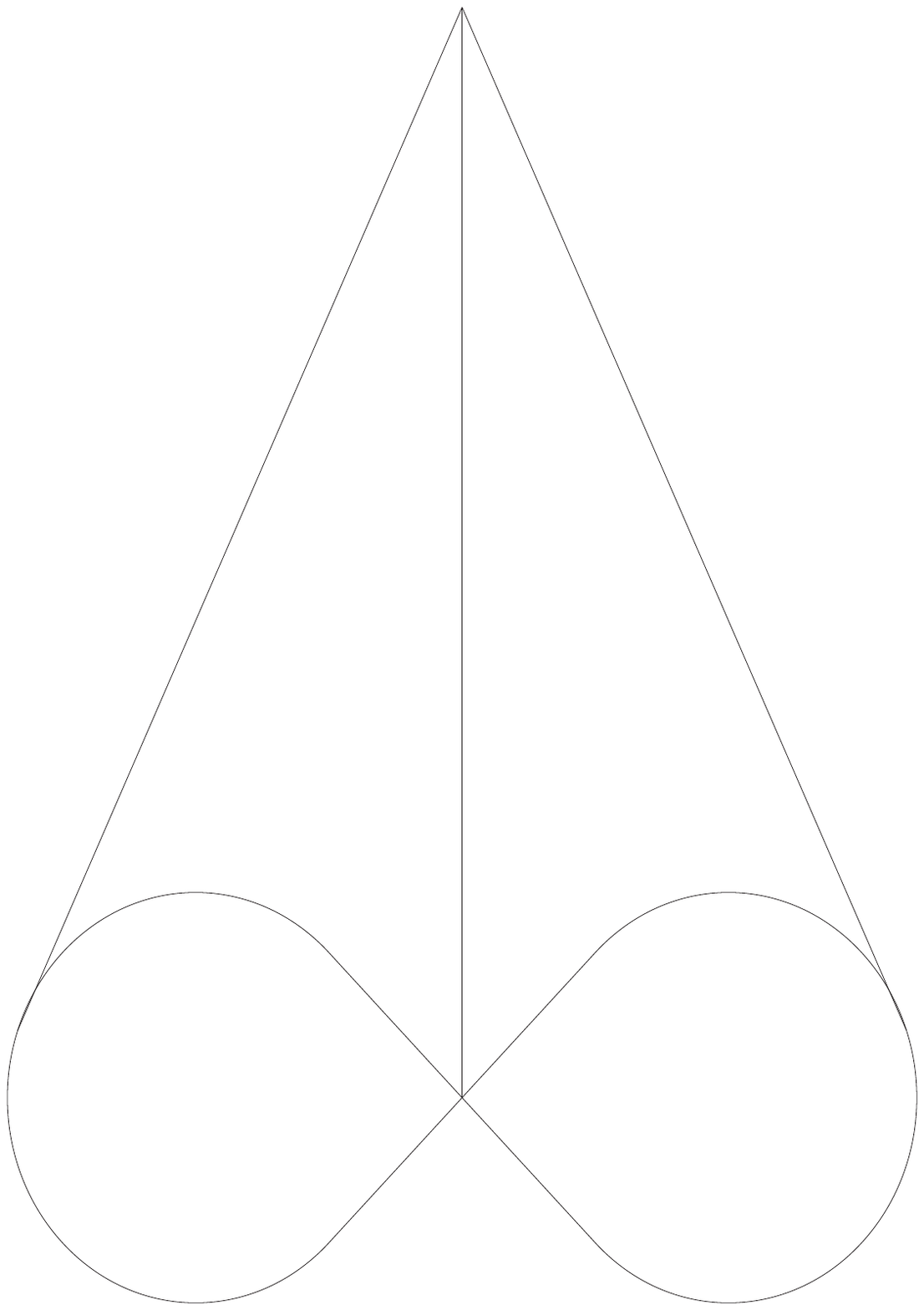}  
\caption{{\bf Whitney-umbrella branch point
}\label{sashimiv}}   
\end{figure}

A quasi-generic map $f$ is {\it generic} if the only branch points are Whitney
branch points.
The three features of a generic map--- 
Whitney branch points, double-point arcs,
and triple points--- 
have slice-histories corresponding respectively to the Reidemeister
$I$-, $II$-, and $III$- moves in 1-knot theory.
\end{defn}

\begin{figure}\bigbreak   \includegraphics[width=150mm]{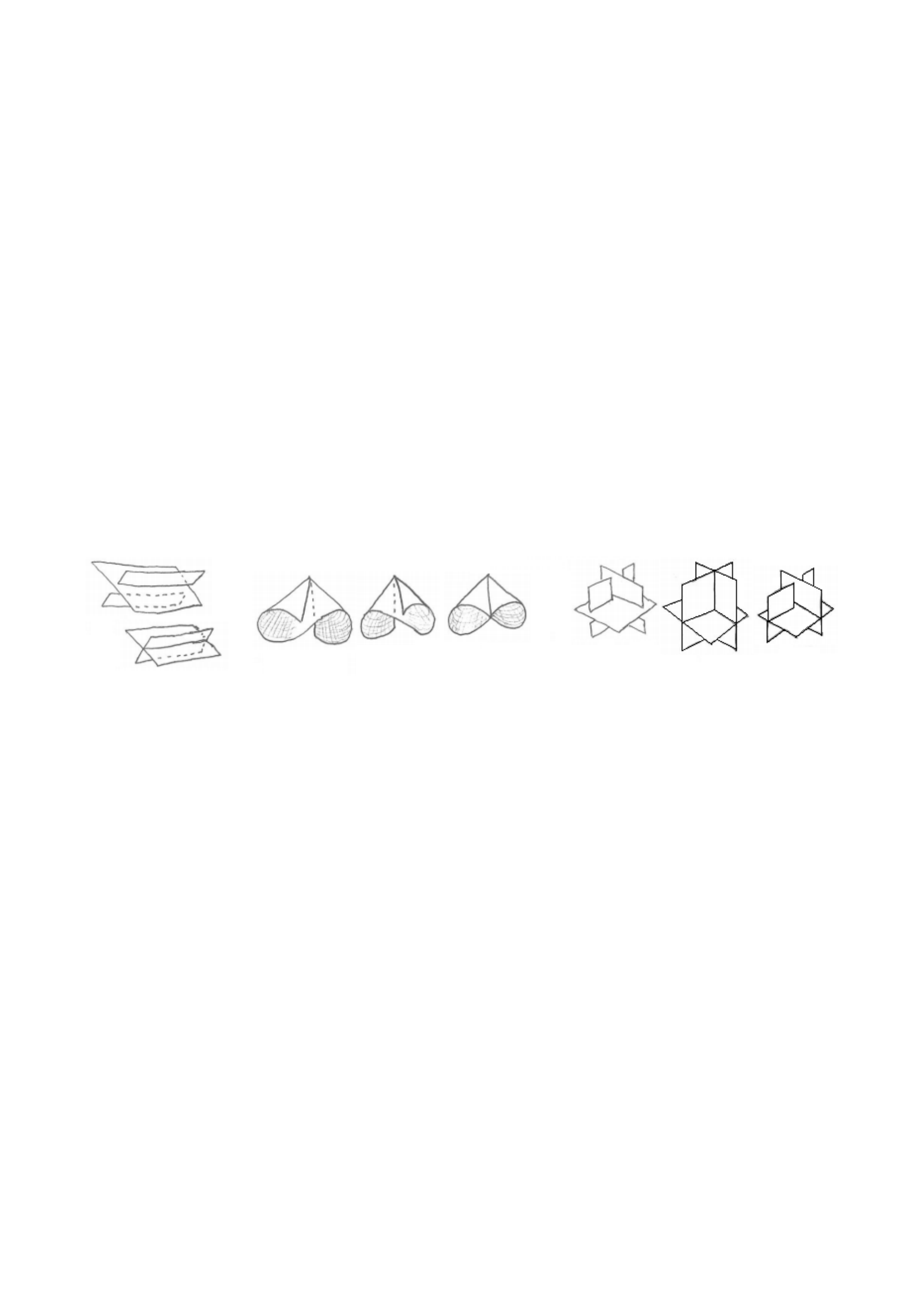}  
\vskip-40mm
\caption{{\bf  The singular point sets of virtual 2-knots}\label{JV1}}   \bigbreak  \end{figure}

\begin{defn}\label{sashimi} 
A {\it virtual 2-knot diagram} consists of a generic map $F$ 
together with classical and virtual crossing data along its double-point arcs. 
Crossing data is representedgraphically as broken and unbroken surfaces: 
See the left two figures of Figure \ref{JV1}.
Branch points can be classical or virtual:
See the middle three figures, 
the figures which are not the above ones nor the following ones, 
in Figure \ref{JV1}. 
At triple points, three crossings meet. Triple points of the following types are
allowed: See the right three figures in Figure \ref{JV1}. 
All other combinations of crossing data are forbidden. Note that the three allowed
triple points have slice-histories corresponding to the Reidemeister $III$-moves 
in Virtual 1-knot theory.
A virtual 2-knot diagram may be reduced to its bare combinatorial structure,
forgetting all but the information that is invariant under isotopies of $\R^3$ and $F$.
In this regard, we do not distinguish diagrams that are related by isotopies of 
$\R^3$ and $F$.
\end{defn}\bigbreak

\begin{defn}\label{JV}{\bf (\cite[section 3.5]{J}.)}     
A virtual 2-knot diagram may be transformed by Roseman moves. There are
seven types of local moves, shown here without crossing data.
When a virtual 2-knot diagram undergoes a Roseman move, its crossing data
carried continuously by the move. Two diagrams related by a series of Roseman
moves are called {\it virtually equivalent}.  
The equivalence classes are {\it Virtual 2-knot types}, or sometimes simply 
{\it Virtual 2-knots}.
\end{defn}

\noindent{\bf Note.} 
The readers need not be familiar with Roseman moves in order to read this paper.\\ 

It is natural to ask whether we can define one-dimensional-higher tubes from  virtual 2-knots 
since we succeed in the virtual 1-knot case 
as written in \S\S\ref{E}-\ref{Proof}. \\

The following facts let it be more natural: 
The one-dimensional-higher tube $\mathcal E(K)$ made from a virtual 1-knot $K$ 
is the spun-knot of $K$ if $K$ is a classical knot (see \cite{Satoh}). 
\cite{Zeeman} defined  spun-knots not only for classical 1-knot but also for classical 2-knots. \\

\begin{que}\label{North Carolina}
Can we define one-dimensional-higher tubes for  virtual 2-knots 
in a consistent way?  Suppose that these tubes are diffeomorphic to 
$F\x S^1$ if the virtual 2-knot is defined by $F$. 
\end{que}\bigbreak

Note that Satoh's method in \cite{Satoh} did not say anything about the virtual 2-knot case.  
In the virtual 1-knot case, in \cite{Rourke},  
Rourke  interpreted Satoh's method  as we reviewed 
in Theorem \ref{Montana} and Definition \ref{Nebraska}. \\  

\begin{note}\label{kaiga}
In the virtual 2-dimensional knot case 
we also 
use the terms `fiber-circle' and `Rourke-fibration' 
in Definition \ref{Nebraska}.   
\end{note}

\begin{defn}\label{suiri}
Let $M$ be a 3-dimensional compact submanifold of $\R^5$. 
Regard $\R^5$ as $\R^3\x\R^2$. 
We say that the submanifold $M$ admits 
{\it Rourke fibration}, or 
that $M$  is embedded {\it fibrewise} 
if    
$M\cap(p\x\R^2)$ is a collection of circles 
for any point $p\in\R^3$. 
We call the circles in $M\cap(p\x\R^2)$, {\it fiber circles}.   
\end{defn}

\begin{figure}
\includegraphics[width=110mm]{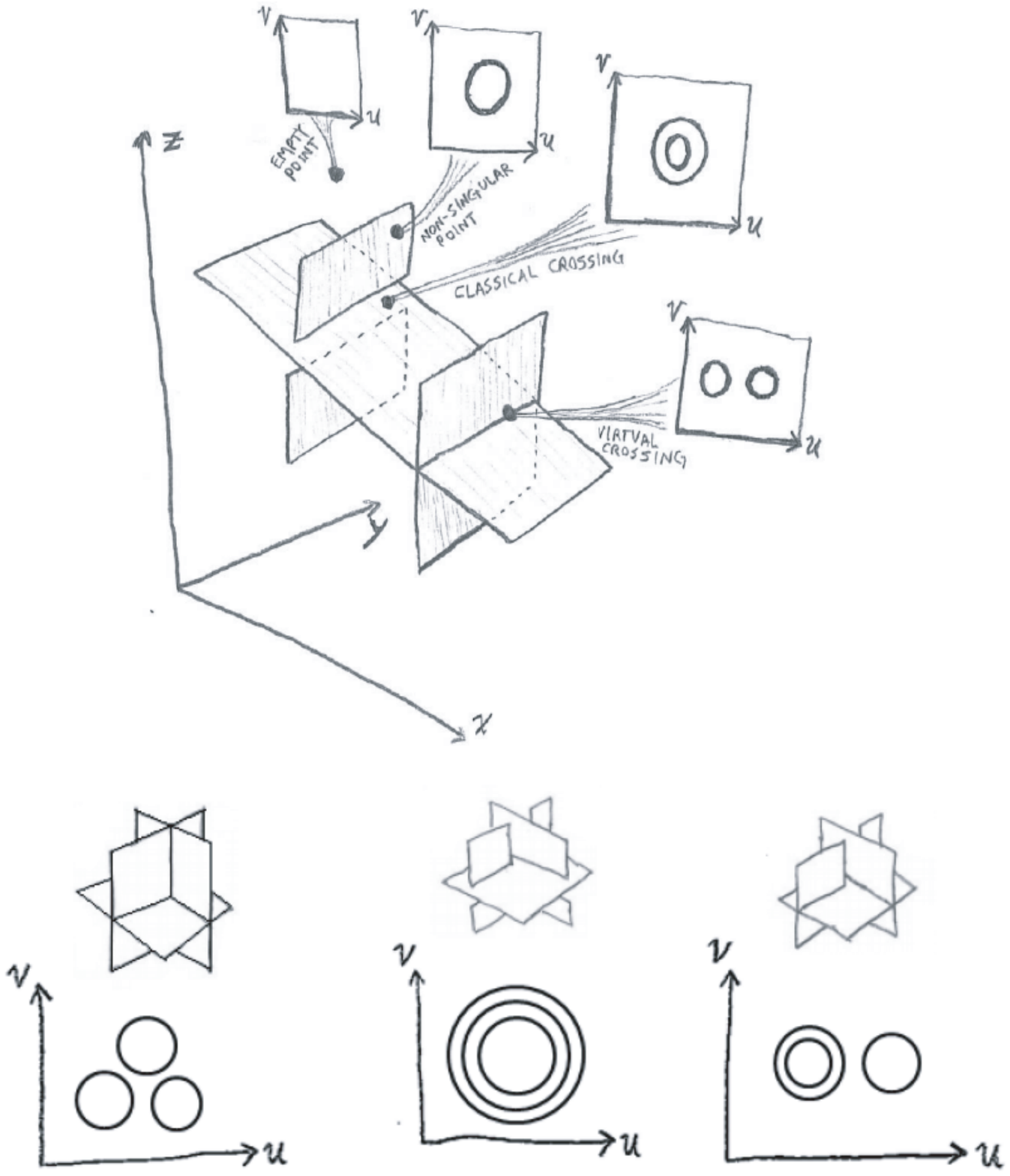}  
\vskip-20mm
\caption{{\bf 
The nest of circles in fibers.
}\label{doremi}}   
\end{figure}

If we try to generalize Rourke's way to the virtual 2-knot case, 
 we will do the following: 
Let $\alpha$ be a virtual 2-knot diagram. 
Let $\mu=0,1,2,3$. 
We give $\mu$-copies of circle to any $\mu$-tuple point in $\alpha$,  
and construct the tube. 
Of course we determine the position of fiber-circles in each fiber plane 
by the property of the $\mu$-tuple point (See \cite[section 3.7.1]{J} for detail). 
See Figure \ref{doremi}. \\

However we encounter the following situation.  
Let $\alpha$ be any virtual 1-knot diagram.  

\smallbreak\noindent(1) 
The case where $\alpha$ has no virtual branch point.  

\smallbreak\noindent(2) 
The case where $\alpha$ has a virtual branch point.  

\smallbreak
In the  (1) case, we can make a tube by Rourke's way. 
See \cite[section 3.7.1]{J}.

In the (2) case, however, \cite{J} found 
it difficult to define a tube near any virtual branch point. 
Thus it is natural to ask the following two questions. \\

\begin{que}\label{North Dakota}
Can we put fiber-circles  over each point of any virtual 2-knot 
in a consistent way as written above, 
and make one-dimensional-higher tube?
\end{que}\bigbreak

\begin{que}\label{South Dakota}
Is there a one-dimensional-higher tube construction which is defined for all virtual 2-knots, and which agrees with the way in the (1) case  written above   
when there are no virtual branch points? 
\end{que}\bigbreak

We generalize  our method in \S\S\ref{E}-\ref{Proof} 
and  
give an affirmative answer to Question  \ref{South Dakota}, 
and hence to Question   \ref{North Carolina}.  
See Theorem \ref{vv}. 
We also use a spinning construction of submanifolds 
explained in Definition \ref{spinningsubmanifold}.  
Theorem \ref{Rmuri} gives a negative answer 
to Question \ref{North Dakota}.

\bigbreak  
We make the virtual 2-knot version of 
representing surfaces, which are defined above Theorem \ref{vk}. 

\begin{defn}\label{Jbase}   {\bf (\cite[section 3.5]{J}.)}
The development of an invariant for virtual 2-knot theory closely parallels that for virtual 1-knot theory.  The idea is to think of a virtual 2-knot diagram as a classical 2-knot diagram ``drawn" on a closed 3-manifold.  We then define an equivalence relation on these objects that extends classical move-equivalence and allows the 3-manifold to vary. Take as input a virtual 2-knot diagram $\alpha$.  
Let $N(\alpha)$ be a neighborhood of the diagram, 
which is a regular neighborhood except at virtual branch points, 
in the following sense:  
$N(\alpha)$ is formed by thickening $\alpha$ everywhere except at virtual branch points; 
as you approach virtual branch points, let the thickening gradually diminish to zero, 
so that near the virtual branch point $N(\alpha)$ looks like the cone over a thickened figure-$\infty$. \\

Along each virtual crossing curve of $\alpha$, double the square-shaped junction of $N(\alpha)$ to create overlapping ``slabs".  
Call this 3-manifold-with-boundary $B(\alpha)$.  
It has a purely classical knot diagram in it.  (To be precise, $B(\alpha)$ is not technically a 3-manifold-with-boundary at virtual branch points, since the ``slab" is pinched to zero thickness at these points.)  
See Figures \ref{Jbase2} and \ref{Jbase3}. \\

\begin{figure}
\includegraphics[width=130mm]{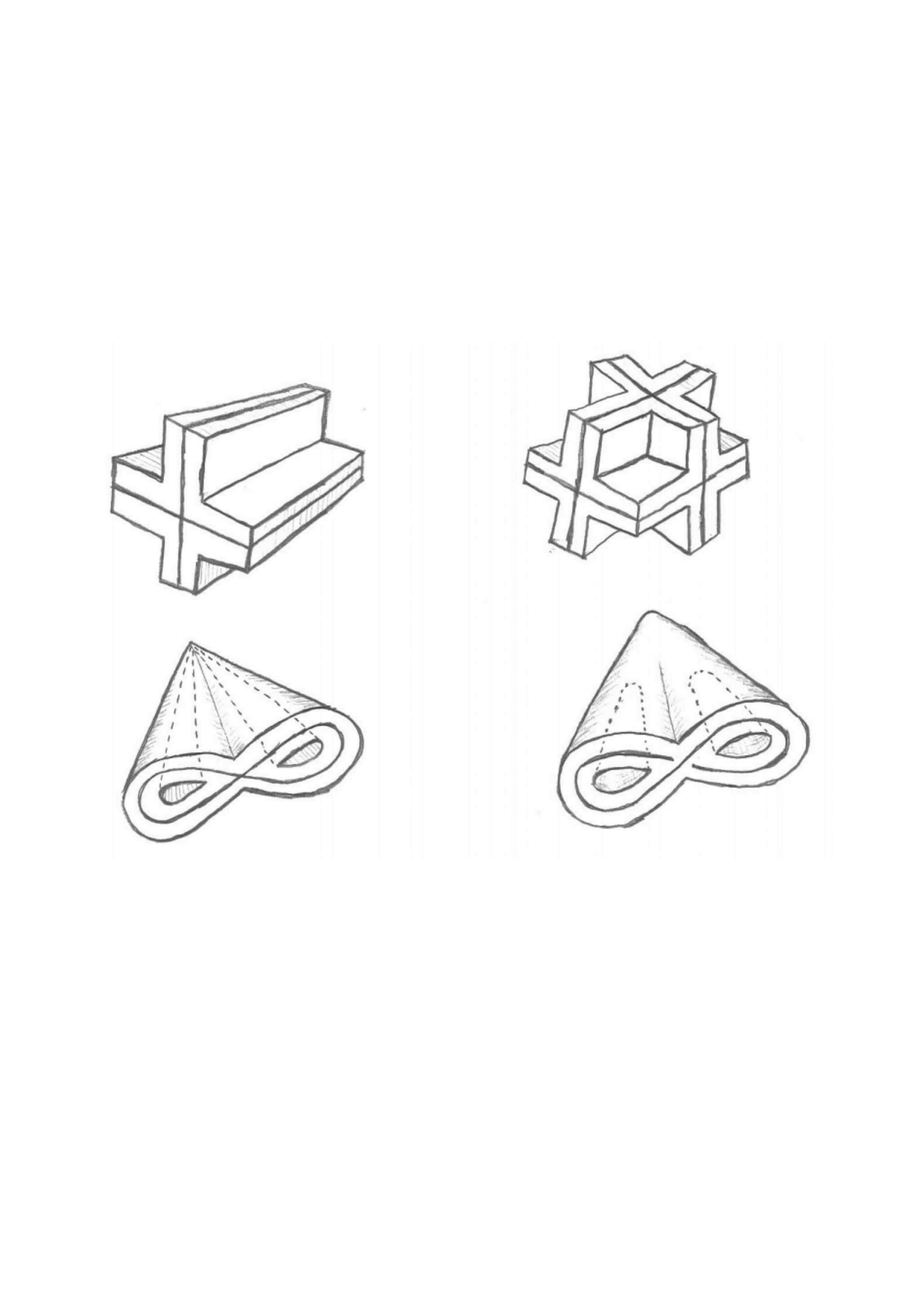}  
\vskip-30mm
\caption{{\bf 
The left upper figure is 
a part of $N(\alpha)$ near a double point curve. 
The right upper figure is that 
 near a triple point. 
The left lower figure is that 
 near a virtual branch point. 
The right lower figure is that 
 near a classical branch point. 
}\label{Jbase2}}   
\end{figure}

Now embed $B(\alpha)$ into any compact oriented 3-manifold (not necessarily connected).  
The resultant compact 3-manifold  
is  called a {\it representing 3-manifold} $M$   
associated with a virtual 2-knot diagram $\alpha$. 
$M$ contains a classical 2-knot diagram $\mathcal I(\alpha)$. 
\\

\begin{figure}
\bigbreak  
\includegraphics[width=130mm]{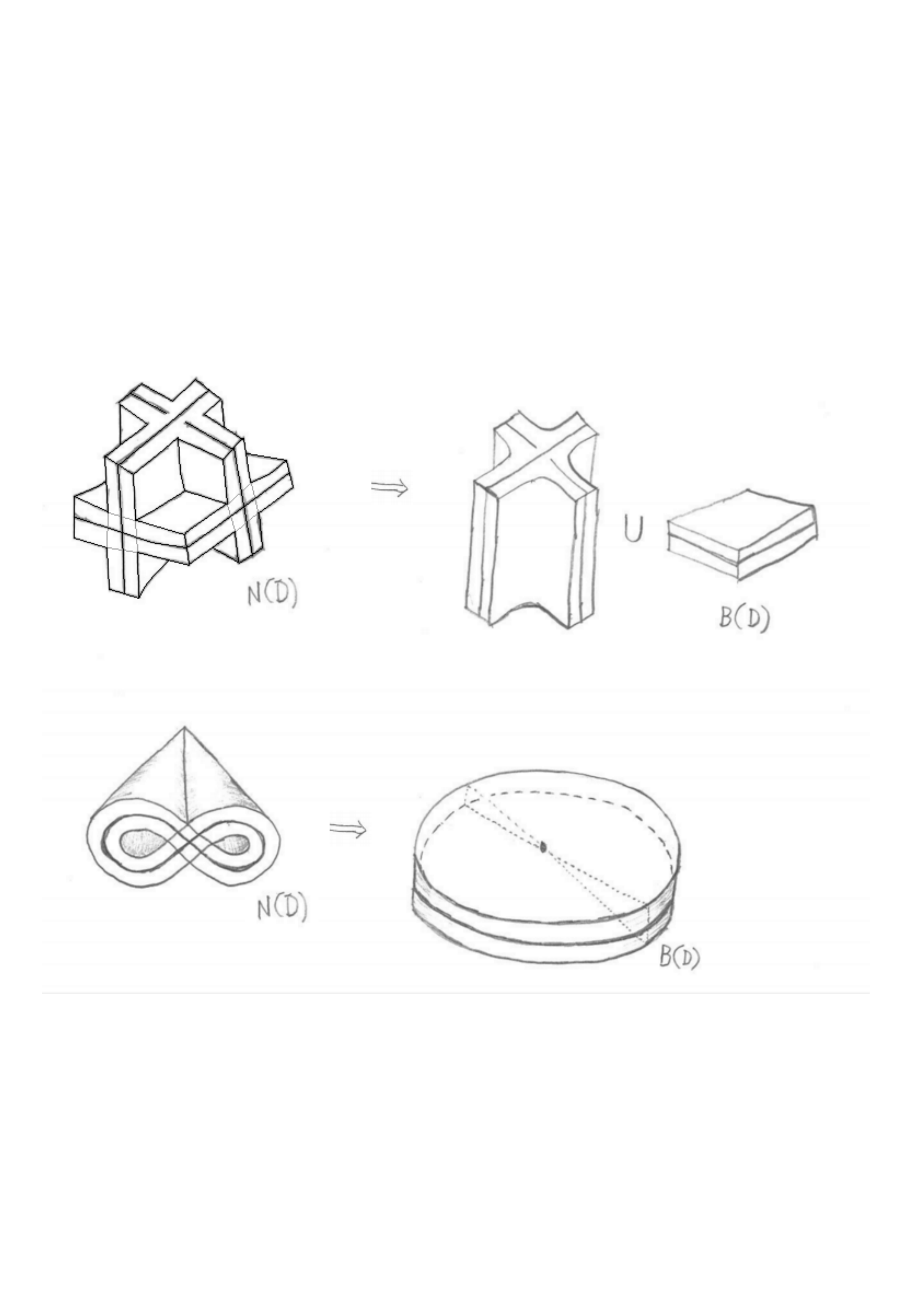}  
\vskip-20mm
\caption{{\bf 
Make $B(\alpha)$ from $N(\alpha)$.
}\label{Jbase3}}   
\bigbreak  
\end{figure}
\end{defn}

\bigbreak
\noindent{\bf Note.} 
Virtual 1-knot has two kinds of equivalent definitions: 
one is defined by using diagrams with virtual points in $\R^2$. 
The other is done by using representing surfaces. 
See Theorem \ref{vk} and \S\ref{K} of this paper, and \cite{Kauffman1, Kauffman, Kauffmani}. 
It is very natural to ask the following question. 

\begin{que}\label{imayatteru}
Do we have the virtual 2-knot version of Theorem \ref{vk} by using representing 3-manifolds in Definition \ref{Jbase}?
\end{que}

This question is open. \cite{J} gave a partial answer. 
We do not discuss it in this paper.

\vskip9mm

Note that $\mathcal I(\alpha)$ is an  immersed surface in an ordinary sense. 
That is, it does not include a virtual point. 
Note that we cannot embed a 
representing 3-manifold 
in $\R^4$ in general. 
We show an example. 
Take the Boy surface in $\R^3$
(see \cite{Boy, OgasaBoy}). 
We can regard it as 
a virtual 2-knot diagram as follows: 
Suppose that the only one immersed crossing curve 
is a virtual one. That is, it consists of 
one virtual triple point and other virtual double points. 
Then no 
representing 3-manifolds for this virtual 2-knot 
can be embedded in $\R^4$. It is 
proved by using obstruction classes of 
the normal bundle of $\R P^2$ in $\R^4$. \\

However, by \cite{Hirsh}, we have the following.

\begin{thm}\label{Hirsh}
Any $M$ in Definition \ref{Jbase} can be embedded in $\R^5$. 
\end{thm}

We will define the virtual 2-knot version of 
 $\mathcal E(\alpha)$ in Definition \ref{wakeru} after some preliminaries.    
Let $X$ be a 3-dimensional closed oriented abstract manifold. 
Let $G_1$ and $G_2$ be  submanifolds of $S^5$ 
which are diffeomorphic to $X$. 
Recall the following fact 

\begin{cla}\label{kantan}
There is a case that 
the submanifolds, $G_1$ and $G_2$, of $S^5$ are non-isotopic. 
\end{cla}

\noindent{\bf Proof of Claim \ref{kantan}.} 
Take a spherical 3-knot $K$ whose Alexander polynomial is nontrivial. 
Let $G_2$ be the knot-sum of $G_1$ and $K$. 
See e.g \cite{Rolfsen} for the Alexander polynomial of 3-knots and the knot-sum. 
\qed
\\

While $G_1$ and $G_2$ may be non-isotopic 
submanifolds, which are diffeomorphic to $X$ above, 
of $S^5$, 
it is the case that they are isotopic 
after removing an open three-ball from each of them.
Let $X^\circ_i$ denote $X-\text{(an open 3-ball)}$. 
Let $G^\circ_i=G_i-\text{(an open 3-ball)}$ be a submanifold of $S^5$ ($i=1,2$).

\begin{cla}\label{wowwow} 
The submanifolds, $G^\circ_1$ and $G^\circ_2$, of $S^5$ are isotopic. 
\end{cla}

\noindent{\bf Note.}
 Claim \ref{wowwow} is  the virtual 2-knot version of Claim \ref{wow}. 
\\

\noindent{\bf Proof of Claim \ref{wowwow}.} 
$X^\circ$ has a handle decomposition 
which consists of one 0-handle, 1-handles, 2-handles and no 3-handle.  
The dimensions of the cores of these handles are 0, 1, or 2. 
Hence the dimensions $\leqq2$. (Here, it is important the dimension $\neq$ 3.)  The dimension of $S^5$ is 5. 
Since $2<\frac{3(2+1)}{2}$, Claim \ref{wowwow} holds by  
 \cite{Haefliger3}.  
\qed

\begin{cla}\label{jimeisoku}
Let $M$ be a compact 3-manifold. By Theorem $\ref{Hirsh}$, $M$ is embedded in $\R^5$.  
The normal bundle $\nu$ of $M$ embedded in $\R^5$ is the trivial bundle for any embedding of $M$ in $\R^5$. 
\end{cla}

\noindent{\bf Proof of Claim \ref{jimeisoku}.} 
If $M$ is closed, $M$ bounds a Seifert hypersurface $V$ in $S^5$ (See \cite[Theorem 2 page 49]{Kirby}). 
Take the normal bundle $\alpha$ of $V$ in $\R^5$. 
Then $\nu$ is a sum of vector bundles $\alpha|_M$ and an orientable $\R$-bundle over $M$.    
Hence Claim \ref{jimeisoku} holds in this case. \\

In the case where $M$ is nonclosed,
take the double $DM$ of $M$ as abstract manifolds. 
Then $DM$ can be embedded in $\R^5$. 
By the previous paragraph, the normal bundle of this embedded $DM$ is trivial.  
Then the restriction of this normal bundle to $M\subset DM$ is trivial.  
By this fact and Claim \ref{wowwow},   Claim \ref{jimeisoku} holds in this case.   
This completes the proof of  Claim \ref{jimeisoku}.
\qed
\\

We introduce the virtual 2-knot version of 
 $\mathcal E(\alpha)$.

\begin{defn}\label{wakeru}
Take an abstract manifold $M$ in Definition \ref{Jbase}, 
where $\mathcal I(\alpha)$ is still contained in $M$. 
Make $M\x[0,1]$.   
We can obtain an embedded surface 
 $\mathcal J(\alpha)$ 
contained in $M\x[0,1]$ 
such that the projection of  $\mathcal J(\alpha)$   
by the projection $M\x[0,1]\to M$ 
is $\mathcal I(\alpha)$. 
We suppose  $\mathcal J(\alpha)\cap(M\x\{0\})=\phi$.    
Take any embedding of $M$  in $\R^5$.   
Define a submanifold  $\mathcal E(\alpha)$ contained in $S^5$ 
to be the spinning submanifold made from $\mathcal J(\alpha)$ around $M$. 
(Recall Claim \ref{jimeisoku}.)
\end{defn}

We prove the virtual 2-knot version of  
Theorem \ref{honto},  
which is Theorem \ref{vv}.     

Theorem \ref{vv} is one of our main results. 
It gives an affirmative answer to Question 
\ref{North Carolina}.

\begin{thm}\label{vv}
 Let $\alpha$ and $\alpha'$ be virtual 2-knot diagrams 
which represent the same virtual 2-knot. 
Make $\mathcal E(\alpha)$ and $\mathcal E(\alpha')$ 
by using a 
representing 3-manifold $M$ $($respectively, $M')$  associated with 
$\alpha$ $($respectively, $\alpha').$  
Then submanifolds, $\mathcal E(\alpha)$ and $\mathcal E(\alpha')$, of $\R^5$ 
are isotopic 
even if $M$ is not diffeomorphic to $M'$.
\end{thm}

\noindent
{\bf Proof of Theorem \ref{vv}.}   
It suffices to prove the following two cases: 

\smallbreak
(i) $\alpha$ is obtained from $\alpha'$ by one of classical moves.  

\smallbreak
(ii) $\alpha$ is obtained from $\alpha'$ by one of virtual moves.  

\smallbreak
In the case (ii), there is a diffeomorphism map $f:M\to M'$ 
such that $f(\alpha)$ is isotopic to $\alpha'$ in $M'$.  
Note that $\alpha\subset M$ and that $\alpha'\subset M'$.\\ 

In the case (i). Take a closed 3-ball $B$ where the classical move is carried out. 
Note that $M\cup B$ (respectively, $M'\cup B$) is a representing 3manifold of 
$\alpha$ (respectively $\alpha'$). 
Note that there is a diffeomorphism map $f:M\cup B\to M'\cup B$ 
such that $f(\alpha)$ is isotopic to $\alpha'$ in $M'\cup B$.  
Note that $\alpha\subset M\cup B$ and that $\alpha'\subset M'\cup B$. 


In both cases, by the following Theorem \ref{bdy}, 
Theorem \ref{vv} holds.   
\qed\\

We prove the following Theorem \ref{ohoh}, 
which is the virtual 2-knot version of Theorem \ref{oh}.  
The key idea of the proof is Claim \ref{wowwow} 
(recall Note below Claim \ref{wowwow}.)  
Let $i=1,2$.  
Take $G_i$ defined in Claim \ref{wowwow}.   
We can regard the tubular neighborhood of $G_i$ in $S^5$ as $G_i\x D^2$. 
Embed a closed oriented surface 
in  $G_i\x [0,1]$, 
where we regard $[0,1]$ as a radius of $D^2$,  
and call the image $J_i$.  
Assume that $J_i\cap(G_i\x\{0\})=\phi.$
Suppose that there is a bundle map 
$\check\sigma:G_1\x D^2\to G_2\x D^2$ 
such that $\check\sigma$ covers an orientation preserving diffeomorphism map $\sigma:G_1\to G_2$ 
and such that $\check\sigma(J_1)=J_2$. 
Define a submanifold $E_i$ contained in $S^5$ 
to be the spinning submanifold made from $J_i$ 
 by the rotation in $G_i\x D^2$. 

\begin{thm}\label{ohoh}
The submanifolds, $E_1$ and $E_2$, of $S^5$ are isotopic. 
\end{thm}

\noindent{\bf Proof of Theorem \ref{ohoh}.} 
We can suppose that $J_i\subset G^\circ_i\x [0,1]$. 
By the existence of $\sigma$,  
there is a bundle map $\check\tau:G^\circ_1\x D^2\to G^\circ_2\x D^2$ 
such that $\check\tau$ covers a diffeomorphism map 
$\tau:G^\circ_1\to G^\circ_2$ 
and such that $\check\tau(J_1)=J_2$.

Note the following: 
Let $f:M^\circ\to S^5$ be an embedding map. 
We can regard $\tau$ as a diffeomorphism map 
$M^\circ\to M^\circ$. 
By Claim \ref{wowwow},  
the submanifolds, $f(M^\circ)$ and $f(\tau(M^\circ))$, of $S^5$ are isotopic. 
Therefore 
the submanifolds, $E_1$ and $E_2$, of $S^5$ are isotopic. 
\qed

\begin{thm}\label{bdy}
Replace 
the condition that $M$ is a closed compact oriented 3-manifold 
with 
the condition that $M$ is a non-closed compact oriented 3-manifold.  
Then Theorem \ref{ohoh} also holds. 
\end{thm}

\noindent{\bf Proof of Theorem \ref{bdy}.} 
The proof of Theorem \ref{bdy} 
is done in a similar fashion to that of Theorem \ref{ohoh}.  
The proof of Theorem \ref{bdy} is easier than that of Theorem \ref{ohoh}.  
\qed\\

We now have completed the proof of Theorem \ref{vv}, and 
answered Question 
\ref{North Carolina}.

We next answer Question \ref{South Dakota}.

We define a consistent way 
to put a 
representing 3-manifold in $\R^5$.

\begin{defn}\label{hatena} 
Let $\alpha$ be a virtual 2-knot diagram 
contained in $\R^3$. 
Regard $\R^3$ as 
$\R^3\x\{0\}\x\{0\}$
$\subset\R^5=\R^3\x\R\x\R$. 
Put `a 
representing 3-manifold 
for $\alpha$' in $\R^5$ as follows.

Take the neighborhood $T$ of $\alpha$ as 
defined in Definition \ref{Jbase}. 
Take a neighborhood of each of classical and virtual branch points 
such that the neighborhood  is diffeomorphic to the closed 3-ball 
and 
such that $\alpha\cap$(the neighborhood) is as drawn in Figure \ref{cvbr}.  \\
\begin{figure}
 \includegraphics[width=40mm]{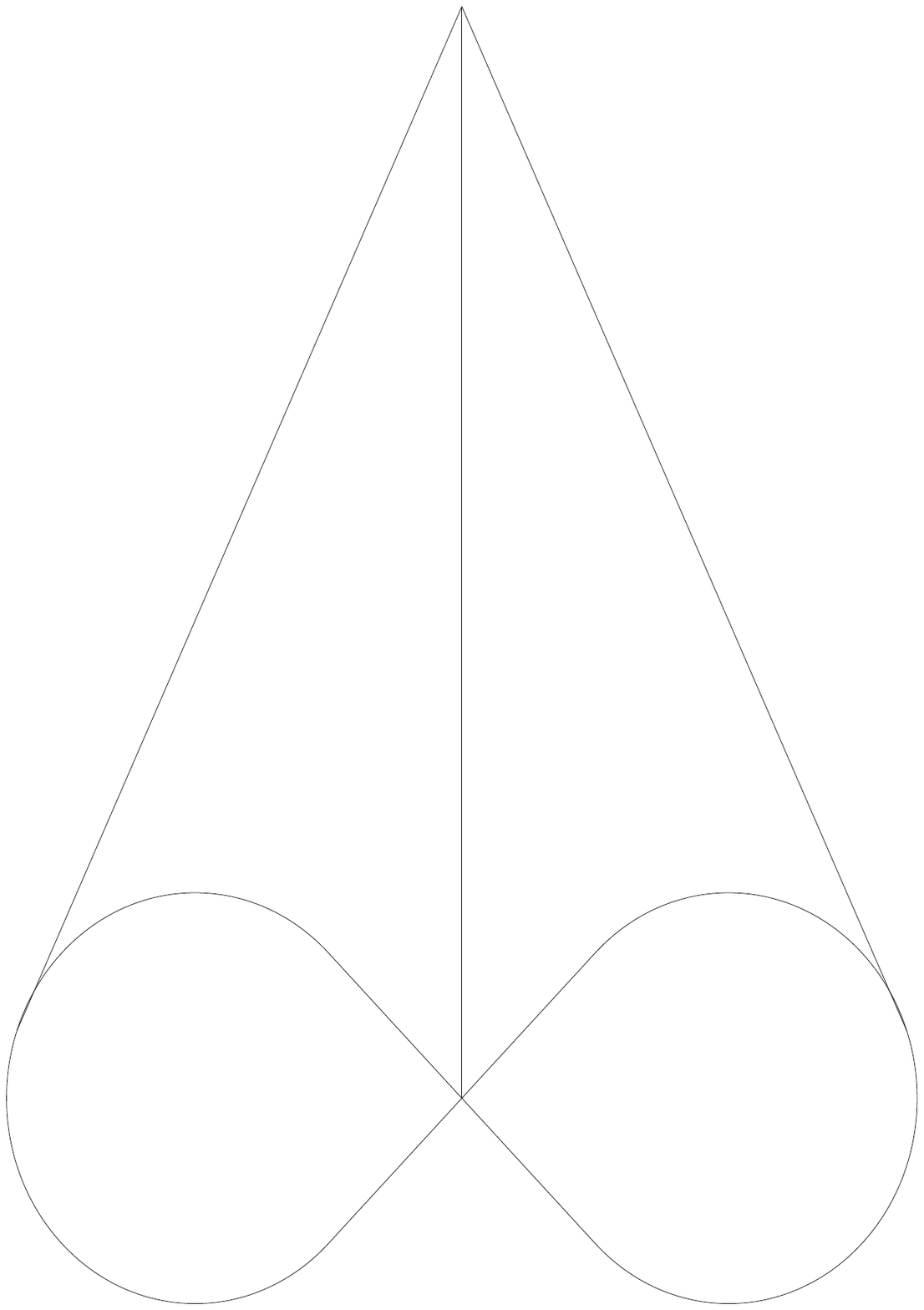}  
\caption{{\bf 
The intersection of 
a classical or virtual branch point
and 
its neighborhood explained in Definition \ref{hatena}
}\label{cvbr}}   
\end{figure}
Let  
$T'
=T-$Int(the neighborhoods of real branch points and those of virtual branch ones).  
Along any virtual crossing line 
we double $T'$ 
as done in Definition \ref{Jbase}. 
Note that 
this operation can be done in 
$\R^5$    
although it cannot be done in $\R^4$ in general. 
Thus we obtain a compact oriented 
3-dimensional submanifold $X\subset\R^5$ from $T'$. \\

For a real branch point, we attach 
`the closed 3-ball which is a neighborhood of the real branch point' 
to $X$. 
Note that
near any virtual branch point, 
 the operation can be done in 
$\R^3\x\R\x\{0\}$. 
For a virtual branch point,  we attach 
`the closed 3-ball which is a neighborhood of the virtual branch point',  as drawn in 
Figures \ref{Delaware}-\ref{Hawaii},  
to $X$. 
Note that 
in Figures \ref{Delaware}-\ref{Hawaii} 
we draw $\R^4=\R^3\x\R\x\{0\}$. 
Note that the virtual branch point vanishes in this closed 3-ball. \\

The resultant compact oriented 3-manifold is 
a 
representing 3-manifold with $\mathcal I(\alpha)$, which is defined in Definition \ref{Jbase}. 
We call it $M_\iota$. 
Recall that $\mathcal I(\alpha)$ has no virtual point and, in particular,  
that $\mathcal I(\alpha)$ has no virtual branch point.
\\

Figure \ref{Delaware}  draws a part of a 
representing 3-manifold $M$ near a virtual branch point. 

Figure \ref{Florida} adds a part of 
$\mathcal I(\alpha)$ 
to Figures \ref{Delaware}.  

Figure \ref{Georgia}  
draws Figures \ref{Delaware} 
by seeing from a different direction.

Figure \ref{Hawaii} 
draws Figures \ref{Florida} 
by seeing from a different direction.

\begin{figure}
\bigbreak
     \includegraphics[width=140mm]{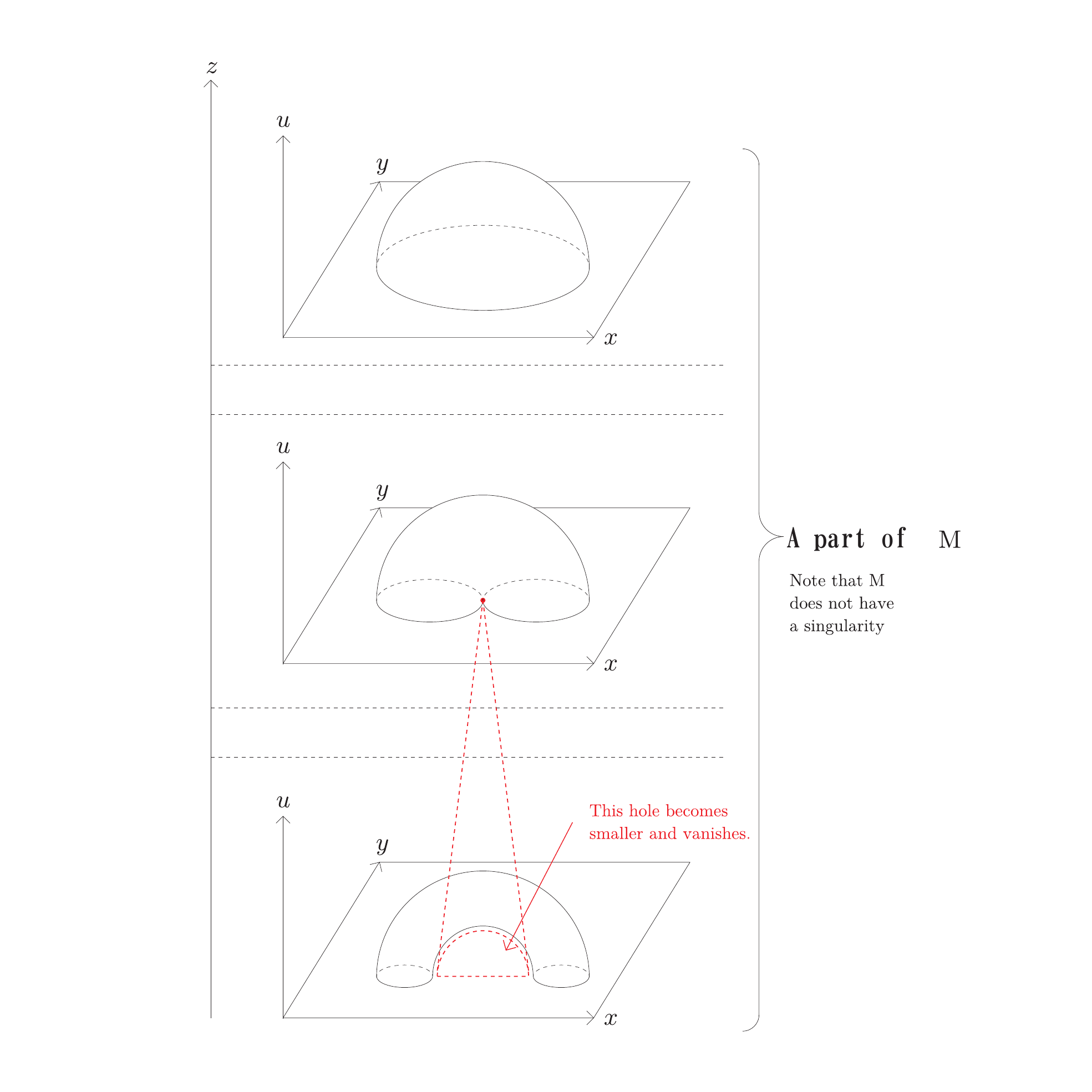}  
\caption{{\bf  
A part of a 
representing 3-manifold near a virtual branch point. We do not draw  
a virtual branch point here. 
In \ref {Florida} we do it. 
}\label{Delaware}}   
\bigbreak  
\end{figure}

\begin{figure}
\bigbreak
   \includegraphics[width=145mm]{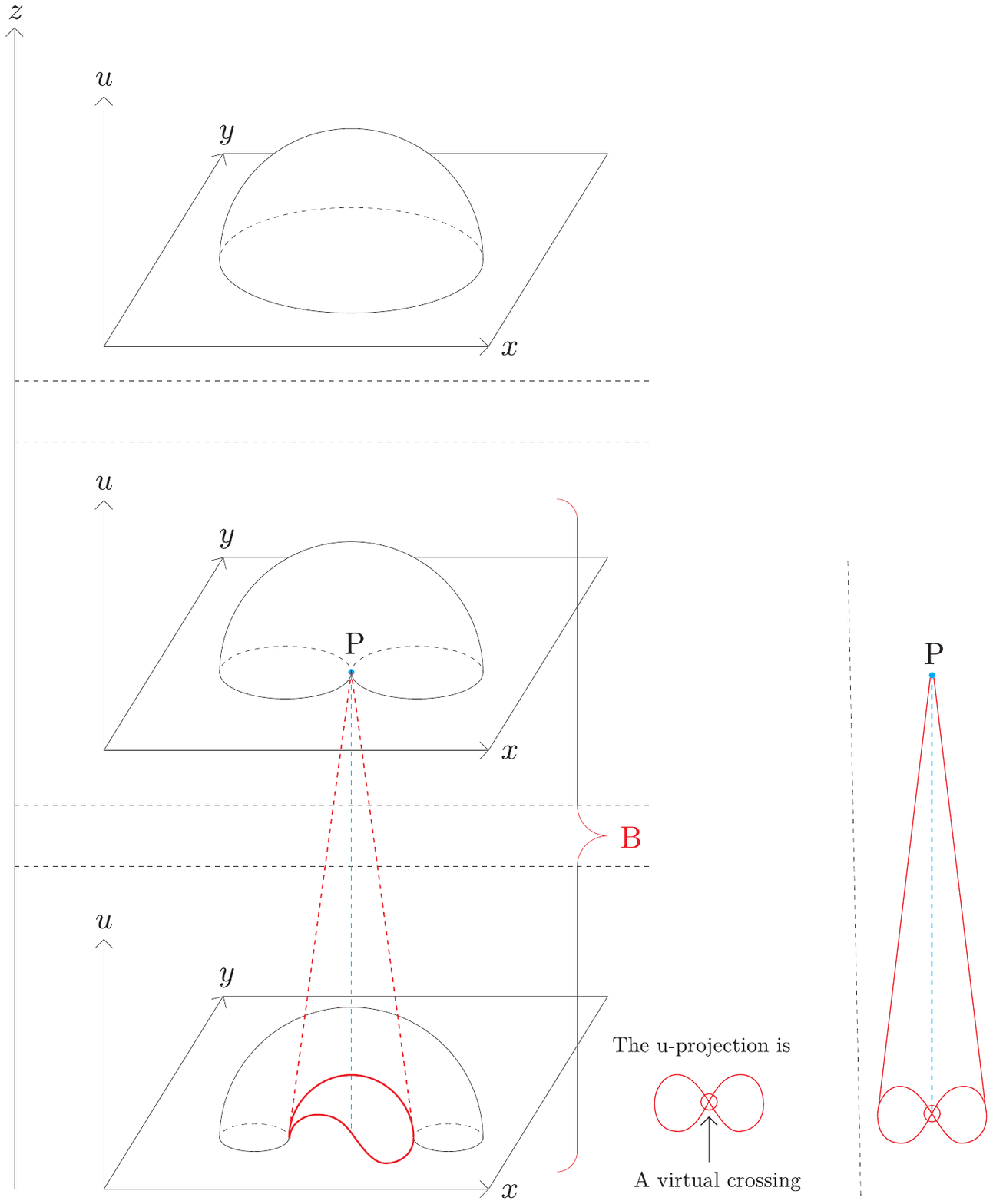}  
\vskip-30mm
\caption{{\bf  
A part of a 
representing 3-manifold near a virtual branch point. We draw  
a virtual branch point here. 
}\label{Florida}}   
\bigbreak  
\end{figure}

\begin{figure}
\bigbreak
     \includegraphics[width=140mm]{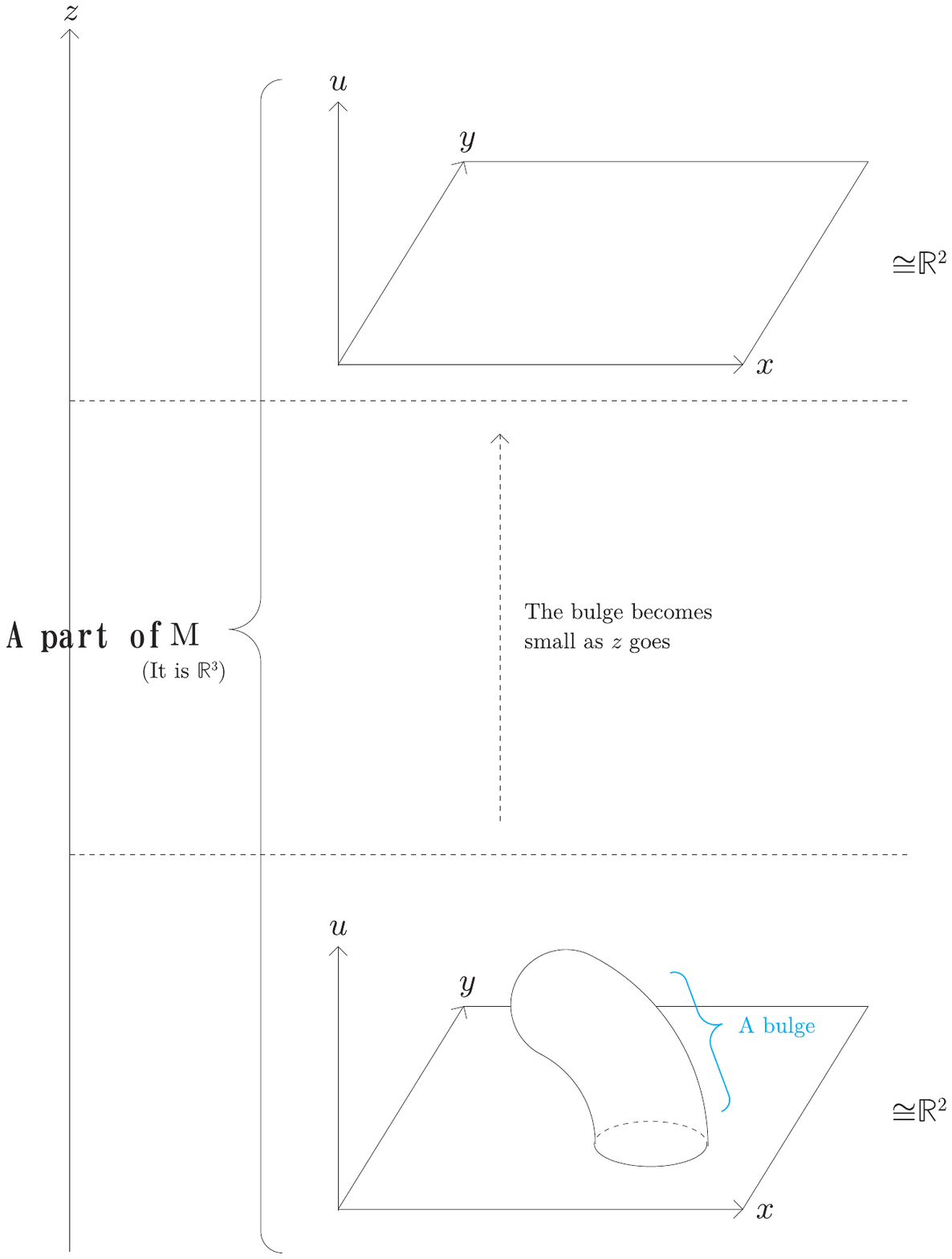}  
\vskip-40mm
\caption{{\bf  
A part of a 
representing 3-manifold near a virtual branch point. We do not draw  
a virtual branch point here. 
In Figure \ref{Hawaii} we will do it. 
}\label{Georgia}}   
\end{figure}

\begin{figure}
     \includegraphics[width=140mm]{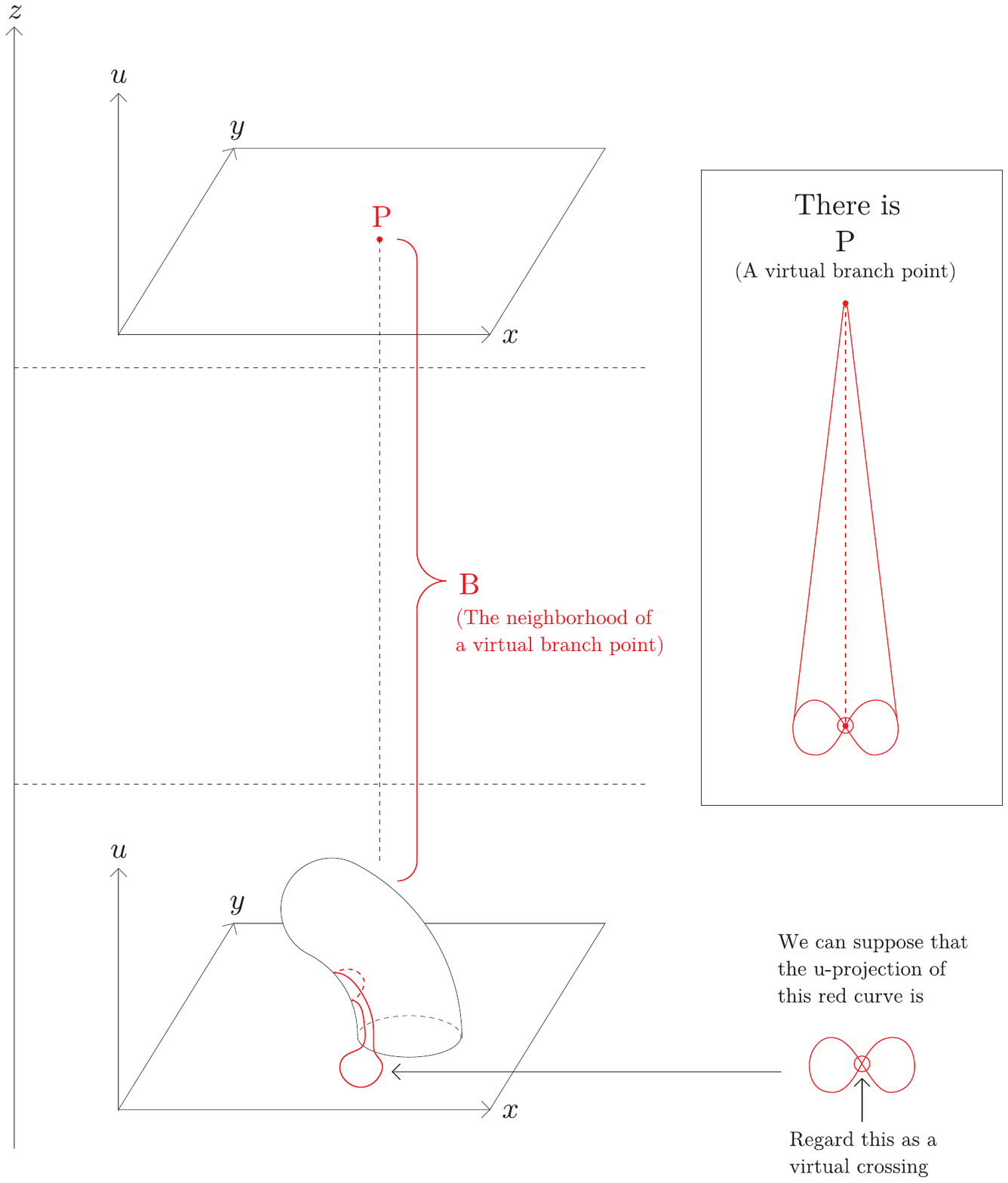}  
\vskip-30mm
\caption{{\bf  
A part of a 
representing 3-manifold near a virtual branch point. We draw  
a virtual branch point here. 
We explain the most lower figure 
in more detail in Figure \ref{hosoku}.    
}\label{Hawaii}}   
\bigbreak   
\end{figure}

\begin{figure}
\includegraphics[width=110mm]{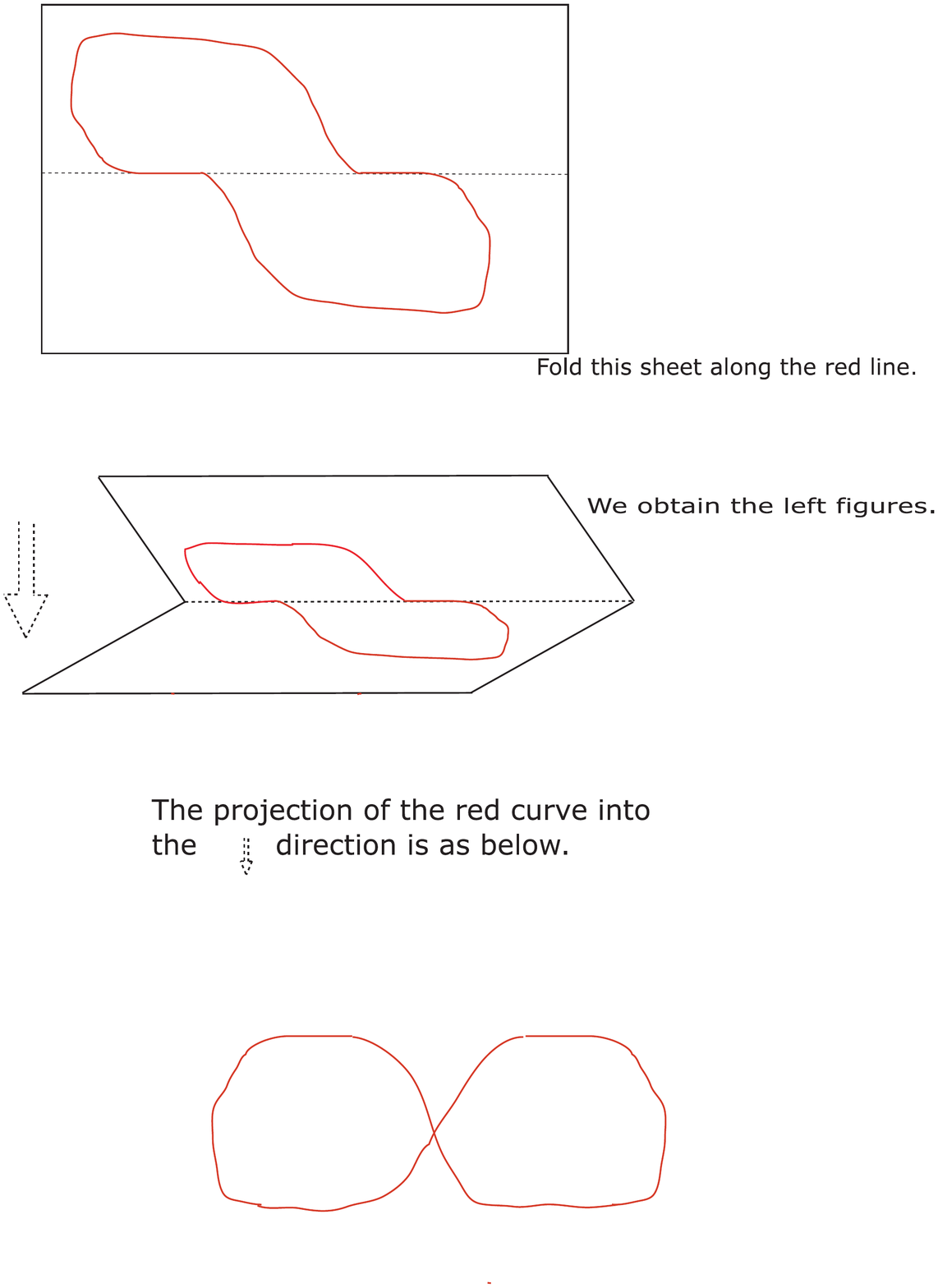}  
\vskip-10mm
\caption{{\bf  
The explanation of the most lower figure of 
Figure \ref{Hawaii}. 
The neighborhood of the red curve of that figure is 
obtained by curving the middle figure of the above figures. 
}\label{hosoku}}   
\end{figure}
 \end{defn}

We prove that we have an affirmative answer to 
 Question \ref{South Dakota}.

\begin{thm}\label{konnyaku}
Let $\alpha$ be a virtual 2-knot diagram. 
Make $\mathcal E(\alpha)$ by using $M_\iota$,   
and call it  $\mathcal E_\iota(\alpha)$. 
If $\alpha$ has no virtual branch point, then 
$\mathcal E_\iota(\alpha)$ admits Rourke fibration. 
\end{thm}

\noindent{\bf Proof of Theorem \ref{konnyaku}.}
Let $p\in\alpha$. 
Regard $\R^5$ in Definition \ref{hatena} 
as $\R^3\x\R\x\R$. 
By the construction of $\mathcal E_\iota(\alpha)$,  
 $\mathcal E_\iota(\alpha)\cap(p\x\R\x\R)$ 
 is the empty set or a collection of circles 
such that this correspondence satisfies Rourke's description. 
Hence Theorem \ref{konnyaku} holds.  
\qed

\begin{figure}
\bigbreak  
 \includegraphics[width=130mm]{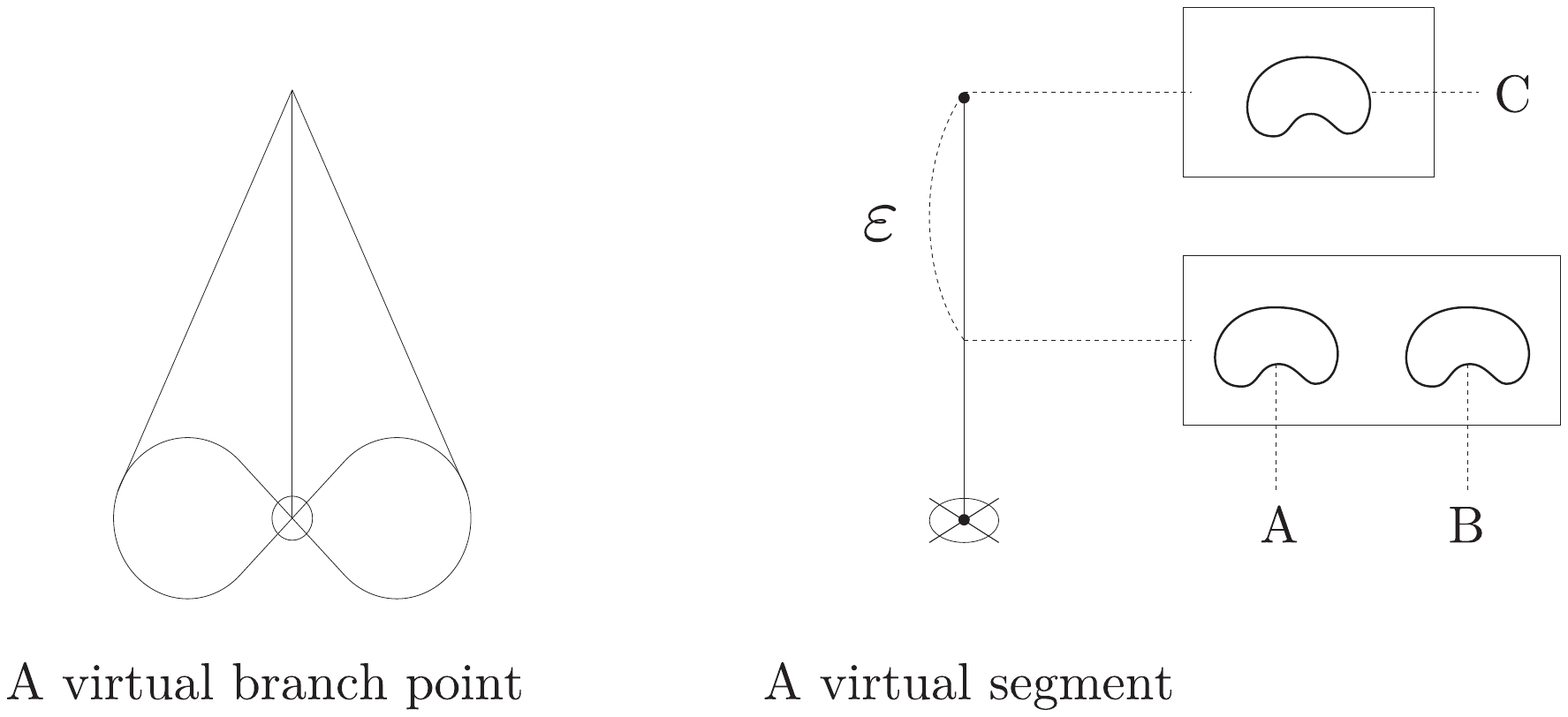}  
\vskip-60mm
\caption{{\bf 
The assumption of `reductio ad absurdum'}\label{Maryland}}   
\end{figure}

\begin{note}\label{vrei} 
It is trivial that if we use another embedding of another $M$, 
$\mathcal E(\alpha)$ associated with the embedding 
may not admit Rourke fibration. 
Such an example exists. 
Let $\xi$ be the trivial 2-knot diagram. 
It is trivial that $\xi$ admits Rourke fibration.
Let $\zeta$ be a virtual 2-knot diagram of the trivial 2-knot. 
Assume that the singular point set of $\zeta$ consists of two virtual branch points 
and one virtual segment.  \\

A {\it virtual segment} is the segment with the following properties. 
It is a segment included in a virtual 2-knot diagram.   
One of the boundary is a virtual branch point.
The points in the interior of the segment are virtual double points. 
It is drawn in Figure \ref{Maryland}. 
It is drawn in Figure \ref{sashimiv} if the branch point there is a virtual branch point.  
See \cite{J}.  \\

$\zeta$ does not admit Rourke fibration by Theorem \ref{Rmuri}. 
\end{note}

Note the following claim. 

\begin{cla}\label{shichi}
Take $\mathcal E_\iota(\alpha)$ in 
Theorem $\ref{konnyaku}$.     
If $\alpha$ includes a virtual branch point, 
$\mathcal E_\iota(\alpha)$  does not admit Rourke's fibration. 
That is,  
$\mathcal E_\iota(\alpha)$ is not embedded fiberwise.
\end{cla} 

\noindent{\bf Note.}
It is trivial that if we use another embedding of another $M$, 
$\mathcal E(\alpha)$ associated with the embedding 
may admit Rourke fibration. 
Such an example exists. 
It is the one in Note \ref{vrei}. 
\\

\noindent{\bf Proof of Claim \ref{shichi}.}   
By Theorem \ref{Rmuri}. \qed \\

We give an alternative proof of Claim \ref{shichi} after Proof of Theorem \ref{Rmuri}.  
\\

Theorem \ref{Rmuri} is an answer to Question \ref{North Dakota},  
and is one of our main results.

\begin{thm}\label{Rmuri} 
The answer to Question $\ref{North Dakota}$  
is negative. 
\end{thm}

\noindent{\bf Proof of Theorem \ref{Rmuri}.}  
We prove by `reductio ad absurdum'. 
We suppose the following assumption, and will arrive at a contradiction.

\smallbreak
\noindent
{\bf Assumption.} The neighborhood of a virtual branch point can be covered by the fiber-circles. 

\smallbreak
Note the fiber over the virtual segment as shown in Figure \ref{Maryland}. 

Give a Euclidean metric to $\R^5$.

By Assumption, 
the circles, 
$A$ and $B$ in Figure \ref{Maryland}, 
meet at the circle 
$C$ when $\varepsilon\to0$. 
Let $s$ be the area of $C$. 
When $\varepsilon\to0$, $A\to C$ and $B\to C$. 
Hence, we have the following. 
\begin{equation}\label{Massachusetts}
{\text{When $\varepsilon\to0$, (the area of $B)\to s$.}} 
\end{equation}
Note that $s$ is a fixed positive real number.

\begin{figure}
\includegraphics[width=120mm]{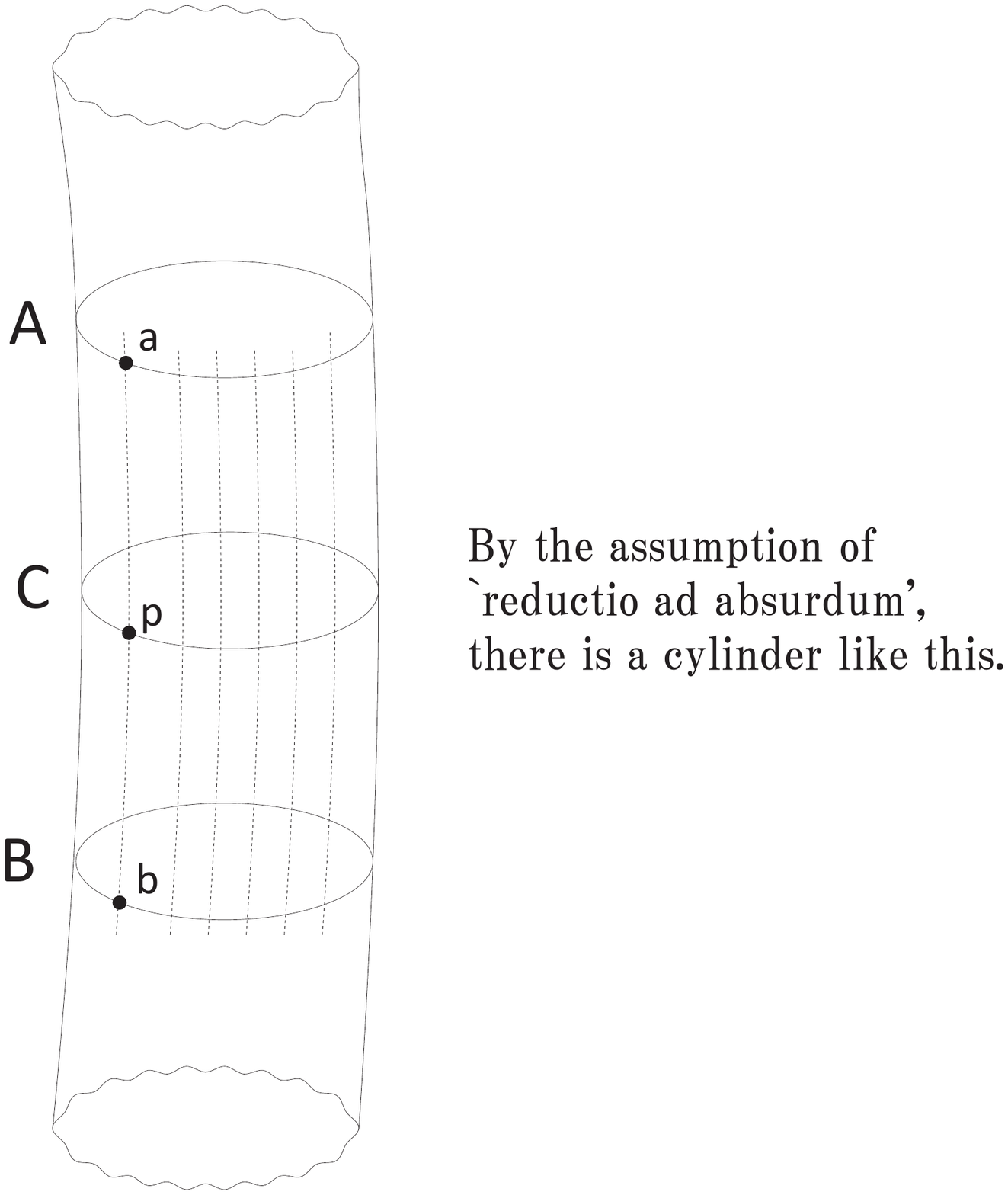}  
\vskip-30mm
\caption{{\bf A one parameter families}\label{Missouri}}   
\end{figure}

\begin{figure}
\includegraphics[width=130mm]{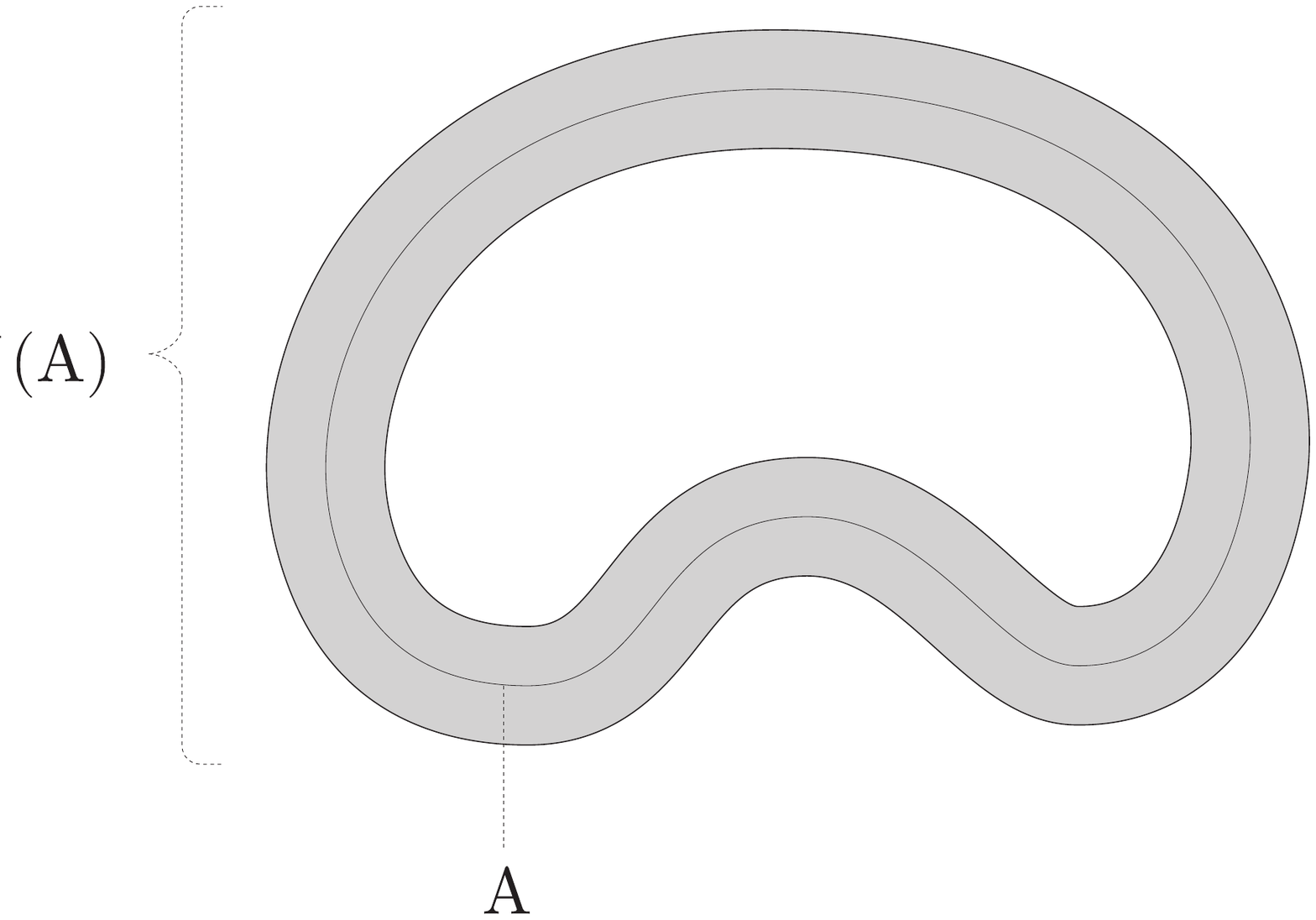}  
\vskip-70mm
\caption{{\bf $N(A)$}\label{Mississippi}}   
\end{figure}

Take a one-parameter-family for each point $p\in C.$ See Figure \ref{Missouri}. 
Suppose that $a$ and $b$ go to $p$ when $\varepsilon\to0$. 
Let $\delta(a,b)$
the distance along the trace of the one-parameter-family  
between $a$ and $b$. 
\begin{equation}\label{Michigan}
{\text{When $\varepsilon\to0$, $\delta(a,b)\to0.$}}  
\end{equation}  
\noindent
In the fiber $\R^2$ which includes $A$ and $B$,  
take any point $x\in A$. 
Suppose that $x$ goes to $y\in B$ by the one-parameter-family. 
In this fiber $\R^2$ 
take a disc of radius $2\delta(x,y)$ whose center is $x\in A$. 
Call the sum of the discs, $N(A)$. See Figure \ref{Mississippi}. 
When $\varepsilon\to0$, \\
(the area of $N(A))\to 0$. 
By (\ref{Michigan}), 
$B\subset N(A)$. 
Note that in this fiber $\R^2$, 
$B$ (respectively, $A$) is not included in the inside of $A$ (respectively, $B$).  
Therefore, by Jordan curve theorem,  
 the inside of $B$ is also included in $N(A)$. 
Hence we have the following.

\begin{equation}\label{Minnesota} 
{\text{When $\varepsilon\to0$, (the area of $B)\to 0$.}}
\end{equation}

By (\ref{Massachusetts}) and  (\ref{Minnesota}), 
we arrived at a contradiction.

This completes the proof of Theorem \ref{Rmuri}. 
\qed\\

We give a direct proof of why $\mathcal E_\iota(\alpha)$ does not admit Rourke fibration,  without using Theorem \ref{Rmuri}.
Note that it is not an alternative one of Theorem \ref{Rmuri}. 

\begin{figure}[H]
     \includegraphics[width=170mm]{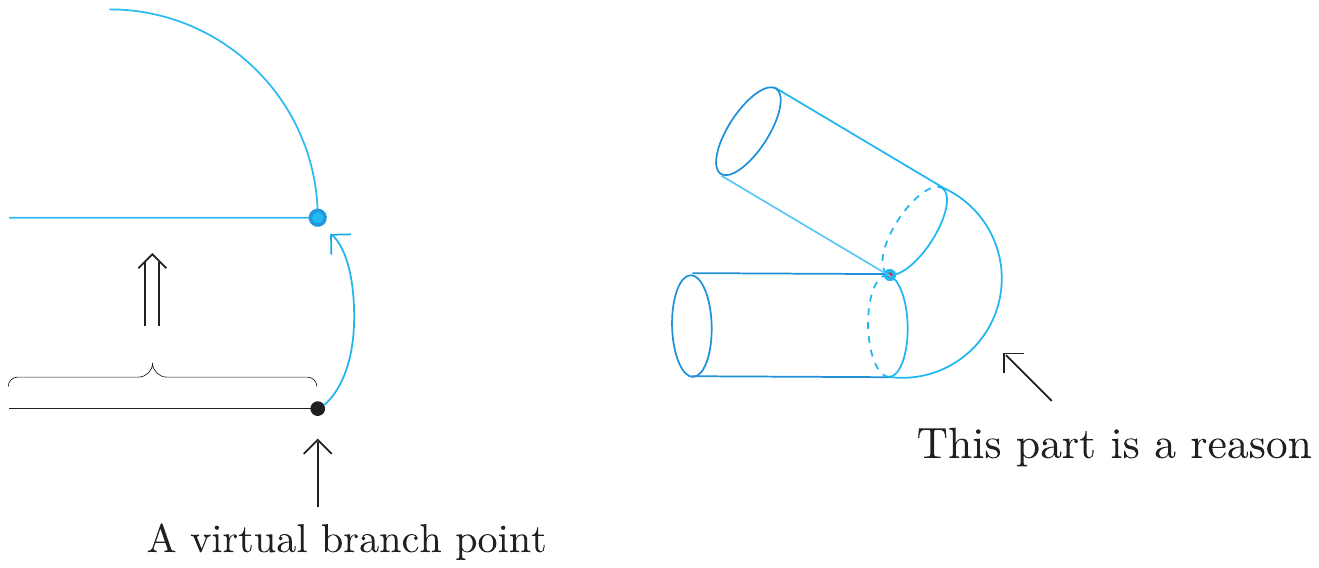}  
\vskip-90mm
\caption{{\bf 
$\mathcal E_\iota(\alpha)$ made by the spinning construction can be embed in $\R^5$ 
but  $\mathcal E_\iota(\alpha)$ does not admit Rourke fibration. 
The reason why we cannot make the fiber over any virtual branch  point a collection of circles is drawn. }\label{Idaho}}   
\end{figure}

\bigbreak
\noindent{\bf Alternative proof of Claim \ref{shichi}.}   
If $p$ is a virtual branch point, 
$p$ is in the boundary of a virtual segment in $\R^3$. 
Take $\mathcal I(\alpha)$ immersed in $M$. 
Let $\kappa:\mathcal I(\alpha)\to\alpha$ be the natural map defined in 
Definition \ref{JV}.  
We have the following.
$\kappa^{-1}$(the virtual segment) 
is a union of two segments, $\Psi$ and $\Phi$.   
A point of $\partial\Psi$ and that of $\partial\Phi$ meet at a point 
as drawn in Figure \ref{Idaho}.  
$\kappa$(this point) is the virtual branch point.  \\

The two segments make an angle. 
See Figures \ref{Delaware}-\ref{Hawaii}. 
The angle is acute. 
Even if we take an arbitrary 
representing 3-manifold of the virtual 2-knot diagram $\alpha$, 
the angle is acute not obtuse.   
Furthermore the angle is put as drawn there. 
The reason for this is that  there is always an acute angle as drawn in 
Figure \ref{Idaho} whichever representing 3-manifolds we take.  \\

As we preannounced in Notes \ref{kabuto} and \ref{kuwagata}, 
we use Figures \ref{Arizona} and \ref{Arkansas}.   
In particular, see the most lower figure of Figure \ref{Arkansas}.   

Therefore  
 $\mathcal E_\iota(\alpha)\cap(p\x\R_u\x\R_v)$ 
 is a bouquet, 
not the empty set or a collection of circles.

\np{\color{white}a}
\vskip-30mm

Therefore Claim \ref{shichi} holds. 
\qed

\bigbreak
\section{The $\mathcal E$-equivalence}\label{vw}
We introduce a new equivalence relation 
of the set of 1-(respectively, 2-)dimensional virtual knots. 

\begin{defn}\label{zoo}
Let $K$ and $J$ be 1-(respectively, 2-)dimensional virtual knots. 
If the submanifolds, $\mathcal E(K)$ and $\mathcal E(J)$, 
of $\R^4$ (respectively, $\R^5$) are isotopic, 
$K$ and $J$ are said to be 
{\it $\mathcal E$-equivalent}. 
See Theorem \ref{vv}, and the line right below Theorem \ref{honto} for $\mathcal E(\quad)$.  
\end{defn}

\begin{thm}\label{milk}
 {\rm {\bf (By \cite[Theorem 2.2]{Rourke} and  \cite[Proposition 3.3]{Satoh}.)}}
 If two virtual 1-knots are welded equivalent, 
then they are $\mathcal E$-equivalent. 
Hence there are two virtual 1-knots, $J$ and $K$,  
such that $J$ is not virtually equivalent to $K$ 
but such that $\mathcal E(J)$ is isotopic 
to $\mathcal E(K)$. 
\end{thm} 

\noindent{\bf Proof of Theorem \ref{milk}.}
By \cite[Theorem 2.2]{Rourke}, there are two virtual 1-knots, $J$ and $K$,  
such that 
$J$ is not virtually equivalent to $K$
but such that 
$J$ is welded equivalent to $K$. 
By \cite[Proposition 3.3]{Satoh},   
$J$ and $K$ are $\mathcal E$-equivalent. 
\qed\\

Thus it is natural to ask whether 
we have the virtual 2-knot version of Theorem \ref{milk}.  
In other words,  are there virtual 2-knots, $J$ and $K$,   
which are $\mathcal E$-equivalent but which are not virtually equivalent? 
We answer this question below.
\\

Let $\alpha$ be a 1-dimensional virtual knot diagram defined in $\R^2$.  
Regard $\R^3$ as the result of rotating $\R^2_{\geqq0}=\R^1\x\{t\vert t\geqq 0\}$ around  $\R^1\x\{t\vert t=0\}$ as the axis.    
Take $\alpha$ in 
$\R^1\x\{t\vert t> 0\}$. 
When we rotate $\R^2_{\geqq0}$, rotate $\alpha$ together. 
Then we obtain a 2-dimensional virtual knot diagram in $\R^3$ naturally, 
and call it $\mathcal O(\alpha)$. 
Note that 
$\mathcal O(\alpha)$ is a virtual 2-knot diagram made from $T^2$. \\

If 1-dimensional virtual knot diagrams, $\alpha$ and $\beta$,  are virtually equivalent, 
it is trivial that 
2-dimensional virtual knot diagrams, $\mathcal O(\alpha)$ and $\mathcal O(\beta)$ 
are virtually equivalent (see Definition \ref{JV}).   
Hence it makes sense that we define an 2-dimensional virtual knot $\mathcal O(K)$ for a 1-dimensional virtual knot $K$.  \\

Let $X$ be a classical surface knot contained in $\R^4=\R^3\x\{t\in\R\}$. Take $X$ in $\R^3\x\{t>0\}$. Regard $\R^5$ as the result of rotating $\R^3\x\{t\geqq0\}$ around $\R^3\x\{t=0\}$ as the axis.  Then we rotate $X$ together. 
Call the resultant 3-dimensional submanifold of $\R^5$, ${\mathcal O}(X)$. 
Note the following: If $X$ is diffeomorphic to 
a closed 
surface $\Sigma_g$, 
then 
${\mathcal O}(X)$ is diffeomorphic to 
$\Sigma_g\x S^1$. \\

\begin{pr}\label{kakan}
Let $K$ be a virtual 1-knot. 
Then 
the submanifolds,  $\mathcal E({\mathcal O}(K))$ and  ${\mathcal O}(\mathcal E(K))$,   of $\R^5$ 
are isotopic. 
\end{pr}  

\noindent{\bf Proof of Proposition \ref{kakan}.}
By the definitions. \qed
\\

The {\it standardly embedded torus} or {\it standard torus} is 
a submanifold of $\R^4$, 
diffeomorphic to the torus, 
and put in the standard position.   
Let $\Sigma_g$ be an oriented closed surface.  
The {\it standardly embedded surface 
$($diffeomorphic to $\Sigma_g)$} 
or 
{\it standard surface $($diffeomorphic to $\Sigma_g)$} 
is defined as well. 
%
%
Note that we can regard classical 1- (respectively, 2-) knots  as virtual knots. \\

\begin{figure}
\includegraphics[width=100mm]{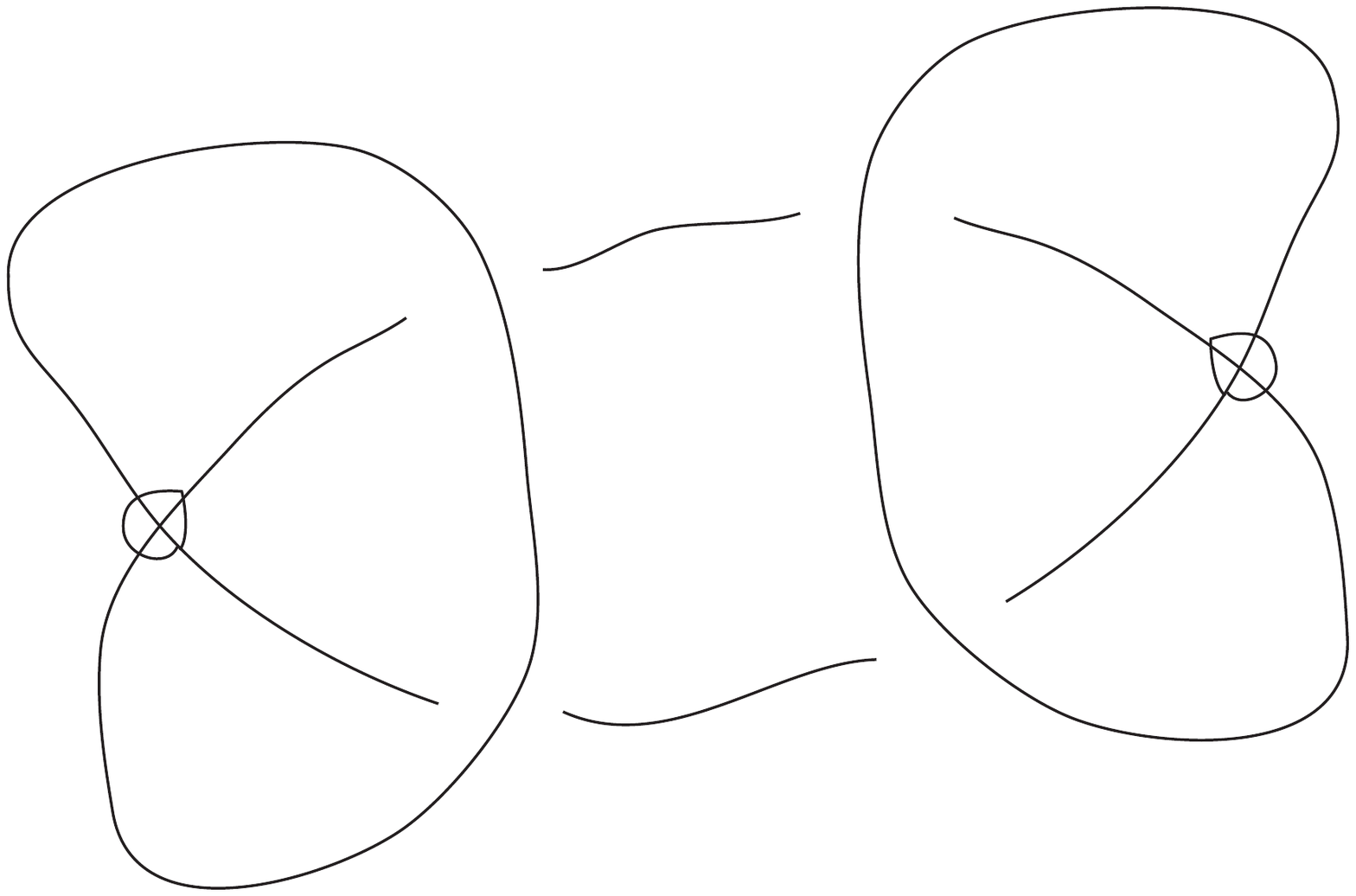}  
\vskip-40mm
\caption{{\bf Virtual reef knot}\label{vr}}   
\end{figure}

Let $R$ be the virtual reef knot  
whose diagram is drawn 
in \cite[Figure 3, section three]{Rourke}. 
We cite the diagram in Figure \ref{vr}. 
As written there, $R$ is a nontrivial virtual 1-knot,  
is welded equivalent to the trivial 1-knot, 
and has the group $\Z$. \\

\begin{cla}\label{theta} 
The virtual 2-knot $\mathcal O(R)$ 
is not virtually equivalent to 
the standard torus.  
\end{cla}

\begin{note}\label{tau}
The submanifold, $\mathcal E(R)$ and the standard torus,  
of $\R^4$ are isotopic  
because the virtual 1-knot $R$ is welded equivalent to the unknot.  
\end{note}

\noindent{\bf Proof of Claim \ref{theta}.}
The proof is done 
in a similar way in \cite{Takeda} and a generalized fashion of the manner in 
\cite[section three]{Rourke}: 
The fundamental group of the virtual reef knot $R$ is $\Z$. 
However the fundamental group of the mirror image of $R$ is 
non-trivial. 
The fundamental group of the mirror image of $R$ is the lower fundamental group of $\mathcal O(R)$ and 
its non-triviality demonstrates the non-triviality of 
$\mathcal O(R)$ 
as a virtual 2-knot.
\qed\\

Theorem \ref{Maine} is one of our main results. \\

\begin{thm}\label{Maine}  
There is a virtual 2-knot $K$ 
with the following conditions. 

\smallbreak\noindent$(1)$ 
The virtual 2-knot $K$ 
is not virtually equivalent to 
the standard surface. 

\smallbreak\noindent$(2)$ 
The virtual 2-knot $K$ 
is $\mathcal E$-equivalent to 
the standard surface. 
\end{thm}

\noindent{\bf Proof of Theorem \ref{Maine}.}  
Let $K$ be the virtual 2-knot $\mathcal O(R)$ in Claim \ref{theta}. 
Claim \ref{theta} implies Theorem \ref{Maine}.(1). 
Let $T$ denote the standard torus. 
By Note \ref{tau},  $\mathcal E(R)=T$. 
Proposition \ref{kakan} implies 
$\mathcal E(K)
=\mathcal E(\mathcal O(R))
=\mathcal O(\mathcal E(R))
=\mathcal O(T)$, 
where $=$ denotes the ambient isotopy  
of submanifolds. 
$\mathcal O(T)$ and $\mathcal E(T)$ are standardly embedded $T^3$ in $\R^5$ by the definition of them.
Hence we have Theorem \ref{Maine}.(2). 
Therefore Theorem \ref{Maine} holds.  \qed\\

We  ask  questions.
\begin{que}\label{France}
(1) Do we have the following? 
Let $\Sigma_g$ be 
a closed oriented genus $g$ surface.  
Let $Q$ (respectively, $Q'$) be  
a virtual surface-knot 
made from $\Sigma_g$. 
If $Q$ and $Q'$ have the group $\Z$, 
then the submanifolds, $\mathcal E(Q)$ and $\mathcal E(Q')$, 
of $\R^5$ are isotopic. 

\smallbreak 
\noindent
(2) 
Is a virtual 1- (respectively, 2-) knot $K$ welded equivalent to the trivial 1-knot 
if $K$ has the group $\Z$? 
\end{que}

\bigbreak
\section{The fibrewise equivalence}\label{New Mexico}
\subsection{
The fibrewise equivalence
is equal to the rotational welded equivalence, 
and is different from the welded equivalence of virtual 1-knots
}\label{sub1}\hskip20mm\\%

We research relations among 
the fiberwise equivalence of virtual 1-knots, 
the welded equivalence of them, 
and 
the rotational welded equivalence of them. 
We mentioned it in the last few paragraphs of \S\ref{i3}. 
See \cite{Rourke, Satoh} for the definition of the welded equivalence, 
and \cite{Kauffman, Kauffmanrw, J} for that of the rotational welded equivalence, 
as we also mentioned them in the last few paragraphs of \S\ref{i3}.\\

We first introduce the definition of the fiberwise equivalence of virtual 1-knots. 
For our purpose (to prove Theorems \ref{smooth} and \ref{Montgomery}), 
we will modify the definition a few times as below. 

\begin{figure}
\includegraphics[width=100mm]{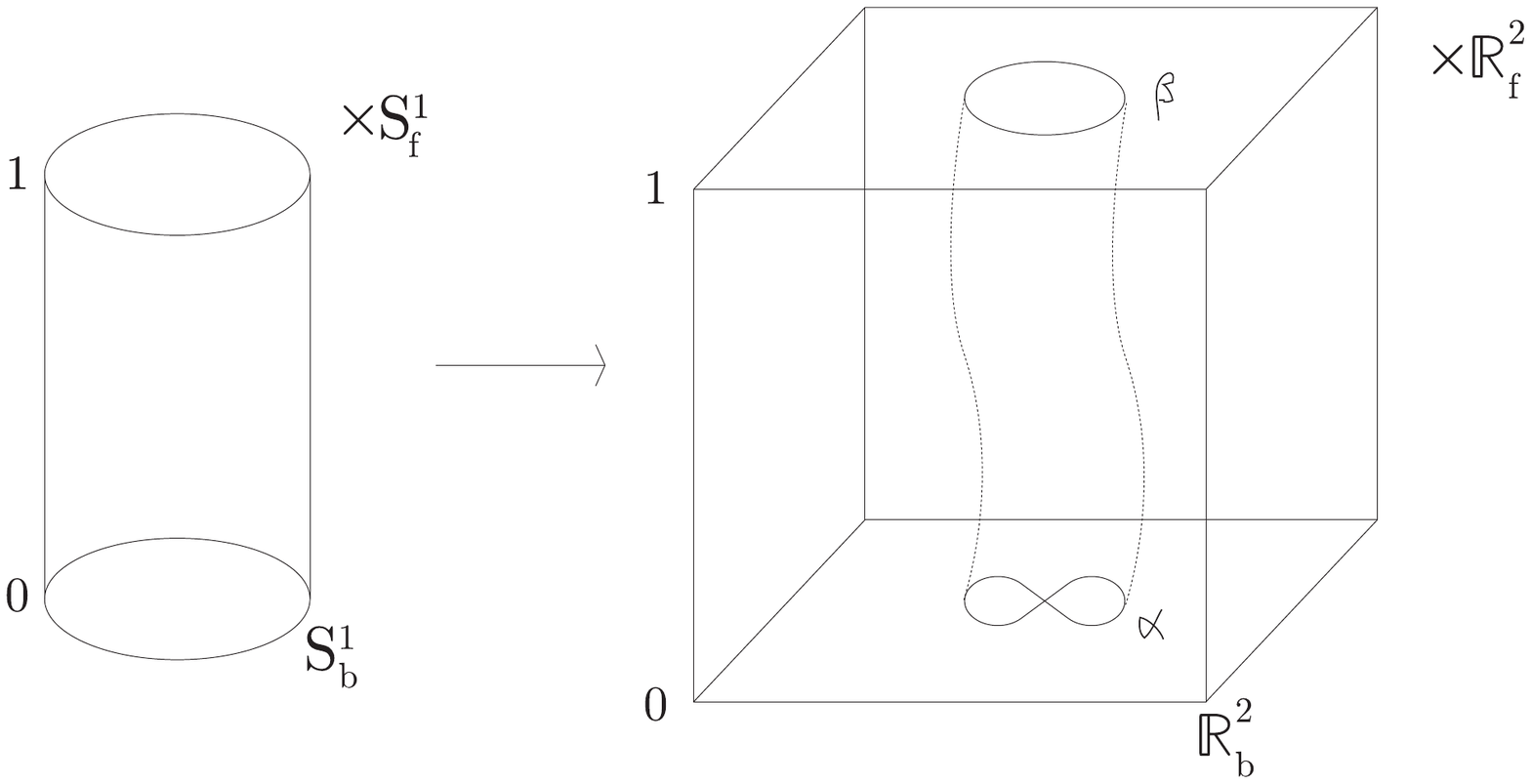}  
\vskip-40mm
\caption{{\bf Fiberwise isotopy}\label{Fib}}   
\end{figure}

\begin{defn}\label{Nevada}  
Let $\alpha$ and $\beta$ be virtual 1-knot diagrams. 
We say that 
$\alpha$ and $\beta$ are {\it fiberwise equivalent} 
if 
 Rourke's description of $\mathcal S(\alpha)$ 
and 
that of $\mathcal S(\beta)$ 
are `fiberwise isotopic'.   
In other words, this means that 
$\alpha$ and $\beta$ satisfy the following conditions. 
There is an  
 embedding map 
$$g:S^1_b\x[0,1]\x S^1_f\hookrightarrow\R^2_b\x[0,1]\x\R^2_f$$  
with the following properties. 
See Figure \ref{Fib}.

\smallbreak\noindent
(1)  
For any fixed $t\in[0,1]$, $g(S_b^1\x\{t\}\x S^1_f)\subset\R^2_b\x\{t\}\x\R^2_f$. 

\smallbreak\noindent
(2) 
For any fixed $p\in S^1_b$ and any fixed $t\in[0,1]$, 
$g(\{p\}\x\{t\}\x S^1_f)$ is contained in 
the same fiber $\{q\}\x\R^2_f$ for a point 
$q\in\R^2_b\x[0,1]$. 

  
\smallbreak\noindent (3) 
Let $\pi:\R^2_b\x[0,1]\x\R^2_f\to\R^2_b\x[0,1]$. 
$(\pi\circ g)(S_b^1\x\{0\}\x S^1_f)$  
(respectively, \\$(\pi\circ g)(S_b^1\x\{1\}\x S^1_f)$) 
$\subset \R^2_b\x\{0\}$ (respectively, $\R^2_b\x\{1\}$)  
is the diagram $\alpha$ (respectively, $\beta$) without information whether 
each crossing point is a classical one or a virtual one. 
This information is given by the fiber-circles over each crossing point as 
in Theorem \ref{Montana} and Definition \ref{Nebraska}.   
$\pi\circ g$ meets $R^2_b\x\{0,1\}$ transversely. 

\bigbreak
In knot theory we usually use an `ambient' isotopy in order to define the equivalence relation of knots as below. 
We impose the following condition (4). 
(See \cite[sections 1.1 and 1.2]{BZ}  for an explanation on this fact 
in the 1-dimensional classical knot case.)  

\smallbreak\noindent (4) 
Let $g_t$ denote 
$$g|_{S^1_b\x\{t\}\x S^1_f}: S^1_b\x\{t\}\x S^1_f\hookrightarrow\R^2_b\x\{t\}\x\R^2_f$$  
for $0\leqq t\leqq1$. 
There is an an isotopy 
$$H_t:\R^2_b\x\{t\}\x\R^2_f\to\R^2_b\x\{t\}\x\R^2_f (0\leqq t\leqq1)$$
such that 
$H_0$ is the identity map 
and such that 
$g_t=H_t\circ g_0$ for any $t\in[0,1]$. 

We call $g$ a {\it special isotopy} between $\alpha$ and $\beta$.  
\end{defn}
\bigbreak

\begin{defn}\label{anko}  
Take $g$ in Definition \ref{Nevada}. 
If $g'$ is obtained by moving $g$ by an ambient isotopy map 
$G_t$, where $0\leq t\leq1$, $g_0=g$, and $g_1=g'$,  
keeping the conditions $(1)$-$(4)$ of Definition \ref{Nevada},  
then we say that $g'$ is {\it level preserving, fiberwise isotopic} or {\it special isotopic} to $g$, 
or 
that we {\it perturb $g$ in the special way} to obtain $g'$.  
We write $g\sim g'$. 
$G_t$ is called a 
{\it level preserving, fiberwise isotopy} 
or
{\it special isotopy} between $g$ and $g'$. 
\end{defn}
\bigbreak

\begin{note}\label{kakanzu}
The following holds. 
Let $\rho:S^1_b\x[0,1]\x S^1_f\to S^1_b\x[0,1]$ be the natural projection.   
Then there is a (not necessarily smooth) continuous map  \\
$\underline{g}:S^1_b\x[0,1]\to \R^2\x[0,1]$  
such that $\pi\circ g=\underline{g}\circ\rho$. 
That is, there is  the following commutative diagram.
$$
\begin{matrix}
S^1_b\x[0,1]\x S^1_f&\stackrel{g}\to&\R^2_b\x[0,1]\x \R^2_f\\
\downarrow_\rho&\circlearrowright &\downarrow_\pi \\
S^1_b\x[0,1]&\stackrel{\underline{g}}\to&\R^2_b\x[0,1]
\end{matrix}
$$
\end{note}

\begin{defn}\label{neba}
Under the above condition, we say that $\underline{g}$ is {\it covered} by $g$. 
\end{defn}

\bigbreak 
The following theorem is one of our main results. 

\begin{thm}\label{smooth}
Two virtual 1-knot diagrams $\alpha$ and $\beta$ are smooth fiberwise equivalent if and only if  
$\alpha$ and $\beta$ are smooth rotational welded equivalent.  
\end{thm}

\noindent{\bf Note.}  See Note \ref{xmikan}. \bigbreak

\noindent
{\bf Proof of Theorem \ref{smooth}.} 
The `if' part is easy. 

We prove the `only if' part.
\\

\noindent{\bf Strategy.} 
See (I) and (II) below. 
We want to prove  (I)$\Leftrightarrow$(II). 
It is easy to prove  (II)$\Rightarrow$(I). 
We will prove (I) $\Rightarrow$ (II) as follows. 

\bigbreak
\noindent(I) Smooth virtual 1-knot diagrams $\alpha$ and $\beta$ are smooth fiberwise equivalent. 

\bigbreak
\noindent(II)  Smooth virtual 1-knot diagrams $\alpha$ and $\beta$ are smooth rotational welded equivalent. 

\bigbreak

See (1) below.
In Claim \ref{xbeef}  
we will prove (I)$\Rightarrow$(1). 

\bigbreak
\noindent(1) There is a PL virtual 1-knot diagram $\alpha'$ (respectively, $\beta'$) 
which is piecewise smooth isotopic to  $\alpha$ (respectively, $\beta$) 
such that $\alpha'$ and $\beta'$ are PL fiberwise equivalent. 
\bigbreak

See (2) below. 
In Theorem \ref{fwrw}, we will prove  (1)$\Rightarrow$(2). 
It will be proved in the text 
which starts from Proposition \ref{polygon}, 
and ends in Claim \ref{takusan}. 

\bigbreak
\noindent(2) $\alpha'$ and $\beta'$ are PL rotational welded equivalent.


\bigbreak
In Lemma \ref{PLtosmooth}, we will prove (2) 
$\Rightarrow$ (II).  Thus we will finish the proof  of  (I) $\Rightarrow$ (II).

\bigbreak
\bigbreak
%
Assume that smooth virtual 1-knot diagrams $\alpha$ and $\beta$ are smooth fiberwise equivalent. 
We do not know whether or not there are two 
special isotopies $g$ and $g'$ between $\alpha$ and $\beta$ 
%
 with the following properties. 
$g$ and $g'$ are not smooth 
special isotopic 
 but piecewise smooth special isotopic 
Although we do not answer this question,   
we accomplish the proof of (I)$\Leftrightarrow$(II).

\bigbreak
Take $g$ 
in Definition \ref{Nevada}.   
We do not know
whether there is a smooth $g'$ with $g'\sim g$ 
with the following properties: 
There is a finite simplicial structure on 
Im $\pi\circ g'$ 
which restricts to 
a finite simplicial structure on the singular subset 
of Im  $\pi\circ g'$. 
One reason is as follows. Im $\pi\circ g$ may be the projection of a wild embedding
for a smooth $g$ even if $g$ is not a wild embedding map. 
Although we do not answer this question, we accomplish the proof of (I)$\Leftrightarrow$(II).

\begin{defn}\label{PLNevada}
Consider the conditions of Definition \ref{Nevada} in the PL category. 
The equivalence relation is called {\it PL fiberwise equivalence}.
\end{defn}

\begin{cla}\label{xbeef}
If virtual 1-knot diagrams $\alpha$ and $\beta$ are smooth fiberwise equivalence, 
then  $\alpha$ and $\beta$ are PL fiberwise equivalence. 
\end{cla}

\noindent{\bf Proof of Claim \ref{xbeef}.}
It is enough to prove that 
the map $g$ in Definition \ref{Nevada} is approximated 
by a fiberwise level-preserving PL embedding map. We prove it below. 

Regard $S^1_b$ as $[0,1]/\sim$, where $0\sim1$. 
Regard $S^1_f$ as $[0,1]/\sim$, where $0\sim1$. 
Hence we can regard 
$S^1_b\x[0,1]\x S^1_f$  
as the one made from $[0,1]\x[0,1]\x[0,1]$ by these equivalence relations. 

Let $n$ be any positive integer. 
Take points 
$(\frac{i}{2n},\frac{j}{2n},\frac{k}{2n})\in[0,1]\x[0,1]\x[0,1]$, 
where $i$ (respectively, $j$, $k$) is any integer with the condition 
$0\leqq$$i$ (respectively, $j$, $k$) $\leqq 2n$. 

Let $l$ be any integer with the condition $0\leqq 2l$(respecctively, $2l+2)\leqq 2n$. 
Take any cube $C$ whose vertices are 
$(\frac{\alpha}{2n},\frac{\beta}{2n},\frac{\gamma}{2n})$, 
where $\alpha$ (respectively, $\beta$, $\gamma$) is any integer in $\{2l,2l+2\}$. 

Take a simplicial division on $S^1_b\x[0,1]\x S^1_f$ as follows. 
  
\smallbreak\noindent(1) 
0-simplices are all $(\frac{i}{2n},\frac{j}{2n},\frac{k}{2n})\in[0,1]\x[0,1]\x[0,1]$ as above.  

\begin{figure}
\includegraphics[width=120mm]{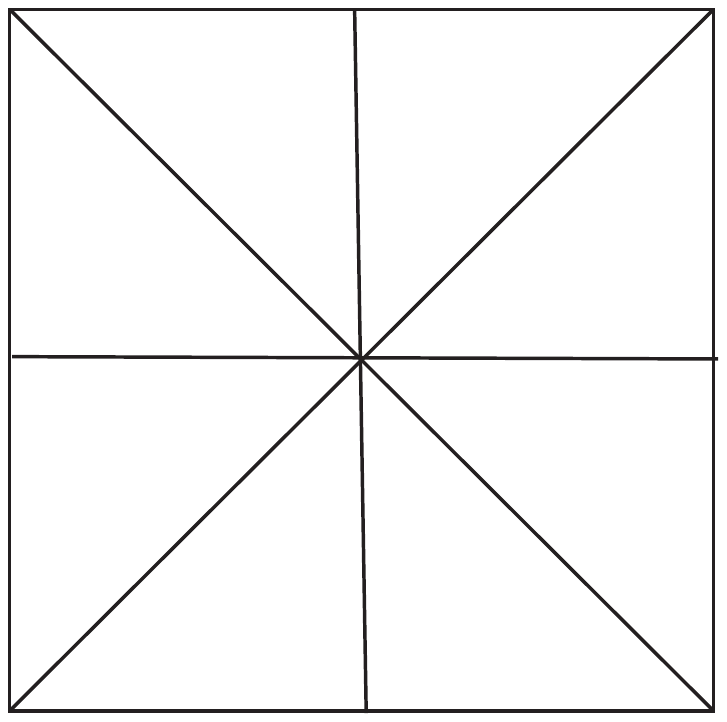}  
\vskip-30mm
\caption{{\bf 16 1-simplices on a site of $C$}\label{menkirikata}}   
\end{figure}

\smallbreak\noindent(2) 
1-simplices are defined as follows. 
Take any cube $C$. 
Note that each of six sites includes nine 0-simlices, 
that the sum of six sites includes 26 0-simlices, 
and 
that $C$ includes 27 0-simlices. 
Take the 0-simplex $P$ in $C$ which is not included in any site. 
Take any segment whose boundary is $P$ and one of the other 28 0-simplices.  
Take 16 segments in each site of $C$ as drawn in Figure \ref{menkirikata}. 
1-simplices are these two kinds of segment.

\smallbreak\noindent(3) 
The set of 1-simplices defines 2-simplices naturally. 

\smallbreak\noindent(4)
The set of 2-simplices defines 3-simplices naturally. 

\bigbreak

Let $n$ be sufficiently large. 
Take the image of all 0-simplices by $g$ in  $\R^2_b\x[0,1]\x\R^2_f$. 
They determine a fiberwise level-preserving PL embedding map of $S^1_b\x[0,1]\x S^1_f$ naturally.  
{\it Reason.} Im $g$ is a smooth regular submanifold. 
Hence it has a tubular neighborhood. 

This completes the proof of Claim \ref{xbeef}. \qed

\begin{note}\label{bango}
If $C=$(Im $g$)$\cap$(a fiber $\R^2_f$) is PL homeomorphic to a circle,  
then $C$ is a polygon. However the number of the vertices of $C$ depends on fibers. 
\end{note}

\begin{note}\label{haruwa}
From here to the end of the proof of Theorem \ref{fwrw}, we work in the PL category 
unless we indicate otherwise. 
After that, we will go back to the smooth category. 

When we move a map by isotopy, we take a PL subdivision if necessary. 
\end{note}

Claim \ref{xbeef} implies the following. 

\begin{pr}\label{polygon}
%
$g$ in Definition \ref{PLNevada} satisfies the condition that a finite simplicial structure on Im $\pi\circ g$ 
which restricts to a finite simplicial structure on the singular subset of Im $\pi\circ g$. 
\end{pr}

\begin{figure}
\bigbreak  
\includegraphics[width=150mm]{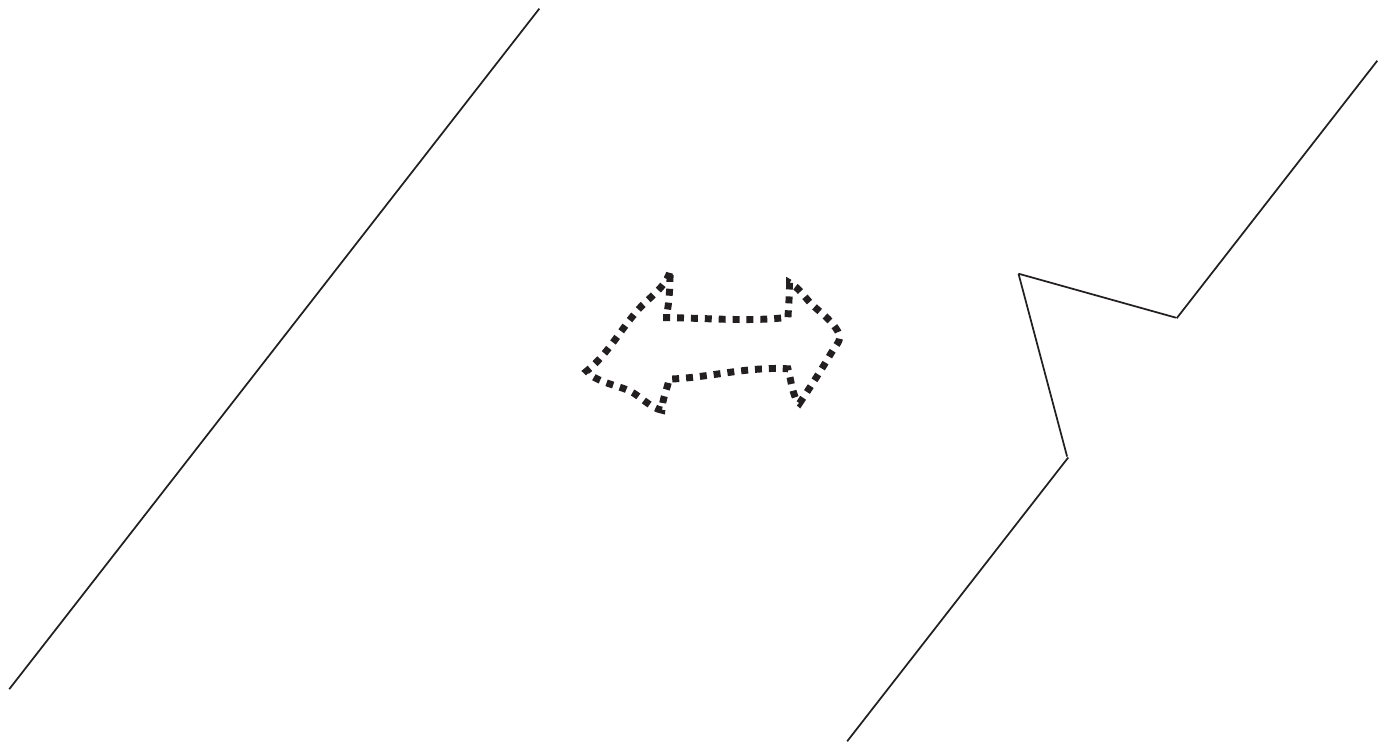}  
\vskip-40mm
\caption{{\bf The $\Delta^1$-move.}\label{zudelta}}   
\end{figure}

We call the operation drawn in Figure \ref{zudelta},  
the $\Delta^1$-move of virtual 1-knot diagrams. 
Note that we do not draw the other part of this diagram. The other part may intersect the part drawn in  Figure \ref{zudelta}. 
By Proposition \ref{polygon},  $\alpha, \beta$ in Definition \ref{PLNevada}  
have the following properties: 

\begin{cla}\label{Delta} 
$\alpha$ $($respectively, $\beta$$)$ is obtained from $\beta$ $($respectively, $\alpha$$)$ by 
a finite step of $\Delta^1$-moves. 
\end{cla}

\begin{defn}\label{PLyugentsuki}
Add the following condition to Definition \ref{PLNevada} without changing the other parts. 
(Note we work in the PL category.)

\smallbreak\noindent$(\ref{PLyugentsuki}.1)$  In each fiber  $\R^2_f$,  there are a finite number of circles. 
$($That is, $<\infty.)$
\end{defn}

\noindent{\bf Note.} 
See Note \ref{xudon}. 
Recall Note \ref{bango}.  
\\

Indeed, the following holds. 

\begin{thm}\label{herasu}
Definitions $\ref{PLNevada}$ and $\ref{PLyugentsuki}$ are equivalent.  
\end{thm}

\begin{figure}
\bigbreak  
 \includegraphics[width=145mm]{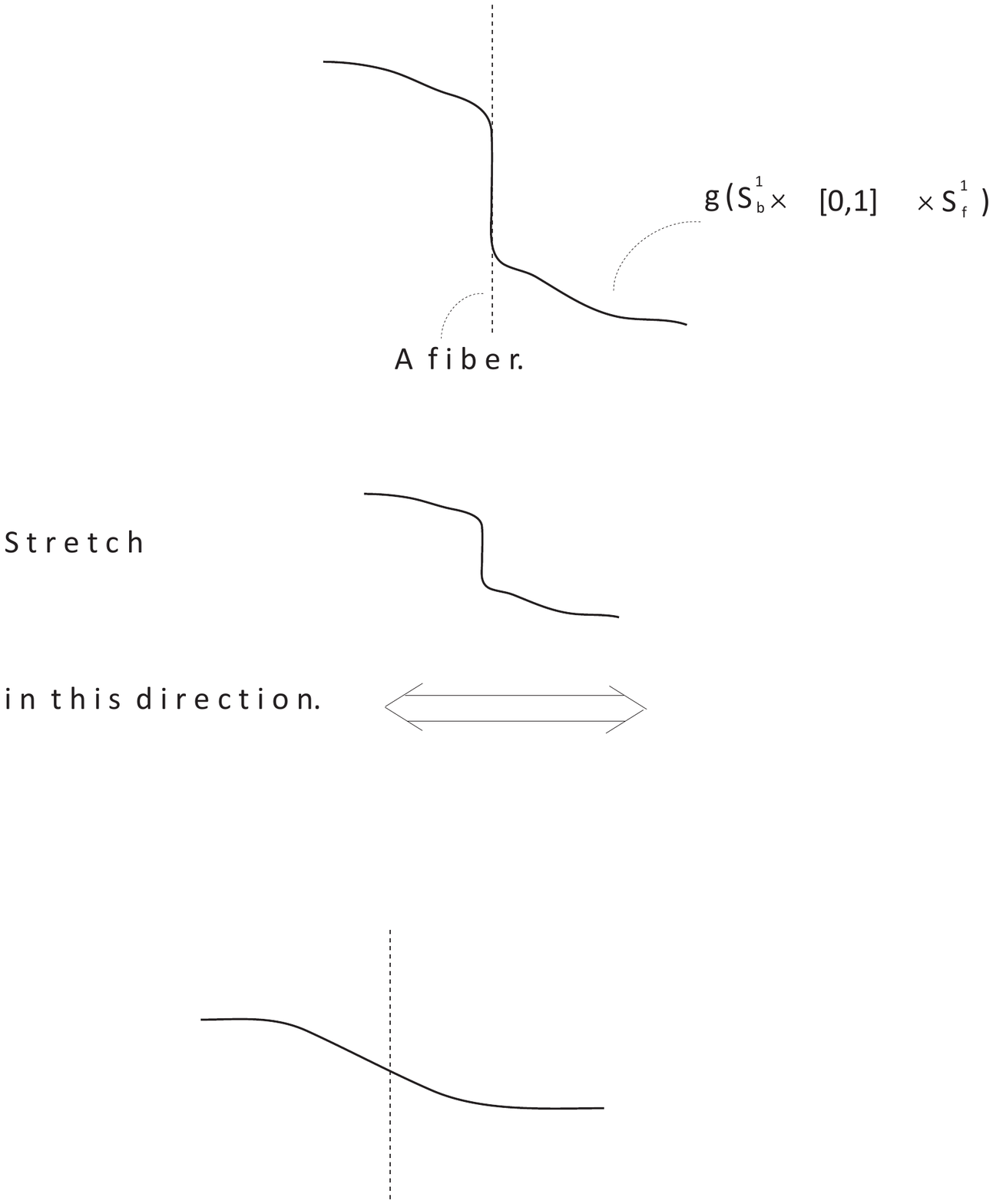}  
\vskip-37mm
\caption{{\bf 
The idea of how we stretch $g(S^1_b\x[0,1]\x S^1_f)\cap N$ 
}\label{hipparu}}   
\end{figure}

\noindent{\bf Proof of Theorem \ref{herasu}.}  
It is trivial that if $g$ satisfies Definition \ref{PLyugentsuki}, then $g$ satisfies Definition \ref{PLNevada}. 
We prove that if $g$ satisfies  Definition \ref{PLNevada},  
then we can perturb $g$ in the special way so that $g$ satisfies Definition \ref{PLyugentsuki}. 
Suppose that $g$ satisfies Definition \ref{PLNevada}. 
Let $q\in\R^2_b$ and $t\in[0,1]$. 
Since Im$g$ is a compact PL regular submanifold, 
Im$g\cap(\{q\}\x\{t\}\x\R_f)$ is 
a disjoint union of a finite number of circles and a finite number of annuli.  
Note that the union of them is a regular submanifold of $\{q\}\x\{t\}\x\R^2_f$.  
Take a tubular neighborhood $N$ of each annulus 
in $\R^2_b\x[0,1]\x\R^2_f$, to be small enough. 
Stretch each annulus into the direction perpendicular to $\{x\}\x\{t\}\x\R^2_f$. 
Then we can obtain a new $g$ which satisfies  Definition \ref{PLyugentsuki}. 
The idea of how we stretch is drawn in Figure \ref{hipparu}
Note that 
Figures \ref{hipparu} draws `figures in PL category' although the figures are smoothened. 
\qed

\bigbreak
A point 
$p\in$Im$\pi\circ g$
$=(\pi\circ g)(S^1_b\x[0,1]\x S^1_f)$
$=(\underline{g}\circ\rho)(S^1_b\x[0,1]\x S^1_f)$ 
$=\underline{g}(S^1_b\x[0,1])$  
is called a {\it multiple point} or  {\it $n$-tuple point} 
if $\underline{g}^{-1}(p)\in S^1_b\x[0,1]$  
consists of $n$ points ($n\geqq2$). 
(Note that in Definition \ref{PLyugentsuki}, $n<\infty$.)
A point $p\in$Im$\pi\circ g$ 
is called a {\it single point} 
if $\underline{g}^{-1}(p)$ consists of a single points. 
The {\it singular point set} of $p\in$Im$\pi\circ g$ consists of branch points and multiple points.  \\

Note the following facts. 
Take $g$ in Definition \ref{PLyugentsuki}, and $\underline{g}$ which is covered by $g$. 
Recall that `cover' is defined in Definition \ref{neba}.  
Suppose that $\underline{g}$ is a generic map. 
Note Im $\underline{g}$. 
We can define whether  each double point is classical or virtual 
by using the information of the fiber-circles over each point as in Theorem \ref{Montana}, 
Definition \ref{Nebraska}, Note \ref{kaiga}, and Definition \ref{suiri}.   
There is a case where a classical (respectively, virtual) double point appears.  
The information of fiber-circles over each branch point determines that the branch point is classical. 
{\it Reason.} By Theorem \ref{Rmuri}, 
there are no virtual branch point.  \\

\begin{figure}
\bigbreak  
 \includegraphics[width=140mm]{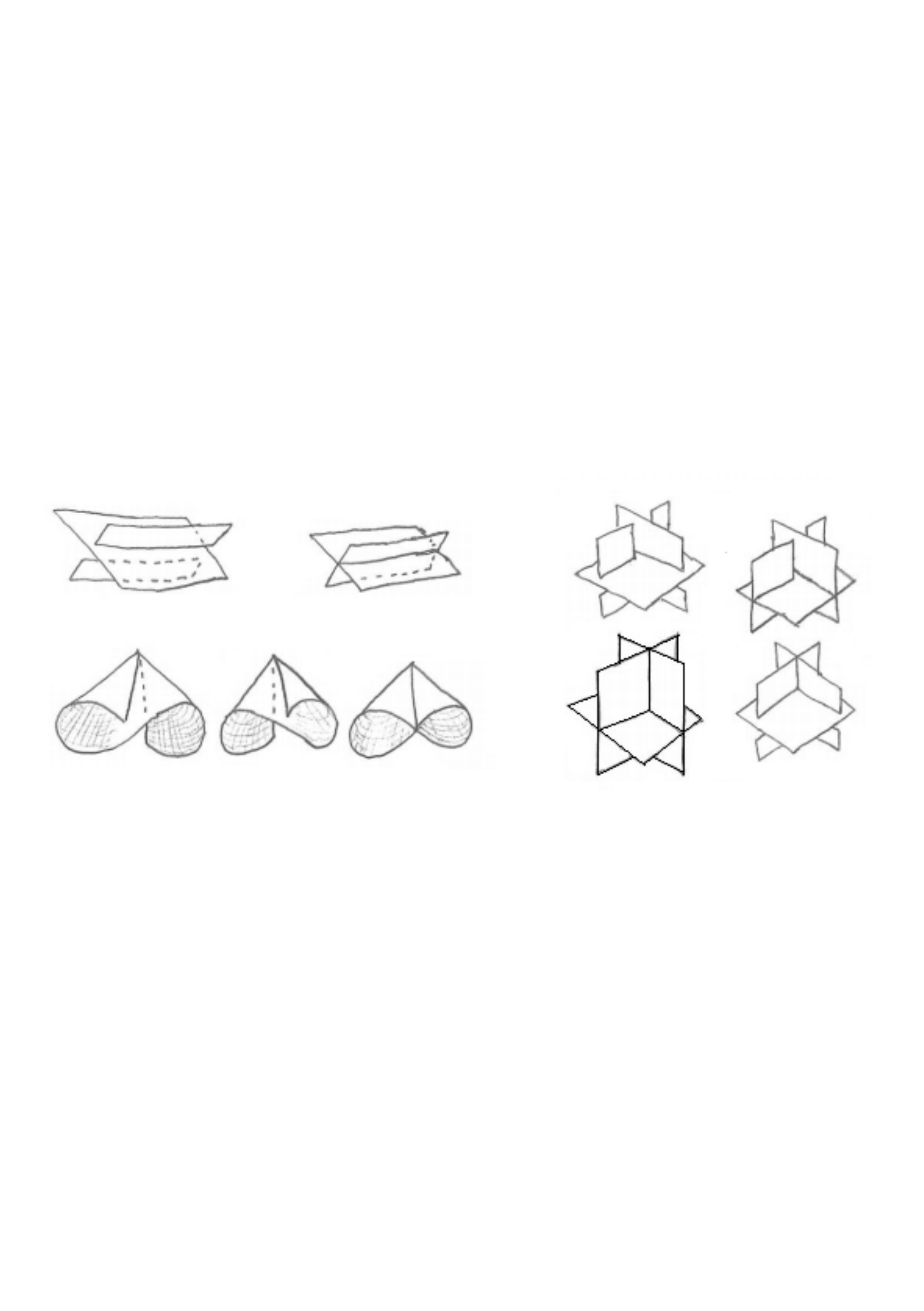}  
\vskip-30mm
\caption{{\bf The singular point sets of welded 2-knots}
\label{JW1}}   
\bigbreak  
\end{figure}

Note each triple point. 
There are three circles in the fiber over each triple point. 
There are four cases how three circles are put in the fiber.  
See Notes \ref{umeboshi} and \ref{faso}, Definnition \ref{JW}, and Figure \ref{JW1}.    
There is a case where each of the four occurs.

\begin{note}\label{umeboshi}
$(\pi\circ g)(S^1_b\x[0,1]\x S^1_f)$  in $\R^2_b\x[0,1]$  
is a welded 2-knot with a fixed boundary  in general, 
and 
is not a virtual 2-knot with a fixed boundary  in general.
See \cite[sections 3.5-3.7]{J} for their definitions 
and their difference.
In the welded 2-knot case we also use the terms, `fiber-circle' and `Rourke-fibration'.  
See Note \ref{faso}.
\end{note}

Here we cite the definition of welded 2-knots from \cite{J}. 



%


Recall that a 2-knot diagram is (the image of) a generic map of a surface in 3-space, with
classical and virtual crossing data along the double-point arcs. 
Also recall that 2-knot diagrams may be transformed by Roseman moves, which preserve the
crossing data locally.

\begin{defn}\label{JW}{\bf (\cite[section 3.6]{J}.)}  
If all triple points of a 2-knot diagram are of the four types shown in Figure \ref{JW1}, 
the diagram is called a Welded 2-knot diagram. If a pair of Welded 2-knot
diagrams are related by a series of Roseman moves, with only Welded diagrams
appearing throughout the process, then the diagrams are Welded equivalent and
belong to the same Welded 2-knot type.
\end{defn}

\noindent{\bf Note.} 
The above definition makes sense both in the smooth and the PL category.  
The readers need not be familiar with Roseman moves in order to read this paper.

\begin{figure} 
\includegraphics[width=100mm]{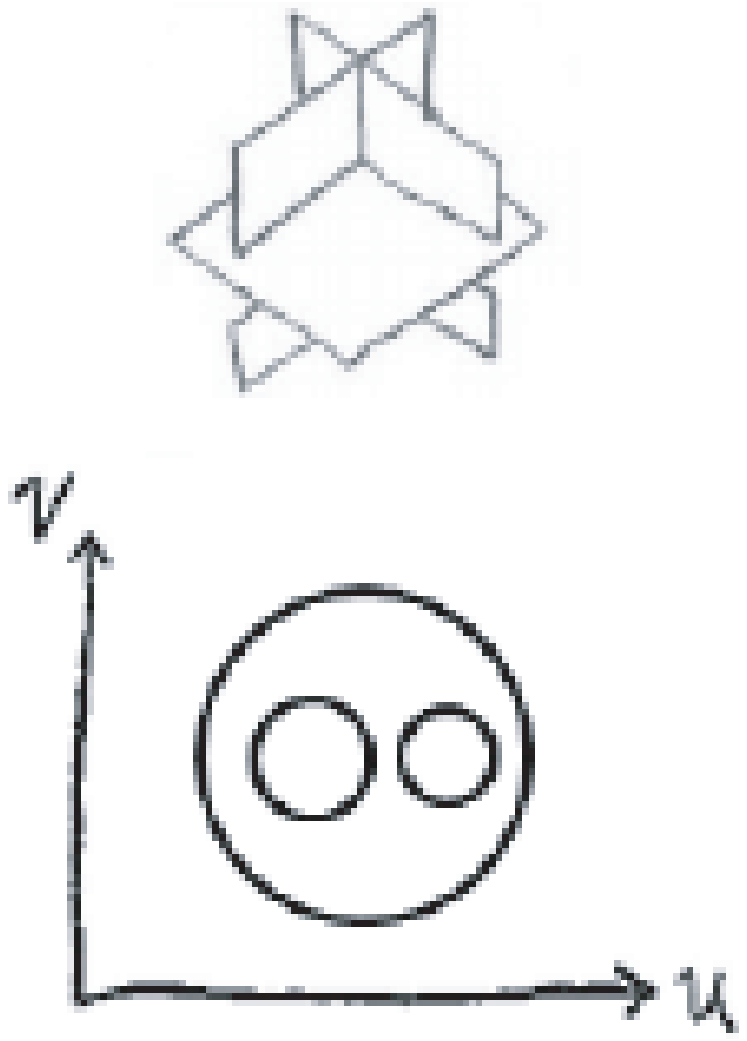}  
\vskip-20mm
\caption{{\bf A new type of a nest of circles.}\label{rashi}}   
\end{figure}

\begin{note}\label{faso}
When we consider circles in fiber $\R^2$ as in Note \ref{umeboshi}, 
there is a new type drawn in Figure \ref{rashi},  
which is not in Figure \ref{doremi}. 
\end{note} 

\begin{note}\label{oshii} 
By Proposition \ref{polygon}, $\underline{g}$ in Definition \ref{PLyugentsuki} 
satisfies the conditions (I)-(III) below, but $\underline{g}$ is not generic. 

\smallbreak\noindent (I)   
$\underline{g}:S^1_b\x[0,1]\to\R^2_b\x[0,1]$
is a continuous map such that 
$\underline{g}(S^1_b\x\{t\})\subset\R^2_b\x\{t\}$ for any $t\in[0,1]$. 

\smallbreak\noindent (II)   
Let $t\in[0,1]$. 
There are closed intervals,  
$I_1,..., I_\mu$ ($\mu\in\N$), 
such that 
$S^1_b\x\{t\}\\=I_1\cup...\cup I_\mu$
and such that  
$\underline{g}|_{I_i}$ is a PL 
embedding for each $i$. 

\smallbreak\noindent (III)   
There are closed 2-discs,  
$D^2_1,..., D^2_\nu$ ($\nu\in\N$), 
such that 
$S^1_b\x[0,1]=D^2_1\cup...\cup D^2_\nu$
and such that  
$\underline{g}|_{D^2_i}$ is a PL 
embedding for each $i$. 
\end{note}

\begin{defn}\label{wasabi} 
If a map 
$\underline{g}:S^1_b\x[0,1]\to\R^2_b\x[0,1]$ satisfies the conditions (I)-(III) in Note \ref{oshii}, 
then $\underline{g}$ is said to be {\it level preserving}. 
If $\underline{g}'$ is obtained by moving $\underline{g}$ by a homotopy 
$\underline{G_t}$, where $0\leq t\leq1$, $\underline{G_0}=g$ and $\underline{G_1}=g'$,  
keeping the conditions (I)-(III) of in Note \ref{oshii},  
then we say that $\underline{g}'$ is {\it level preserving homotopic} to $\underline{g}$ 
or 
that we {\it perturb $\underline{g}$ in the special way} and obtain $\underline{g}'$. 
We write $\underline{g}\sim\underline{g}'$. 
$\underline{G_t}$ is called a {\it level preserving homotopy} or a {\it special homotopy}.  
Let $g: S^1_b\x[0,1]\x S^1_f\to\R^2_b\x[0,1]\x \R^2_f$ 
be a map in Definition \ref{PLyugentsuki}. 
Take a special homotopy 
$\underline{G_t}$ of $\underline{g}$, and 
a special isotopy 
$G_t$ of $g$ 
where $0\leqq t\leqq1.$
If $\underline{G_t}$ 
is covered by $G_t$  
for any element $t$ in $\{t|0\leqq t\leqq1\}$, 
then we say that 
$\underline{G_t}$ is {\it covered} by $G_t$. 
\end{defn}

\begin{defn}\label{gene}
Add the following condition to 
Definition \ref{PLyugentsuki} without changing the other parts. 
\noindent$(\ref{gene}.1)$ 
We can perturb $g$ in  Definition \ref{PLyugentsuki} in the special way 
so that $g$ covers a PL level preserving, generic map $S^1_b\x[0,1]\to\R^2_b\x[0,1]$. 
\end{defn}

We prove the following theorem. 

\begin{thm}\label{Wyoming}  
Definition $\ref{gene}$ is equivalent to Definition $\ref{PLyugentsuki}$  
$($and, by Theorem $\ref{herasu},$  is equivalent to Definition $\ref{PLNevada}.)$  
\end{thm} 

\begin{note}\label{shio} 
Even if we perturb $\underline{g}: S^1_b\x[0,1]\to\R^2_b\x[0,1]$, which is covered by $g$, 
in the special way 
by a special homotopy $\underline{G_t}$,  
$\underline{G_t}$ is not covered by a special isotopy $G_t$ of $g$ 
in general.  
We must make 
$\underline{G_t}$ 
under the condition that 
$\underline{G_t}$ is covered by $G_t$. 
\end{note}

\noindent
{\bf Proof of Theorem \ref{Wyoming}.} 
It is trivial that 
if $g$ satisfies Definition \ref{gene}, then $g$ satisfies Definition \ref{PLyugentsuki}. 
We prove the following. 

\begin{cla}\label{kamen}
If $g$ satisfies Definition \ref{PLyugentsuki}, 
we can perturb $g$ so that $g$ satisfies Definition \ref{gene}. 
\end{cla} 

\noindent{\bf Note}. 
Recall that $\pi\circ g$ does not cover a generic  map $\underline{g}$ in general. 
\\

\noindent
{\bf Proof of Claim \ref{kamen}.} 

\noindent
{\bf The first step.} 
Recall that by Definition \ref{PLyugentsuki}, for each $t\in[0,1]$, (Im$(\pi\circ g)$)$\cap(\R^2\x\{t\})$ is an immersed circle. We prove the following. 

\begin{cla}\label{Oregon}  
We can perturb $g$ 
in the special way 
so that the singular point set of \\
$(${\rm Im}$(\pi\circ g))\cap(\R^2\x\{t\})$ 
is a finite number of points 
except for a finite number of levels  $t\in[0,1]$. 
In other words, we can do so that for only a finite number of levels $t\in[0,1]$, 
the singular point set of $($Im$(\pi\circ g))\cap(\R^2\x\{t\})$ includes a finite number of segments. 
\end{cla}

\noindent
{\bf Proof of Claim \ref{Oregon}.}
Let $I$ denote the interior of a 1-simplex which is in a simplicial complex structure of 
the singular point set of (Im$(\pi\circ g)$)$\cap(\R^2\x\{\gamma\})$
 for a real number $\gamma\in[0,1]$. 
Assume that $I$ consists of multiple PL points.

\begin{figure}
\vskip-30mm
\includegraphics[width=170mm]{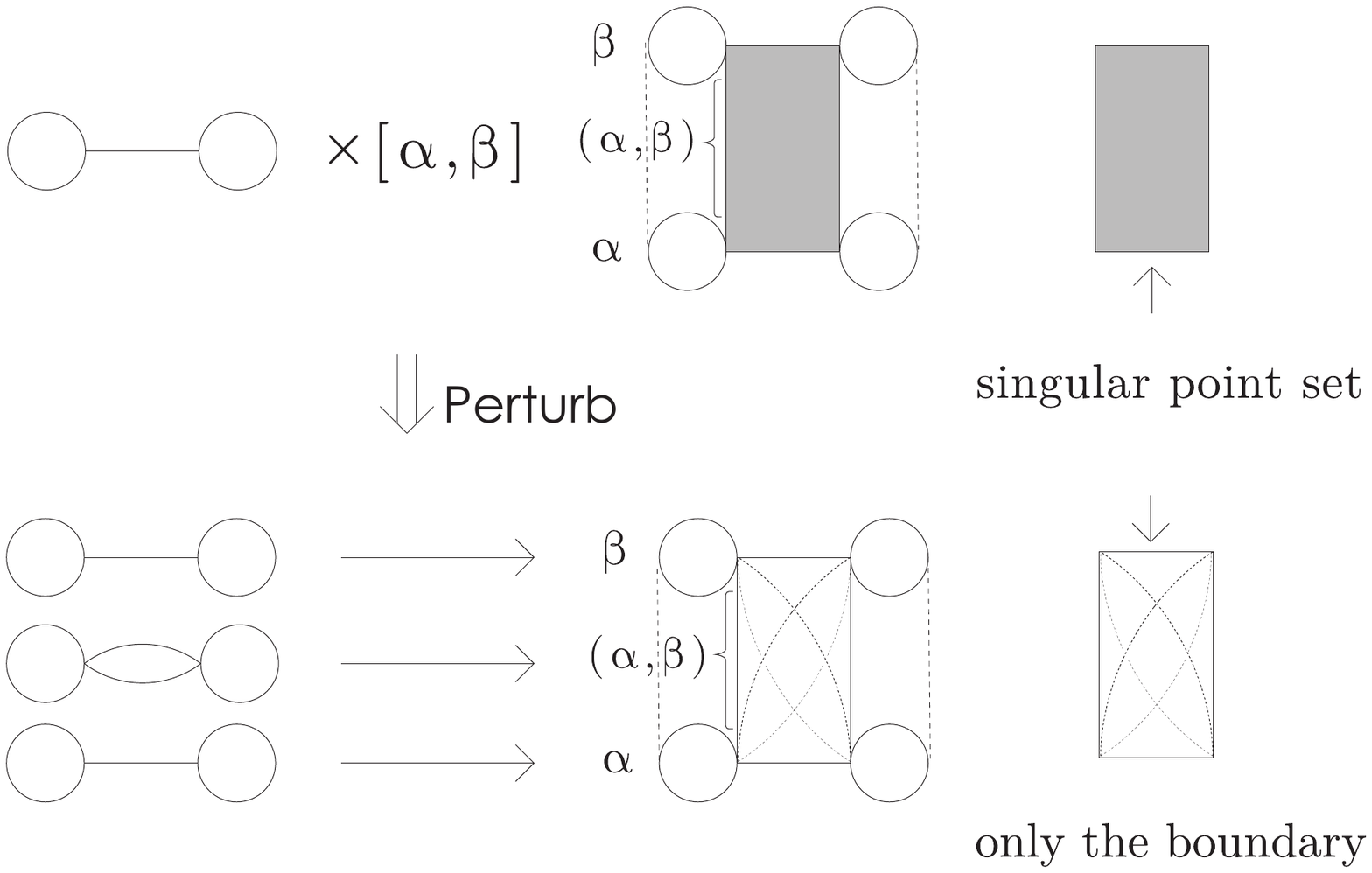}  
\vskip-50mm
\caption{{\bf An example of fiberwise isotopy. }\label{step1}}   
\end{figure}

Suppose that there are real numbers 
 $\alpha,\beta\in[0,1]$ with the following properties: \\
 $\alpha<\gamma<\beta$. 
For $\alpha<t<\beta$, 
$h_t:R^2_b\x\{\gamma\}\to\R^2_b\x(\alpha,\beta)$ is an isotopy 
($t$ runs in $(\alpha,\beta)$) 
 such that 
$h_t(({\rm Im}(\pi\circ g))\cap(\R^2\x\{\gamma\}))=$ 
$({\rm Im}(\pi\circ g))\cap(\R^2\x\{t \})$
for all $t\in(\alpha, \beta)$,  
%
%
%
and such that 
$h_t$ preserves   
the information of fiber-circles over  
 two immersed circles 
which are put in the both sides of $=$.
(Here, 
the information of fiber-circles means what we define in  
Theorem \ref{Montana},  Definition \ref{Nebraska}, Note \ref{kaiga}, and Definition \ref{suiri}.)
Note that $(\underline{g})^{-1}(I)$ 
is a disjoint union of $n$ open segments $I_1,...,I_n$ in $S^1_b\x[0,1]$. 
Note that $\underset{\alpha<t<\beta}{\cup} h_t(I)$ 
consists of $n$-tuple points, 
is an open set, and 
is a discrete submanifold of $\R^2_b\x[0,1]$.  
We can perturb $g$ in the special way 
so that $\underset{\alpha<t<\beta}{\cup} h_t(I)$
separates $n$ copies of 
$\underset{\alpha<t<\beta}{\cup} h_t(I)$
and 
so that we keep the boundary of the closure of 
$\underset{\alpha<t<\beta}{\cup} h_t(I)$
since 
there does not appear a new singularity of the immersed annulus.  
Figure \ref{step1} is an example. \\

\begin{figure}
\includegraphics[width=150mm]{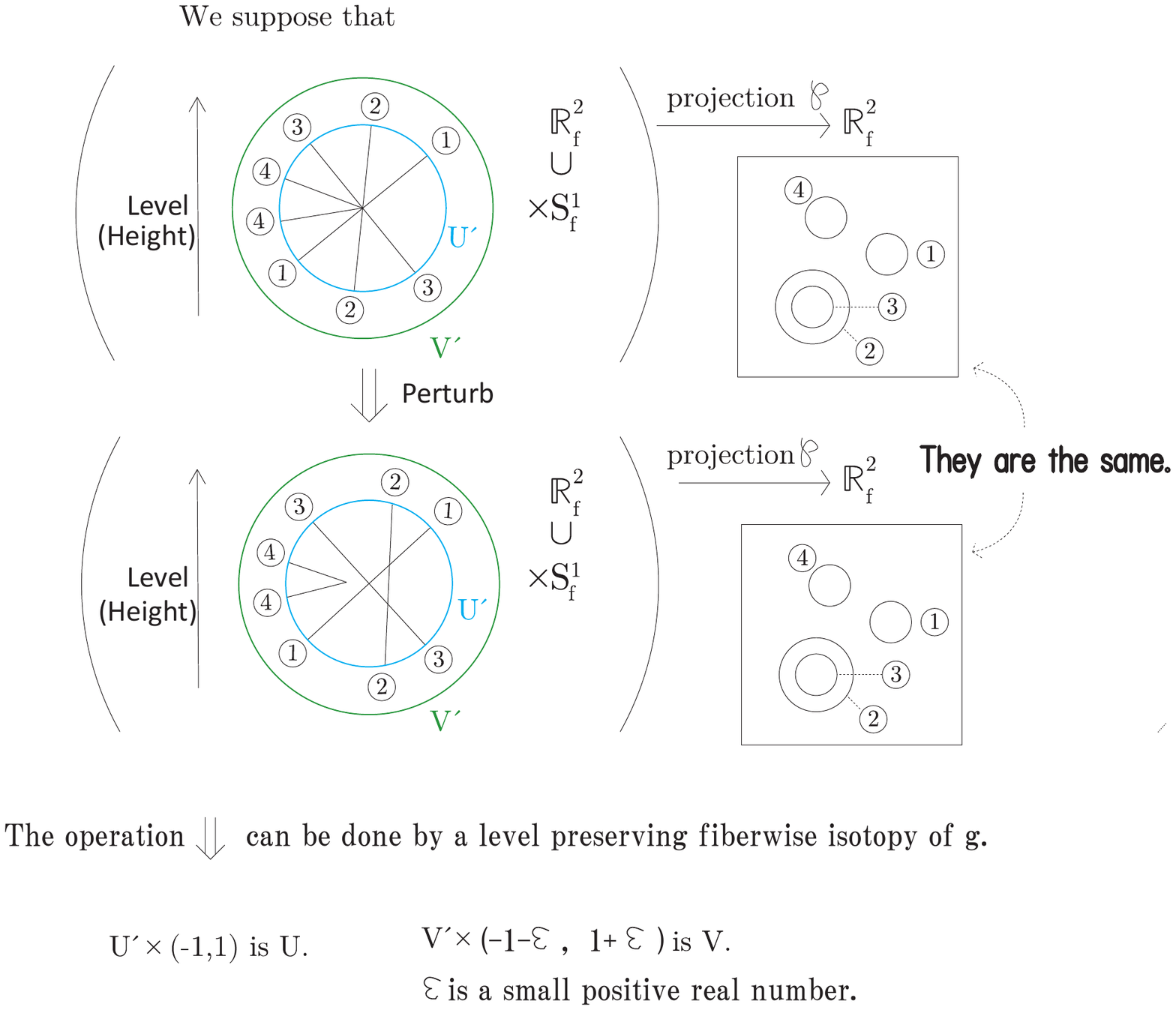}  
\vskip-35mm
\caption{{\bf A special 
isotopy of $g$. The intersection of four sheets in the upper figure 
is perturbed and is changed into the one in the lower figure.}\label{Ohio}}   
\end{figure}

\begin{figure}\includegraphics[width=140mm]{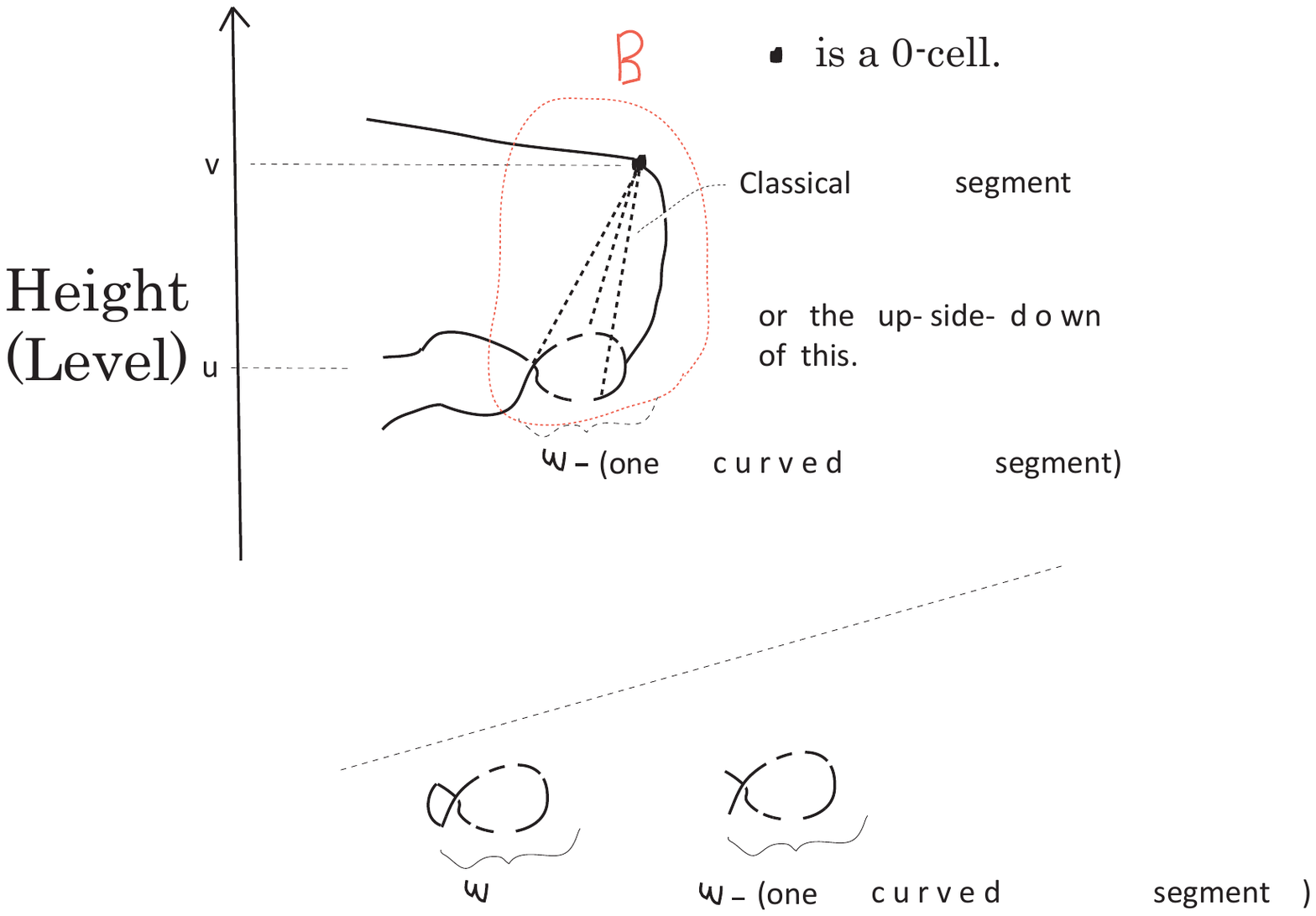}  
\vskip-40mm
\caption{{\bf 
A branch point moved by a special 
isotopy of $g$ }}\label{todome} 
 \end{figure}

\begin{figure}
    \includegraphics[width=140mm]{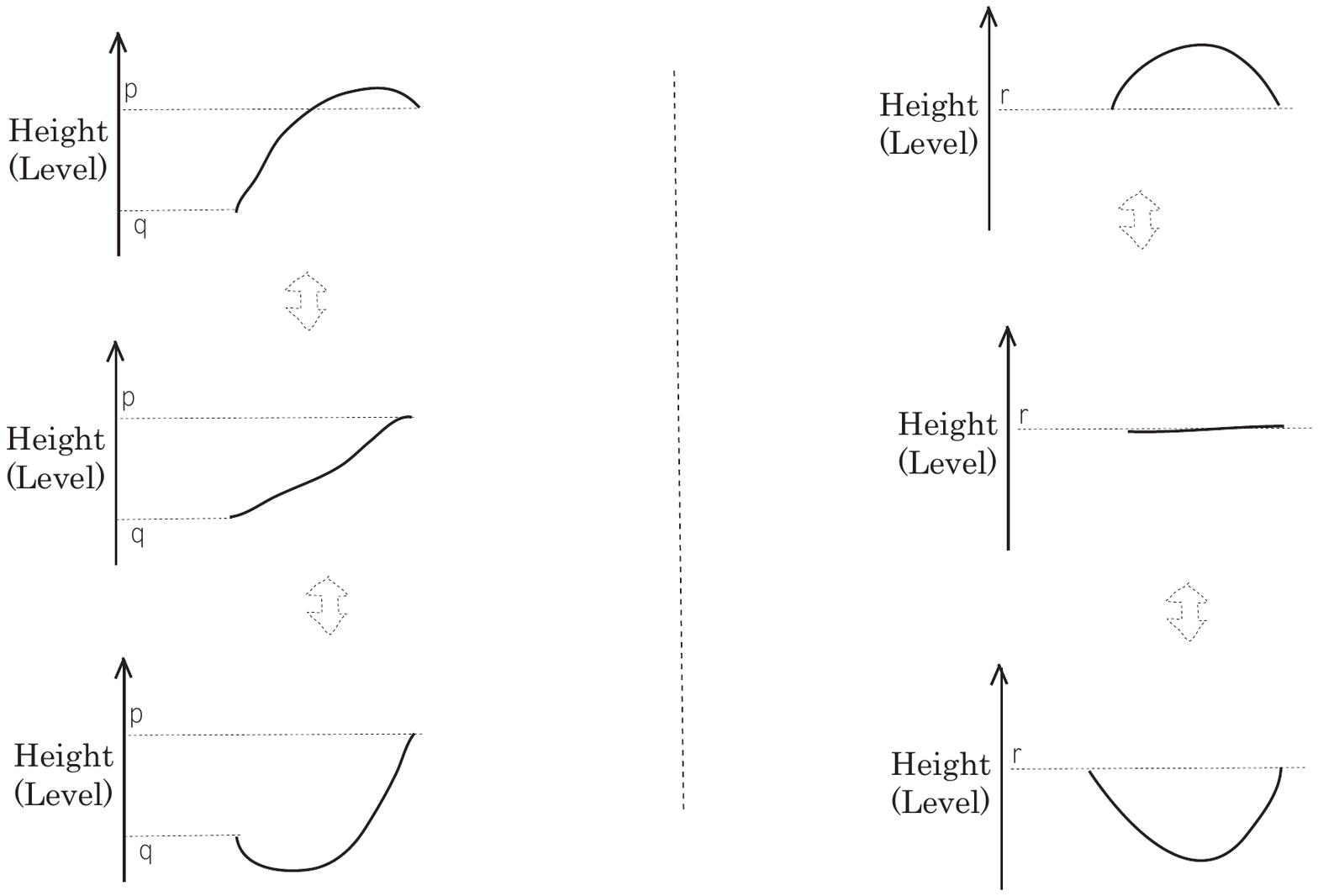}  
\vskip-20mm
\caption{{\bf Changing $X$.}\label{igaini}}   
\end{figure}

\begin{figure}
  \includegraphics[width=120mm]{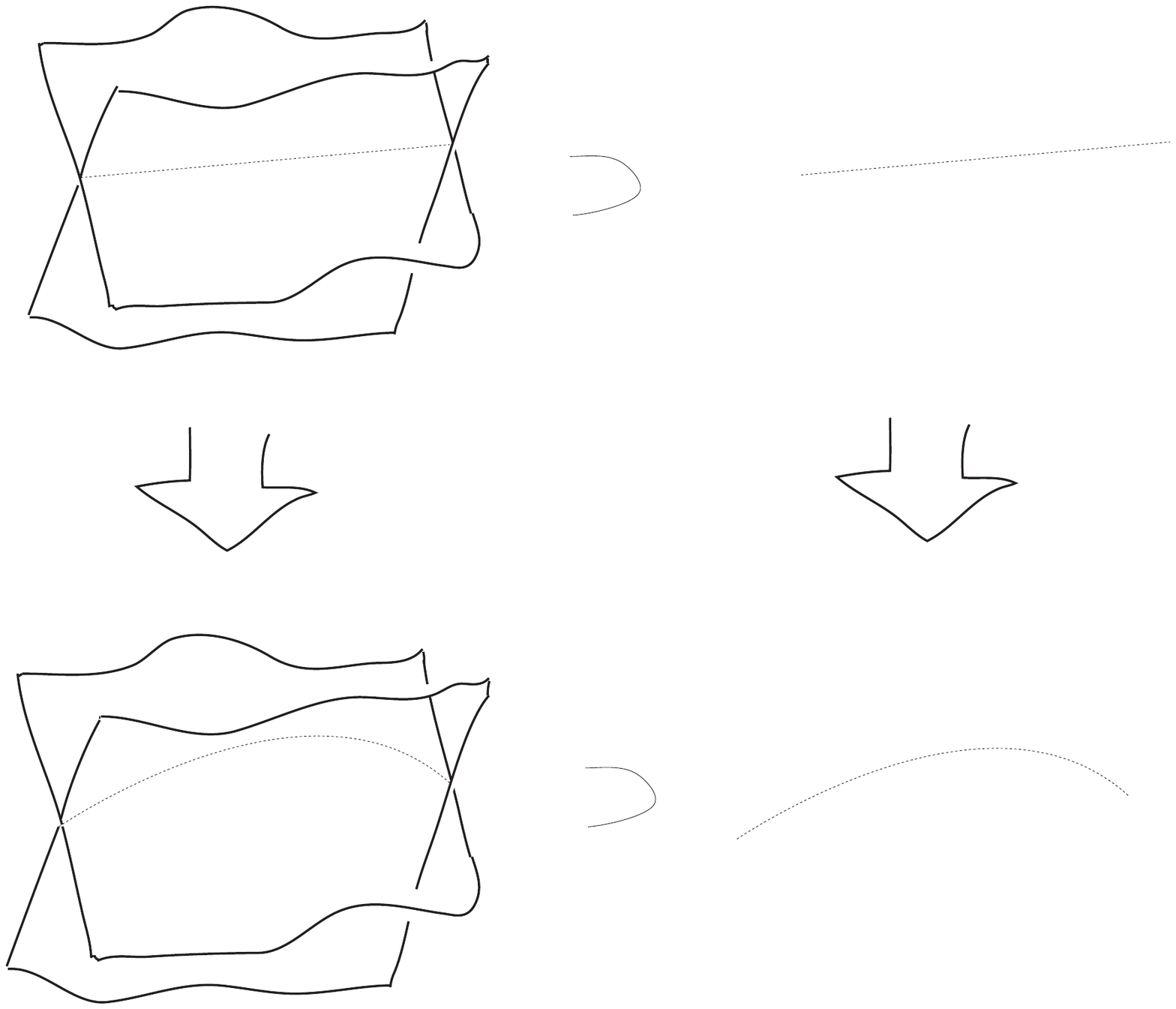}  
\caption{{\bf 
While the middle part of two sheets approaches by 
a special 
isotopy of $g$, 
the intersection $X$ in Lemma \ref{sukoshi} changes.  
}\label{beefsteak}}   
\end{figure}

\begin{note}\label{ultra}
Figures \ref{step1}-\ref{dia1} draw `figures in PL category' although the figures are smoothened. 
When we move a map by isotopy, we take a PL subdivision if necessary. 
\end{note}
%
%

Note that 
the boundary of the closure of each $I$ may have a singular point set. 
The repetition of this procedure and 
Proposition \ref{polygon} 
imply 
Claim \ref{Oregon}. \qed

\begin{note}\label{honmaya}
Note each point in the resultant part which is made 
from $\underset{\alpha<t<\beta}{\cup} h_t(I)$　by the separation.  
By the definition of $I$,  it is a single point. 
\end{note}

\noindent
{\bf The second step.}  We prove the following.

\begin{cla}\label{Pennsylvania} 
Suppose that $g$ satisfies 
the condition of Claim {\rm \ref{Oregon}}.  
We can perturb $g$ in the special way 
so that $\pi\circ g$ covers a level preserving transverse immersion $\underline{g}$ 
except for a finite number of points contained in $S^1_b\x[0,1]$ 
with the following property: 
Let $P$ be an exceptional point. 
Then $\underline{g}^{-1}(\underline{g}(P))$ may be more than one point. 
\end{cla}

\noindent
{\bf Proof of Claim \ref{Pennsylvania}.} 
Since $g$ satisfies the condition of Claim  \ref{Oregon}, 
the singular point set of Im$(\pi\circ g)$ is a 1-dimensional finite simplicial complex.  
Recall Proposition \ref{polygon} 
and Note \ref{honmaya}. 
Take the interior $I$ of a 1-simplex in the singular point set of the simplicial complex structure 
with the following property: 

\smallbreak\noindent(1)
$\underline{g}^{-1}(I)$ is disjoint $n$ open segments $I_1,...,I_n$ in $S^1_b\x[0,1]$ ($n\in\N$). 
$\underline{g}\vert_{I_i}$ is an embedding map.

\smallbreak\noindent(2)
There is an open neighborhood $U$ of $I$ in $\R^2_b\x[0,1]$
with the following property:   
There are open discs 
$D^2_i$ embedded in $S^1_b\x[0,1]$ 
 each of which is a tubular neighborhood of $I_i$ in  $S^1_b\x[0,1]$ 
for each $i$. 
$D^2_i\cap D^2_j=\phi$ for each distinct $i,j$.  
$\underline{g}|_{D^2_i}$ is an embedding map. 
$U\cap\underline{g}(S^1_b\x[0,1])$ is 
$\underline{g}(D^2_1)\cup...\cup\underline{g}(D^2_n)$. 
$\underline{g}(D^2_i)\cap\underline{g}(D^2_j)=I$ for each distinct $i,j$. 


\smallbreak 
We perturb $\underline{g}|_{(D^2_1\cup...\cup D^2_n)\cap(\underline{g}^{-1}(U))}$ 
in the special way 
below 
but we must remember Note \ref{shio}.
\smallbreak 

Let $V$ be an open neighborhood of $U$ in $\R^2_b\x[0,1]$.  
Hence $\overline{U}\subset V$. 

Let $\wp:\R^2_b\x[0,1]\x\R^2_f\to\R^2_f.$  
Combine this map $\wp$ and the diagram in Note \ref{kakanzu}:  
$$
\begin{matrix}
S^1_b\x[0,1]\x S^1_f&\stackrel{g}\to&\R^2_b\x[0,1]\x \R^2_f
&\stackrel{\wp}\to\R^2_f\\
\downarrow_\rho&\circlearrowright &\downarrow_\pi&& \\
S^1_b\x[0,1]&\stackrel{\underline{g}}\to&\R^2_b\x[0,1]&&
\end{matrix}
$$
\begin{cla}\label{shoga}
We can perturb $g$ in the special way, 
keeping out $V$ $($not $U)$,  
with the following properties:  
The image $\wp(g(\rho^{-1}(D^2_i)))$ 
is a circle $C_i$. 
We have $C_i\cap C_j=\phi$ for each distinct $i, j$. 
The map $\wp\vert_{g(\rho^{-1}(D^2_i))}$ is the projection.  
\end{cla}

\bigbreak\noindent{\bf Proof of Claim \ref{shoga}.} 
Take a point $\sigma\in I$. 
Let $\underline{g}^{-1}(\sigma)=\{\sigma_1,...,\sigma_n\}$ and $\sigma_i\in I_i$. 
Then the image of $\wp(g(\rho^{-1}(\sigma_i)))$ is a circle $C'_i$,  
and  we have $C'_i\cap C'_j=\phi$ for each distinct $i, j$. 
We can take $g$ so that the circle $C_i$ which we want is this circle $C'_i$ for each $i$. 
Then Claim \ref{shoga} holds. \qed 

\bigbreak\noindent{\bf Note.}
The reason why we prepare $V$ is 
as follows: 
Before the perturbation, the map $\wp\vert_{g(\rho^{-1}(\partial D^2_i))}$ is not a projection.  
Note $\partial D^2_i\subset\overline U$. 
However 
$\wp\vert_{g(\rho^{-1}(\partial D^2_i))}$ is the projection 
after the perturbation.  
\\

We next make $\underline{g}|_{(D^2_1\cup...\cup D^2_n)\cap(\underline{g}^{-1}(U))}$ 
 a level preserving transverse immersion 
since we can perturb $g$ in the special way, 
keeping out $U$,  with the following properties: 
For any point $q\in D^2_i$ 
and any point $r$ in the circle $\rho^{-1}(q)$,   
$\wp(g(r))\in\R^2_f$ is fixed while we perturb $g$. 
Claim \ref{shoga} ensures that while we perturb $g$ in this way, 
we keep a property that $g$ is an embedding map.  
Figure \ref{Ohio} is an example of this 
procedure. The repetition of this procedure implies Claim \ref{Pennsylvania}. \qed\\

\noindent
{\bf The third step.}  We prove the following.

\begin{cla}\label{Rhode Island}
Suppose that 
$g$ satisfies the condition of Claim {\rm\ref{Pennsylvania}}.  
We can perturb $g$ 
in the special way 
so that $\pi\circ g$ covers a level preserving transverse immersion $\underline{g}$ 
except for a finite number of points contained in $S^1_b\x[0,1]$ 
with the following property: 
Let $P$ be any exceptional point. 
The set 
$\underline{g}^{-1}(\underline{g}(P))$ consists of only one point. 
\end{cla}

\noindent
{\bf Proof of Claim \ref{Rhode Island}.} 
Assume that  
$f^{-1}(f(P))$ consists of $m$ points  
$P_1,...,P_m$ ($m\in\N$)  in $S^1_b\x[0,1]$.  
Take an open neighborhood $U$ of $f(P)$ in $\R^2_b\x[0,1]$ with the following properties: 
There are open discs 
$D^2_i$ in $S^1_b\x[0,1]$ 
which is a tubular neighborhood of $P_i$ in  $S^1_b\x[0,1]$ 
for each $i$. 
$D^2_i\cap D^2_j=\phi$ for each distinct $i,j$.  
$U\cap\underline{g}(S^1_b\x[0,1])$ is  
$\underline{g}(D^2_1)\cup...\cup\underline{g}(D^2_m)$. 
$\underline{g}(D^2_1))\cap...\cap\underline{g}(D^2_m)=f(P)$. 
For a pair $(i,j)$,  we may have 
$\underline{g}(D^2_i))\cap\underline{g}(D^2_j)\underset{\neq}{\supset}f(P)$. 

\smallbreak
We perturb $\underline{g}|_{(D^2_1\cup...\cup D^2_n)\cap(\underline{g}^{-1}(U))}$ 
in the special way 
 below 
but we must remember Note \ref{shio}.
Take an open neighborhood $V$ of $U$ in $\R^2_b\x[0,1]$ 
such that $\overline{U}\subset V$. 
 We can perturb $g$ in the special way, 
keeping out $V$ (not $U$),  
with the following properties:  
$\wp(g(\rho^{-1(}D^2_i)))$ is a circle $C_i$. 
$\wp\vert_{g(\rho^{-1(}D^2_i))}$ is the projection.  \\

We can make $\underline{g}|_{(D^2_1\cup...\cup D^2_n)\cap(\underline{g}^{-1}(U))}$ 
 a level preserving transverse immersion except for a finite number of points 
since we can perturb $g$ in the special way, 
keeping out $U$,  
with the following properties: 
For any point $q\in D^2_i$ 
and any point $r$ in the circle $\rho^{-1}(q)$,   
$\wp(g(r))\in\R^2_f$ is fixed while we perturb $g$. 
(Note that while we perturb $g$ in this way, 
we keep a property that $g$ is a 
embedding map.)  
The repetition of this procedure 
and Note \ref{honmaya} imply 
Claim \ref{Rhode Island}. \qed\\

\noindent{\bf The fourth step.}   
Take $\pi\circ g$ and $\underline{g}$ in Claim \ref{Rhode Island}. 
Let $P$ be any exceptional point. Recall that $P\in S_b^1\x[0,1]$.  
Let $N(P)$ be the tubular neighborhood of $P$ in $S_b^1\x[0,1]$.  
Take the tubular neighborhood $B$ of $\underline{g}(P)$ in $\R^2_b\x[0,1]$.  
We can suppose that 
$\underline{g}(N(P))\subset  B$ 
and that 
$\underline{g}(\partial N(P))\subset  \partial  B$.  
The image $\underline{g}(N(P))$ makes $\underline{g}(P)$ 
a branch point 
(recall Definition \ref{oyster}). 
Here we ignore the information 
of fiber circles over $P$.  
The information of Rourke fiber makes $\underline{g}(\partial N(P))\subset  \partial  B$, 
 a virtual 1-knot diagram $\omega$ in $\partial  B-$(a point). 
Note that  $\partial  B-$(a point) is the 2-space and 
that the point is not included in $\omega$.  
Recall virtual segments defined in Note \ref{vrei}. 
A {\it classical segment} is the segment with the following properties. 
It is a segment included in a virtual 2-knot diagram.   
One of the boundary is a classical branch point.
The points in the interior of the segment are  classical double points. 
An example is drawn in Figure \ref{sashimiv} if the branch point there is a classical branch point.  \\

\begin{cla}\label{mochi} 
We can assume that all branch points of Im $\underline{g}$ are classical Whitney branch points.  
\end{cla} 

\noindent{\bf Proof of Claim \ref{mochi}.}
Since $\underline{g}(P)$ is a branch point, 
$n$ virtual segments and $m$ classical segments 
meet at 
$\underline{g}(P)$, 
where $\{n,m\}\subset\N\cup\{0\}$ and $n+m\geqq2$.  
We can prove that there is no virtual segment 
in the same fashion as the one in the proof of Theorem \ref{Rmuri}. 
(Note that in \S\ref{v2} we proved Theorem \ref{Rmuri} in the smooth category 
but we can prove the PL version of Theorem \ref{Rmuri} in the same way.)  
Therefore 
more than one classical segment meet at $\underline{g}(P)$.  
Hence $\omega$ is a classical diagram and 
determines a classical 1-knot. 
\\

In order to complete the proof of Claim \ref{mochi},  we will prove Claim \ref{koma}. 
In order to prove Claim \ref{koma}, we prove the following Claim \ref{tantei}.

\begin{defn}\label{prod}
Let $u,v\in[0,1]$. Let $u\leq t\leq v$. 
The map $\underline{g}|_{S^1_b\x[u,v]}$ is called a {\it product map} 
if there is an isotopy $\iota_t$
  of $\R^2$ from the identity map 
such that  $\iota_t:\R^2\x\{u\}\to\R^2\x\{t\}$  
carries $({\rm Im}\underline{g})\cap(\R^2\x\{u\})$ to $({\rm Im}\underline{g})\cap(\R^2\x\{t\})$. 
Let $B^3$ be an embedded closed 3-ball in $\R^2_b\x[0,1]$. 
The map $\underline{g}|_{S^1_b\x[u,v]}$ is called a {\it product map out $B^3$} 
if there is an isotopy $\iota_t$
   of $\R^2$ from the identity map 
such that  $\iota_t:\R^2\x\{u\}\to\R^2\x\{t\}$  
carries $({\rm Im}\underline{g}-B^3)\cap(\R^2\x\{u\})$ 
to $({\rm Im}\underline{g}-B^3)\cap(\R^2\x\{t\})$. 
\end{defn}

We have the following.

\begin{cla}\label{tantei}
By using a special 
isotopy of $g$,  
any branch point is moved 
as drawn in Figure $\ref{todome}$: 
Let $\alpha_u$ $($respectively, $\alpha_v)$ be an immersed circle determined by  
$\underline{g}(S^1_b\x\{u\})\subset\R^2_b\x\{u\}$ 
$($respectively, $\underline{g}(S^1_b\x\{v\})\subset\R^2_b\x\{v\})$
with the information of Rourke fiber determined by 
$g(S^1_b\x\{u\}\x S^1_f)\subset\R^2_b\x\{u\}\x\R^2_f$ 
$($respectively, \\$g(S^1_b\x\{v\}\x S^1_f)\subset\R^2_b\x\{v\}\x\R^2_f).$     
The map $\underline{g}|_{S^1_b\x[u,v]}$ is a product map out $B$. 
Hence $\alpha_u=\alpha_v\#\omega$, 
where $\#$ denotes the connected sum of immersed circles into $\R^2$ 
and $=$ means that there is an orientation preserving diffeomorphism 
of $\R^2$ which carries the left hand side to the right side one. 
\end{cla}

\noindent{\bf Proof of Claim \ref{tantei}.} 
For each $t$, $(\R^2_b\x\{t\})\cap$Im $\underline{g}$ is connected. 
Hence 
the branch point is not a local maximal (respectively, minimal) point of 
the restriction of the height function $\R^2_b\x[0,1]\to[0,1]$ to  Im $\underline{g}$. 
\\

Claim \ref{sukoshi} follows from Claim \ref{tako}.

\begin{cla}\label{sukoshi}
Let $X$ be the closed interval contained in $\mathcal S$. Assume that $X$ does not have self-intersection.  
Then, by using a special 
isotopy of $g$, 
we can move {\rm Int}$X$ as drawn in Figure \ref{igaini} 
with the following properties: 
We move {\rm Int}$X$ by an isotopy of embedding of  {\rm Int}$X$,  
keeping  the position of $\partial X$ in $\R^2_b\x[0,1]$.  
We keep the position of $\overline{\mathcal S-X}$ in $\R^2_b\x[0,1]$. 
We keep the condition $X\cap(\mathcal S-X)=\phi.$
\end{cla}

In Figure \ref{beefsteak} there is an example of Claim \ref{sukoshi}.  
Threre is drawn 
how $X$ changes by a special 
isotopy of $g$ 
in the case of the upper two figures in 
the right column of Figure \ref{igaini}. 
Note that Int$X$ consists of double points. 
Each point of $\partial X$ is a branch, double or triple one.
Let $B$ be an open disc contained in $S^1_b\x[0,1]\x S^1_f$. By Definition \ref{PLNevada}.(1), $\pi\circ g(B)$ is not parallel to $\R^2_b\x\{0\}$. 
Note that if $\pi\circ g(B)$ is parallel to $\R^2_b\x\{0\}$, 
the phenomenon in the right column of Figure \ref{beefsteak} does not occur. \\

\bigbreak

\begin{cla}\label{tako} 
Let $B^3$ be a closed $($respectively, open$)$ 3-ball embedded in $R^2_b\x[0,1]$.  
Take any orientation preserving isotopy of 
diffeomorphism of $B^3$ fixing $\partial B^3$ from the identity map. 
We can give a coordinate $(x,y,t)$ to $p\in B^3\subset\R^2_b\x[0,1]$.  
Suppose that this isotopy carries $p$ to a point whose coordinate is $(x,y,t')$, 
where $t'\neq t$ or $t'=t$ holds. 
Use this isotopy and make a homotopy of $\underline{g}.$  
Suppose that this homotopy is a special 
homotopy of $\underline{g}$. 
Then this homotopy of $\underline{g}$ can be covered by a 
special 
isotopy of $g$. 
\end{cla}

By Claim \ref{sukoshi}, 
we can let the interior of all classical segments exist below (respectively, over) the branch point 
with respect to the height 
as drawn in Figure \ref{todome}.

By the first, second and third steps, 
the singular point set of Im $\underline{g}$
is a finite simplicial complex. 
Hence we have the following. 

\begin{cla}\label{koma}
There are only finite number of $t\in[0,1]$ with the following properties: 
There is no real number $\varepsilon$ such that 
the map $\underline{g}|_{S^1_b\x[t-\varepsilon, t+\varepsilon]}$ 
is a product map out $B$. 
\end{cla}

By Claims \ref{tako} and \ref{koma},  we have Claim \ref{tantei}.   
\qed\\

\begin{cla}\label{kinako}
$\omega$ defines the trivial knot. 
Hence we obtain $\alpha_v$ from $\alpha_u$ 
by using only classical Reidemeister moves. 
\end{cla} 

\noindent{\bf Proof of Claim \ref{kinako}.}
By the map $g|_{S^1_b\x[u,v]\x S^1_f}$, 
$\alpha_u$ and $\alpha_v$ are fiberwise equivalent. 
Therefore 
the submanifolds, $\mathcal S(\alpha_u)$ and $\mathcal S(\alpha_v)$, of $\R^4$ 
are isotopic. 
Hence 
$\pi_1(\R^4-\mathcal S(\alpha_u))\\\cong
\pi_1(\R^4-\mathcal S(\alpha_v))$. 
Hence 
the group of $\alpha_u$ and that of $\alpha_v$ are isomorphic. 
Since $\alpha_u\\=\alpha_v\#\omega$, 
the group of $\alpha_u$ is 
the free product of that of $\alpha_v$ and that of $\omega$.  
Hence 
the group of $\omega$ is $\Z$. 
Since $\omega$ defines a classical 1-knot,  
 $\omega$ defines the trivial 1-knot. 
Since $\omega$ is a classical 1-knot diagram and represents the trivial 1-knot, 
 $\omega$ is changed into the trivial 1-knot diagram by using only classical Reidemeister moves. 
Hence  Claim \ref{kinako} holds. 
\qed
\\

It is easy to prove that 
if two virtual 1-knot diagrams are obtained each other 
by using only classical Reidemeister moves, 
they are fiberwise equivalent. 
Therefore, 
we change $\underline{g}|_{S^1_b\x[u,v]}$  in $B$ 
so that we let 
the map $\underline{g}|_{S^1_b\x[u,v]}$ be a level preserving, 
generic map. 
Hence the following holds: If $\underline{g}(S^1_b\x[u,v])$ includes a branch point, 
it is the classical Whitney branch point. 
These classical Whitney branch points appear 
when we carry out classical Reidemeister $I$ move 
while we change $\alpha_u$ into $\alpha_v$. 
After repeating this procedure, 
all branch points of Im $\underline{g}$ are classical Whitney branch points.  
This completes the proof of Claim \ref{mochi}. \qed\\

This completes the proof of Claim \ref{kamen}. \qed\\

This completes the proof of Theorem \ref{Wyoming}. \qed \\

Claim \ref{kinako} implies the following Proposition \ref{amakara}. 

\begin{defn}\label{shuza}
Virtual 1-knot diagrams $\alpha$ and $\beta$ are said to be 
{\it strongly fiberwise equivalent} 
if $\alpha$ and $\beta$ satisfy the conditions  
which are made by replacing the phrase 
`level preserving generic map' with  
`level preserving transverse immersion' 
without changing other parts in  Definition \ref{gene}. 
\end{defn}

\begin{pr}\label{amakara} 
If virtual 1-knot diagrams $\alpha$ and $\beta$ are fiberwise equivalent, 
there is a sequence of virtual 1-knot diagrams, 
$\alpha=\alpha_1,\beta_1,\alpha_2,\beta_2,...,
\alpha_{k-1},\beta_{k-1},\alpha_k,\beta_k=\beta$ 
$(k\in\N)$,  
such that 
$\alpha_i$ and $\beta_i$ 
are strongly fiberwise equivalent
$(1\leqq i\leqq k)$ 
and such that 
$\beta_i$ and $\alpha_{i+1}$  $(1\leqq i\leqq k-1)$ 
are classical move equivalent $($and therefore rotational welded equivalent$)$. 
\end{pr}

We prove the following theorem. 

\begin{thm}\label{ike}
If $g$ satisfies Definition $\ref{shuza}$, 
then the following hold.  
%
%
Let $\mathcal S$ be the singular point set of 
$(\pi\circ g)(S^1_b\x[0,1]\x S^1_f)$ in Definition $\ref{shuza}.$   
Note that $\mathcal S$ is a finite 1-dimensional simplicial complex. 

\smallbreak\noindent{\rm (i)} 
$\mathcal S\cap(\R^2_b\x\{0$ $($respectively, $1)\})$ is 
a set of virtual and classical crossing points of $\alpha$ $($respectively, $\beta),$   
and therefore is a set of double points. 
It consists of 0-simplices. 
Only one 1-simplex 
is attached to each of these 0-simplices. 
These 1-simplices 
meet \\
$\R^2_b\x\{0$ $($respectively, $1)\}$ transversely. 

\smallbreak\noindent{\rm (ii)} 
Triple points are 0-simplices. 
$($Recall 
Notes $\ref{umeboshi}$ and $\ref{faso}$, Definnition $\ref{JW}.)$

\smallbreak\noindent{\rm (iii)} 
The restriction of 
`the height function $\mathfrak h:\R^2_b\x[0,1]\to[0,1]$'  
to the interior of any 1-simplex 
in $\mathcal S$ 
has no critical point. 
$($Hence we have the following:   
For each $\zeta\in[0,1]$, $\mathcal S\cap(\R^2_b\x\{\zeta\})$ is a finite number of points. 
No 1-simplex is parallel to  $\R^2_b\x\{0\}.)$  

\smallbreak\noindent{\rm (iv)} 
Let $\zeta\in(0,1).$
$\mathcal S\cap(\R^2_b\x\{\zeta\})$ includes no or only one 0- simplex. 

\smallbreak\noindent{\rm (v)} 
In $\R^2_b\x(0,1)$,  0-simplices 
appear only in the two 
cases of Figure $\ref{zero}$.

\begin{figure}
\vskip-30mm
\includegraphics[width=135mm]{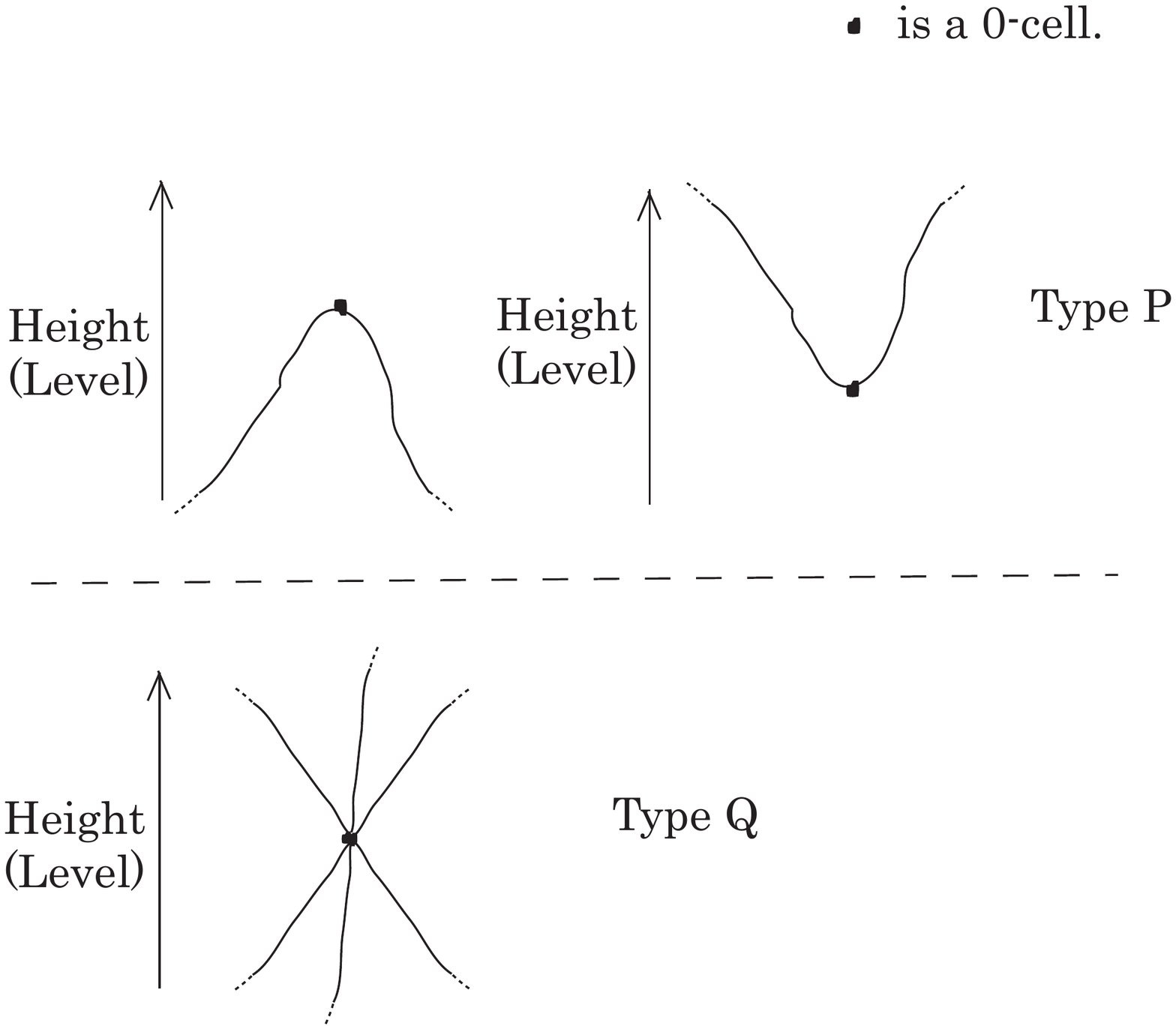}  
\vskip-30mm
\caption{{\bf 0-simplices in $\mathcal S$ }}
\label{zero} 
\end{figure}

\end{thm}

\noindent
{\bf Proof of Theorem \ref{ike}.} 
Theorem \ref{ike}.(i)  
follows from Definition \ref{shuza}.(5) for any simplicial complex structure. 
Theorem \ref{ike}.(ii) 
holds for any simplicial complex structure 
by the definition of simplicial complex structure. 
\\

Proof of Theorem \ref{ike}.(iii).  
Suppose that there is an $X$ 
as is an example explained in Claim \ref{sukoshi} and Figure \ref{beefsteak}.  
Repeating this procedure 
we can take a simplicial complex structure 
such that 
(any 1-simplex)$\cap(\R^2_b\x\{t\})$ 
for any $t\in[0,1]$ is a finite number of points.
Therefore 
the restriction of $\mathfrak h$ to the interior of any 1-simplex of this simplicial complex structure 
has a finite number of critical points. 
Make a new simplicial complex structure 
so that the critical points are new 0-simplicies 
so that we keep  the condition of Theorems \ref{ike}.(i) and (ii). 
Suppose that 
there is a 0-simplex 
$e^0$ to which only two 1-simplices 
$e^1_1$ and $e^1_2$, attach,  
and that 
$e^0$ is not a critical point of the restriction of  $\mathfrak h$ to 
(Int$e^1_1)\cup e^0\cup$(Int$e^1_2$)
=Int($e^1_1\cup e^0\cup e^1_2$). 
Make a new simplicial complex structure 
such that 
$e^1_1\cup e^0\cup e^1_2$  
is changed into a new 1-simplex 
without changing other simplicial complex structure. 
This completes the proof of Theorem \ref{ike}.(iii).  
\\

Theorem \ref{ike}.(iv) 
holds because, by Claims \ref{koma} and \ref{tako}, 
we can change the height of any 0-simplex 
so that we keep the condition of Theorems \ref{ike}.(i)-(iii). 
\\

Proof of Theorem \ref{ike}.(v).  
There are only two 
cases: 
(P) Only two 1-simplices attach to a 0-simplex. 
(Q) Only six 1-simplices attach to a 0-simplex. 
Note that, by Theorem \ref{ike}.(iii), each 1-simplex is attached to two different 0-simplices.
In the case (P), 
by Theorem \ref{ike}.(iii), 
the 0-simplex 
exists as drawn in Type P of Figure \ref{zero}.  
In the case (Q), 
as drawn in  Figure \ref{beefsteak} associated with Claim \ref{sukoshi}, 
we can move 1-simplex  
so that we have the condition as drawn in Type Q of Figure \ref{zero}, and  
so that we keep the condition of Theorems \ref{ike}.(i)-(iv). 
See Figure \ref{tsuika2} for an example of this move.

\begin{figure}
\includegraphics[width=150mm]{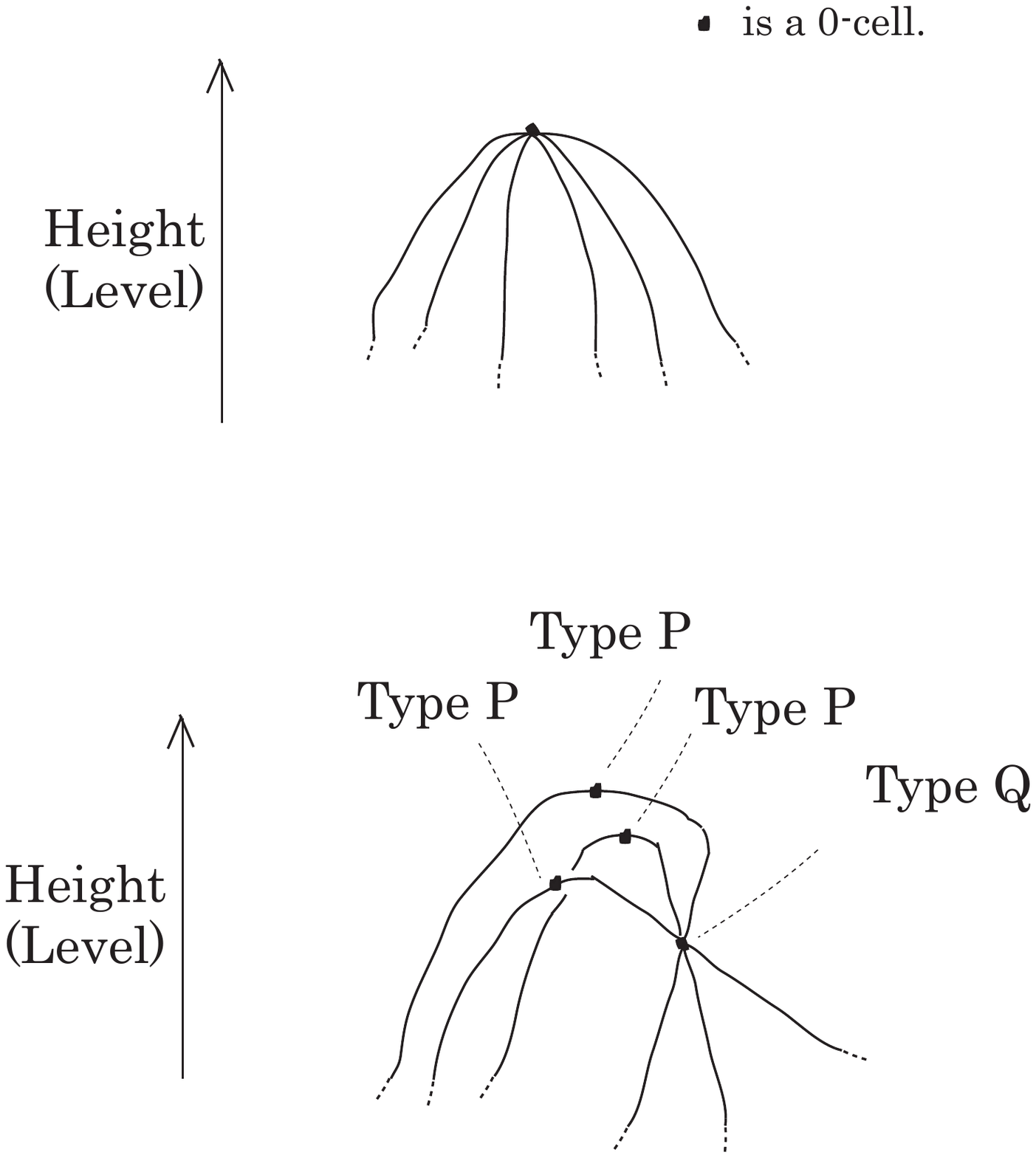}  
\vskip-40mm
\caption{{\bf 
The singularity in the upper figure is made into the one in the lower which consists of one Type Q and three Type P.}}
\label{tsuika2} \end{figure}

This completes 
the proof of Theorem \ref{ike}.(v).\\  


This completes the proof of Theorem \ref{ike}. 
\qed
\\

We have the following theorem.

\begin{thm}\label{fwrw} 
Two virtual 1-knot diagrams $\alpha$ and $\beta$ are PL fiberwise equivalent if and only if  
$\alpha$ and $\beta$ are PL rotational welded equivalent.  
\end{thm}

\noindent
{\bf Proof of Theorem \ref{fwrw}.} 
The `if' part is easy. 

We prove the `only if' part. 
By Proposition \ref{amakara}, it suffices to prove Claim \ref{takusan}

\begin{cla}\label{takusan}
Two virtual 1-knot diagrams $\alpha$ and $\beta$ are 
PL strongly fiberwise equivalent only if  
$\alpha$ and $\beta$ are PL rotational welded equivalent.  
\end{cla}

\noindent{\bf Proof of Claim \ref{takusan}.}
Let $C_\zeta=\underline{g}(S^1_b\x\{\zeta\})=$
$\underline{g}(S^1_b\x[0,1])
\cap(\R^2_b\x\{\zeta\})\\=$
$(\pi\circ g(S^1_b\x[0,1]\x S^1_f))
\cap(\R^2_b\x\{\zeta\})$.  
By Theorem \ref{ike},  
$C_\zeta$ is 
an immersed circle in $\R^2_b\x\{\zeta\}$ and 
its singular point set is a finite number of points.  
$C_\zeta$ changes from $\alpha$ to $\beta$ step by step 
as $\zeta$ runs from 0 to 1.  
If, for a $\zeta_p$, $C_{\zeta_p}$ includes a 0-simplex  
of the simplicial complex structure 
in Theorem \ref{ike}.  
A classical or virtual Reidemeister move is done there. 
We do any of them only there. 
If $\zeta_q<\zeta<\zeta_r$, 
  $C_\zeta$ includes no 0-simplex. 
Then $C_\zeta$ is not changed while $\zeta$ runs from $\zeta_q$ to $\zeta_r$. 
We investigate how $C_\zeta$ changes in detail. 
Near a 0-simplex  
in $\R^2_b\x[0,1]$, 
Im $\underline{g}$ is drawn as in Figure \ref{tsuika}  
since  $\underline{g}$ is a transverse immersion.  
Here, note that we can move $\mathcal S$ by using 
a special 
isotopy of $g$.   \\

\begin{figure}
\includegraphics[width=110mm]{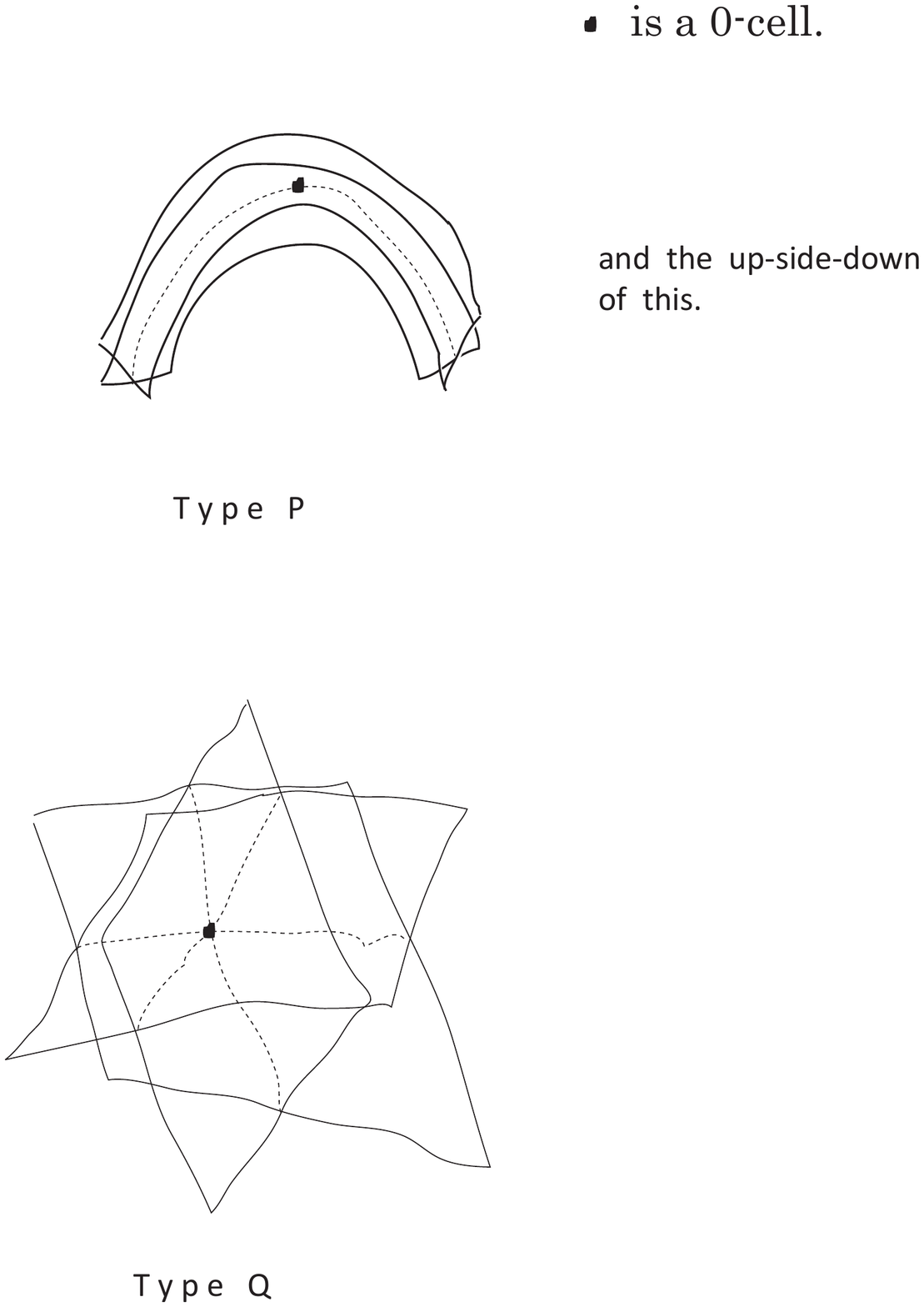}  
\vskip-10mm
\caption{{\bf How sheets intersect near Types P and Q.}}
\label{tsuika} \end{figure}

Therefore 
we have only the following two facts on $\mathcal S$ and local moves on the knot diagrams.

\bigbreak
\noindent (i)  
%
Let $\sigma, \tau\in[0,1]$. 
Suppose that $\mathcal S\cap(\R^2_b\x(\sigma,\tau))$ includes no 0-simplex. 
It holds that $\underline{g}|_{S^1_b\x[\sigma,\tau]}$ is a product map. 
Then $(\pi\circ g)(S^1_b\x[0,1]\x S^1_f)\cap(\R^2_b\x\{\sigma\})$ can be obtained from  $(\pi\circ g)(S^1_b\x[0,1]\x S^1_f)\cap(\R^2_b\x\{\tau\})$ by an isotopy of $\R^2_b$.

\bigbreak
\noindent (ii)
Let $\xi\in[0,1].$
Suppose that $\mathcal S\cap(\R^2_b\x\{\xi\})$ includes only one 0-simplex. 
Suppose that 
$\mathcal S\cap(\R^2_b\x(\xi, \xi+\varepsilon])$ 
(respectively, 
$\mathcal S\cap(\R^2_b\x[\xi-\varepsilon, \xi))$) 
includes no 0-simplex. 
Let 
$D=(\pi\circ g)(S^1_b\x[0,1]\x S^1_f)\cap(\R^2_b\x\{\xi-\varepsilon'\})$ 
and 
 $U\\
=(\pi\circ g)(S^1_b\x[0,1]\x S^1_f)\cap(\R^2_b\x\{\xi+\varepsilon'\})$. 
If the 0-simplex 
is put in Type P or Q, 
then $U$ is obtained from $D$ by 
one welded move other than a virtual Reidemeister $I$ move.  
(Note. 
Type P causes classical and virtual Reidemeister $II$ moves.  
Type Q causes classical and virtual Reidemeister $III$ moves.  
Four types of triple points correspond to four types of classical and virtual Reidemeister $III$ moves.) 
Therefore 
$\alpha$ is changed into $\beta$ by welded moves other than the virtual Reidemeister $I$ move. 
Hence $\alpha$ is rotational welded equivalent to $\beta$.
This completes the proof of Claim \ref{takusan}.  \qed\\

This completes the proof of Theorem \ref{fwrw}. 
\qed\\

We will complete the proof of Theorem \ref{smooth} and go back to the smooth category.
We said that the `if' part of  Theorem \ref{smooth} is easy. 
We will prove the `only if' part of  Theorem \ref{smooth}  by using the following lemma. 

\begin{lem}\label{PLtosmooth}
Let $\alpha$ and $\beta$ be smooth virtual 1-knot diagrams. 
Let $\alpha'$ $($respectively, $\beta')$ be a PL virtual 1-knot diagram  
which is piecewise smooth, planar, ambient isotopic to $\alpha$ $($respectively, $\beta)$. 
Then we have the following. 
$\alpha$ and $\beta$ are smooth rotational welded equivalent if and only if 
$\alpha'$ and $\beta'$ are PL rotational welded equivalent.
\end{lem}

\noindent{\bf Proof of Lemma \ref{PLtosmooth}.}
Let $\xi$ and $\zeta$ be smooth virtual 1-knot diagrams. 
If $\xi$ and $\zeta$ are PL, planar, ambient isotopic to a PL virtual 1-knot diagram $\gamma$,  
then  $\xi$ is smooth, planar, ambient isotopic to $\zeta$. 
{\it Reason.} Smoothen the corner of $\gamma$.  
Each of PL rotational welded Reidemeister moves is regarded as smooth rotational welded Reidemeister move. 
This completes the proof of Lemma \ref{PLtosmooth}.
\qed\\

Assume that two virtual 1-knot diagrams $\alpha$ and $\beta$ are smooth fiberwise  equivalent.  
By Claim \ref{xbeef},  they are PL fiberwise  equivalent. 
By Theorem \ref{fwrw}, they are PL rotational welded equivalent.  
By Lemma \ref{PLtosmooth}, they are smooth rotational welded equivalent.  
Therefore the `only if' part of Theorem \ref{smooth} is true. 

This completes the proof of Theorem \ref{smooth}. \qed\\

We are now back to the smooth category. 
\\

\noindent {\bf Note.}  
Figure \ref{dia1} explains Figure \ref{tsuika2} in more detail.

\begin{figure}  \vskip-10mm  \includegraphics[width=120mm]{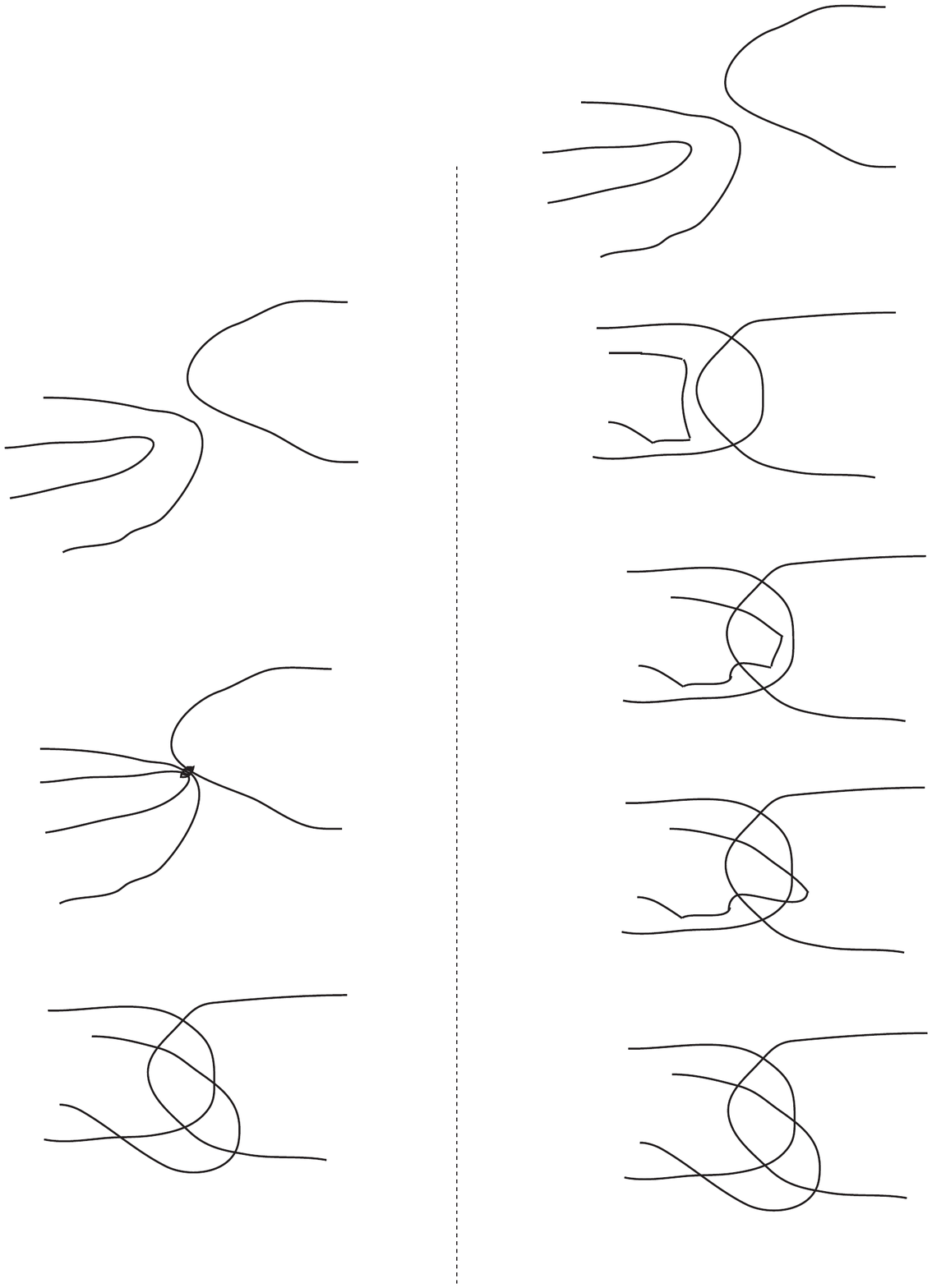}  
\caption{{\bf 
This pair of the left figure and the right one  is an example of a pair of the figures of Figure \ref{tsuika2}. The left (respectively, right) figure is an example of a sequence of diagrams associated with the upper (respectively, lower) figure of Figure \ref{tsuika2}. The left sequence is perturbed and is made into the right sequence.These diagrams are drawn  without information of virtual multiple points and classical ones.  
}}\label{dia1} \end{figure}


\begin{note}\label{xudon}
In \cite{Rourke}, the fiberwise equivalence is defined by the following definition. 
We call the equivalence relation the {\it $f$-fiberwise equivalence} in this paper.  
Note that we work in the smooth category.  

\begin{defn}\label{yugentsuki}
Add the following condition to Definition \ref{Nevada} without changing the other parts. 
We call the equivalence relation the {\it $f$-fiberwise equivalence}.  
Note that we work in the smooth category. 

\smallbreak\noindent$(\ref{yugentsuki}.1)$  In each fiber  $\R^2_f$,  there are a finite number of circles. 
$($That is, $<\infty.)$
\end{defn}

We said that it is easy to prove the following (i). It is also easy to prove the \\following (ii). 

\smallbreak\noindent 
(i) If virtual 1-knot diagrams $\alpha$ and $\beta$ are rotational welded equivalent, 
then  $\alpha$ and $\beta$ are fiberwise equivalent.

\smallbreak\noindent 
(ii) If virtual 1-knot diagrams $\alpha$ and $\beta$ are rotational welded equivalent, 
then  $\alpha$ and $\beta$ are $f$-fiberwise equivalent.

\bigbreak
 Theorem \ref{smooth} and the above (ii) imply the following (iii). 

\smallbreak\noindent 
(iii) If virtual 1-knot diagrams $\alpha$ and $\beta$ are fiberwise equivalent, 
then  $\alpha$ and $\beta$ are $f$-fiberwise equivalent. 
(Note that we work in the smooth category.)
\\

The converse of (iii) is trivial. Hence we have the following: 
Virtual 1-knot diagrams $\alpha$ and $\beta$ are fiberwise equivalent  
if and only if $\alpha$ and $\beta$ are $f$-fiberwise equivalent
\end{note}

\bigbreak
\begin{note}\label{xmikan}
Although Rourke claimed in \cite[Theorem 4.1]{Rourke}  that 
two virtual 1-knot diagrams $\alpha$ and $\beta$ are 
fiberwise equivalent if and only if  $\alpha$ and $\beta$ are welded equivalent 
in the PL (respectively, smooth) category,  
we state that this claim is wrong, as we mentioned it 
in the last few paragraphs of \S\ref{i3}.  
The reason for the wrongness is Theorems \ref{smooth} and \ref{fwrw} and Claim \ref{panda}. 
\end{note}

We introduce a new equivalence relation of the set of virtual 1-knot diagrams. 

\begin{defn}\label{parity}  
Let $\alpha$ and $\beta$ be virtual 1-knot diagrams. 
We say that 
$\alpha$ and $\beta$ are {\it virtually parity equivalent} 
if $\alpha$ and $\beta$ have the same parity of virtual crossing points. 
\end{defn} 

We prove several results associated with the virtual parity.

\begin{cla}\label{hirumeshi}
If two virtual 1-knot diagrams 
$\alpha$ and $\beta$ are 
rotational welded equivalent, 
then $\alpha$ and $\beta$ are virtually parity equivalent. 
\end{cla} 

\noindent{\bf Proof of Claim \ref{hirumeshi}.}
We can obtain $\alpha$ from $\beta$
by some welded-moves 
other than virtual Reidemeister $I$ move. 
\qed

\begin{cla}\label{panda}
The welded equivalence does not imply the rotational welded equivalence.  
\end{cla} 

\noindent{\bf Note.} 
By their definitions, the rotational welded equivalence implies the welded equivalence.  
\\

\begin{figure}
     \includegraphics[width=80mm]{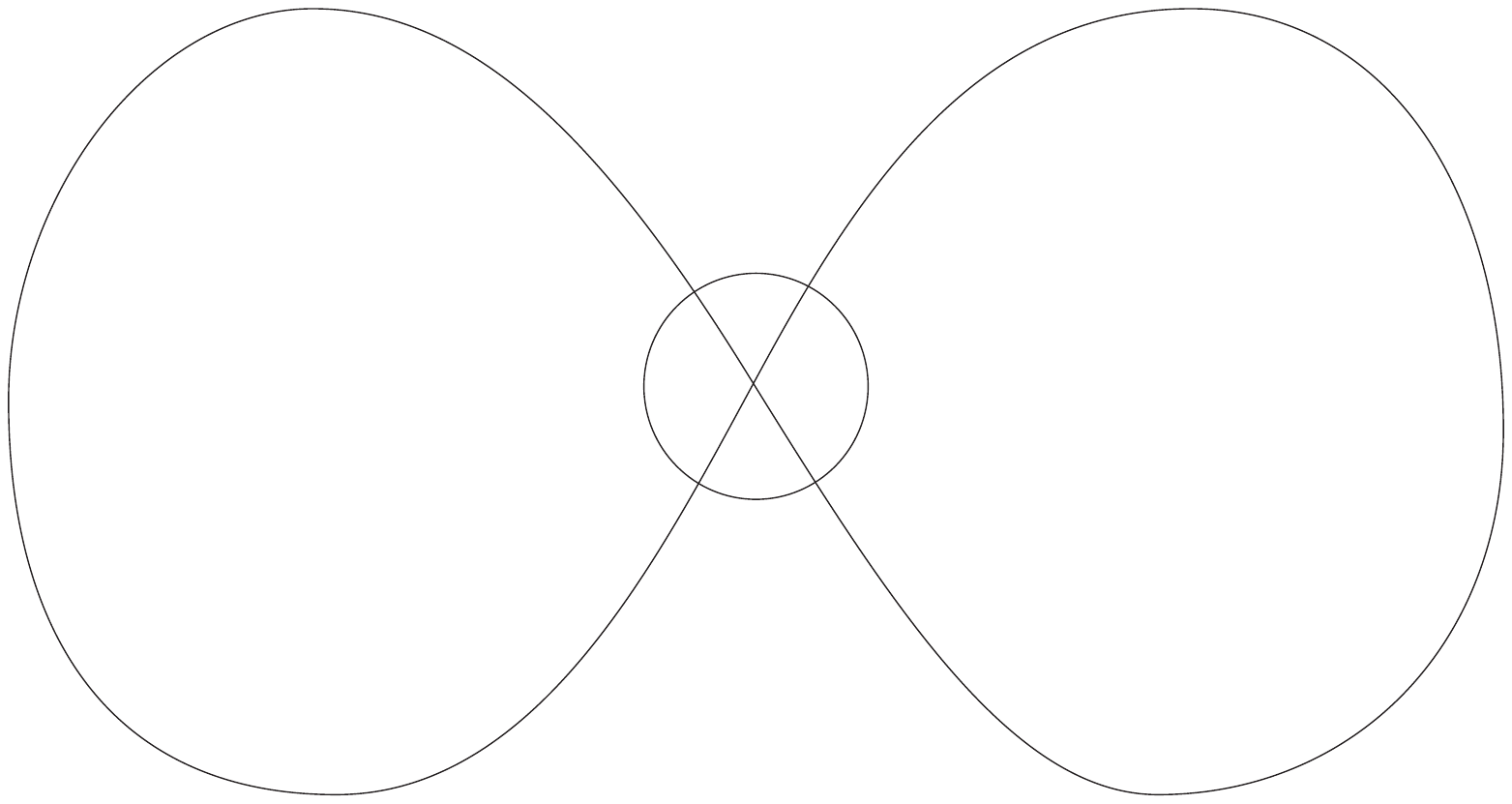}  
\vskip-30mm
\caption{{\bf The virtual figure $\infty$ knot diagram}\label{New Hampshire}}   
\end{figure}

\noindent{\bf Proof of Claim \ref{panda}.}
Call the virtual 1-knot diagram in Figure \ref{New Hampshire},    
the {\it virtual figure $\infty$ knot diagram}. 

The virtual figure $\infty$ knot diagram 
and  the trivial 1-knot diagram 
are welded equivalent by the definition  
but are not rotational welded equivalent 
by Claim \ref{hirumeshi}.

This completes the proof of Claim \ref{panda}.
\qed\\

\vskip10mm
\subsection{Related topics}\label{sub2}\hskip10mm\\%
Theorem \ref{Montgomery} is one of our main results. 

\begin{thm}\label{Montgomery}
If two virtual 1-knot diagrams $\alpha$ and $\beta$ are fiberwise equivalent, 
then $\alpha$ and $\beta$ are virtually parity equivalent. 
\end{thm}

\noindent{\bf Proof of Theorem \ref{Montgomery}.}
Theorem \ref{smooth} 
and Claim \ref{hirumeshi} imply Theorem \ref{Montgomery}.  \qed\\


It is known  that 
the usual trefoil knot diagram is not welded equivalent 
to the trivial knot diagram 
(see 
\cite{Kauffman, Kauffmanrw, Rourke, Satoh, J}).      
Hence these two diagrams are not also rotational welded equivalent. 
Hence we have the following. 

\begin{cla}\label{ice}
The converse of Theorem $\ref{Montgomery}$ is not true in general. 
\end{cla}

We have the following. 

\begin{cla}\label{cream}
The number of virtual crossing points of virtual 1-knot diagrams is not 
 an invariant of the fiberwise equivalence 
$($respectively, the rotational welded equivalence$)$  
in general.   
\end{cla}  
 
\begin{figure}
     \includegraphics[width=70mm]{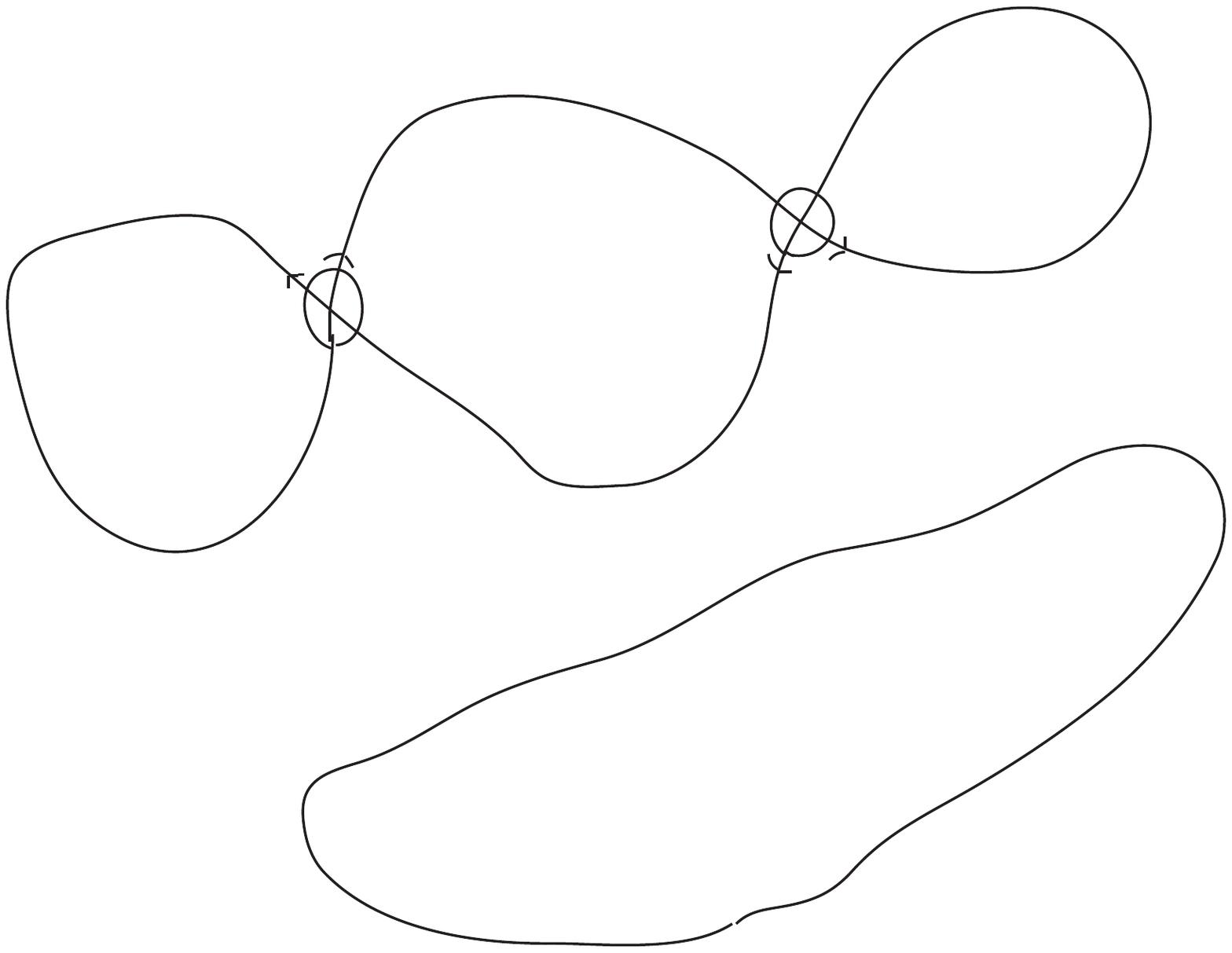}  
\caption{{\bf  
Two virtual knot diagrams which are rotational welded equivalent 
(respectively, fiberwise equivalent).
}}\label{oldWestVirginia}   
\end{figure}

\noindent{\bf Proof of Claim \ref{cream}.}
The two virtual knot diagrams in Figrue \ref{oldWestVirginia}
are fiberwise equivalent (respectively, rotational welded equivalent). \qed\\

We introduce a `weaker' equivalence relation than the fiberwise equivalence defined by 
Definition \ref{Nevada}.  
We want to replace `level preserving embedding of $S^1_b\x[0,1]$'  in 
Definition  \ref{Nevada}  
with an oriented compact surface 
which is not necessarily 
`level preserving embedding of $S^1_b\x[0,1]$',   
and loose a few conditions there. 
We prove 
in Theorem \ref{parityhozon} 
that 
this equivalence relation is equivalent to 
the virtual parity equivalence relation. 

\begin{defn}\label{sigma}  
Let $\alpha$ and $\beta$ be virtual 1-knot diagrams. 
We say that 
$\alpha$ and $\beta$ are {\it weakly fiberwise equivalent} 
if $\alpha$ and $\beta$ satisfy 
the following conditions.

\smallbreak\noindent$(1)$
There is 
a compact generic oriented surface $F$ with boundary whose boudary is a disjoint union of two circles, which is contained in $\R^2_b\x[0,1]$,  and 
$F$ is covered by Rourke's fibration.   
Note that thus there is a submanifold of $\R^2_b\x[0,1]\x\R^2_f$ which is diffeomorphic 
to $F\x S^1$.

\smallbreak\noindent$(2)$  
The in-out information of fiber circles gives  $\alpha$ and $\beta$ 
the information whether each multiple (respectively, branch) point is virtual or classical 
as in Theorem \ref{Montana}.   

\smallbreak\noindent$(3)$  
$\partial F$ is $\alpha$ and $\beta$. 
$F$ meets $\R^2_b\x\{0\}$ $($respectively,  $\R^2_b\x\{1\})$  
at $\alpha$ 
$($respectively, $\beta)$  
transversely.    

\smallbreak
If $F$ above is an annulus, 
we say that 
$\alpha$ and $\beta$ are {\it fiberwise cobordant}.  

\end{defn}

\begin{thm}\label{parityhozon}
Let $\alpha$ and $\beta$ be virtual 1-knot diagrams. 
 $\alpha$ and $\beta$ are weakly fiberwise equivalent 
if and only if 
 $\alpha$ and $\beta$ are virtually parity equivalent. 
\end{thm}

\begin{figure}
\includegraphics[width=80mm]{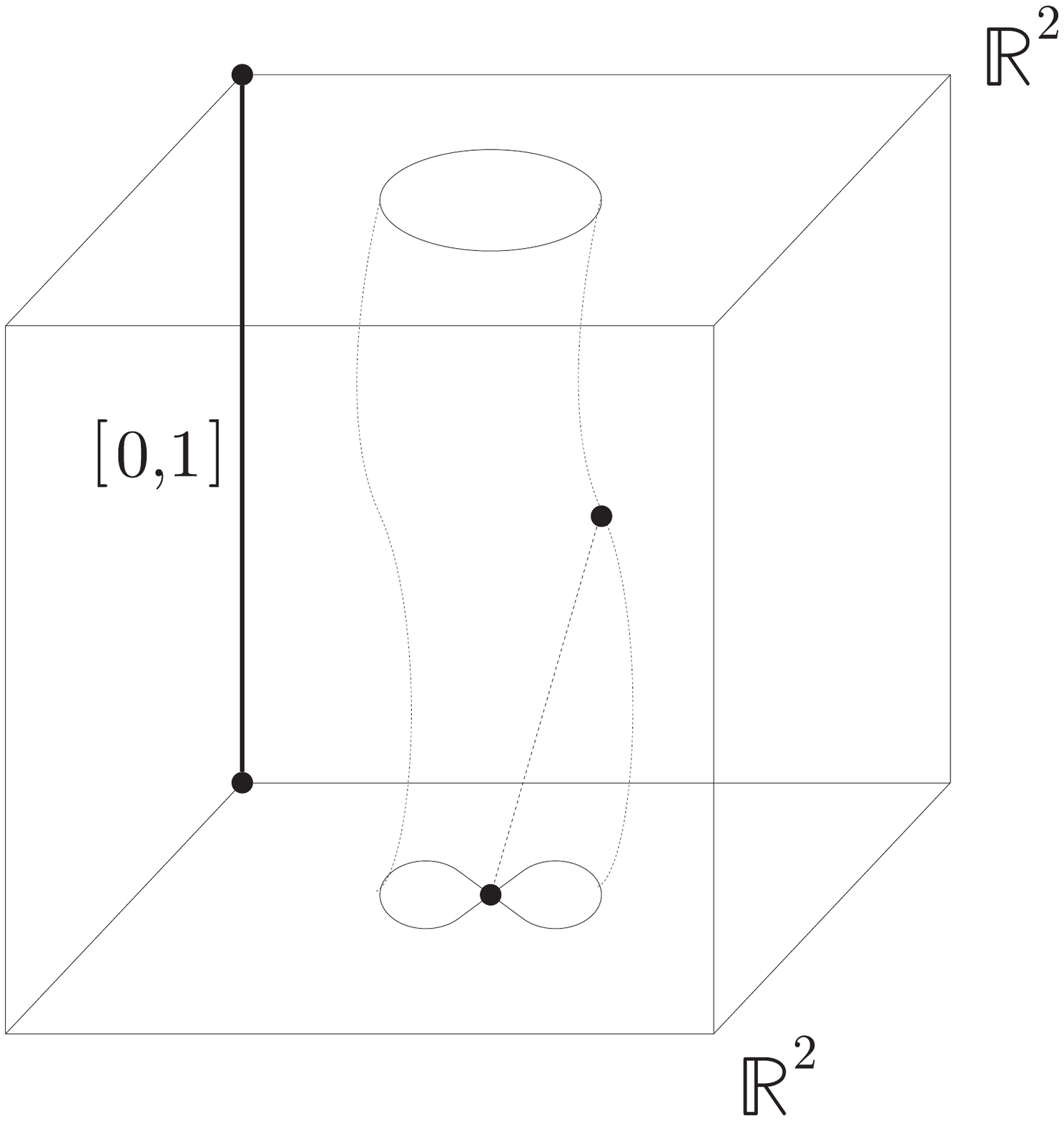}  

\caption{{\bf Let $\alpha$ be  the virtual figure $\infty$ knot diagram  and  $\beta$,   
the trivial 1-knot diagram. 
For example, $(\pi\circ g)(S_b^1\x[0,1]\x S^1_f)$ cannot be realized as drawn above, 
by Theorem \ref{Rmuri}.  
}}\label{Texas}   
\end{figure}

\noindent{\bf Proof of Theorem \ref{parityhozon}.} 
The `only if' part: 
We use  `reductio ad absurdum'.    
We suppose an assumption:  
 $\alpha$ and $\beta$ are not virtually parity equivalent. 
Take a generic surface which connects $\alpha$ and $\beta$  
as in Definition \ref{sigma}. 
Then this generic surface must have at least one 
virtual branch point because 
the union of 
$\alpha$ and $\beta$  
has 
an odd number of 
virtual crossing point.  
By Theorem \ref{Rmuri}, this generic surface never exists. 
(See Figure \ref{Texas}.)   
We arrived at a contradiction. 
Hence the above assumption is false and the `only if' part is true. 
\\

\begin{figure}\bigbreak
\includegraphics[width=140mm]{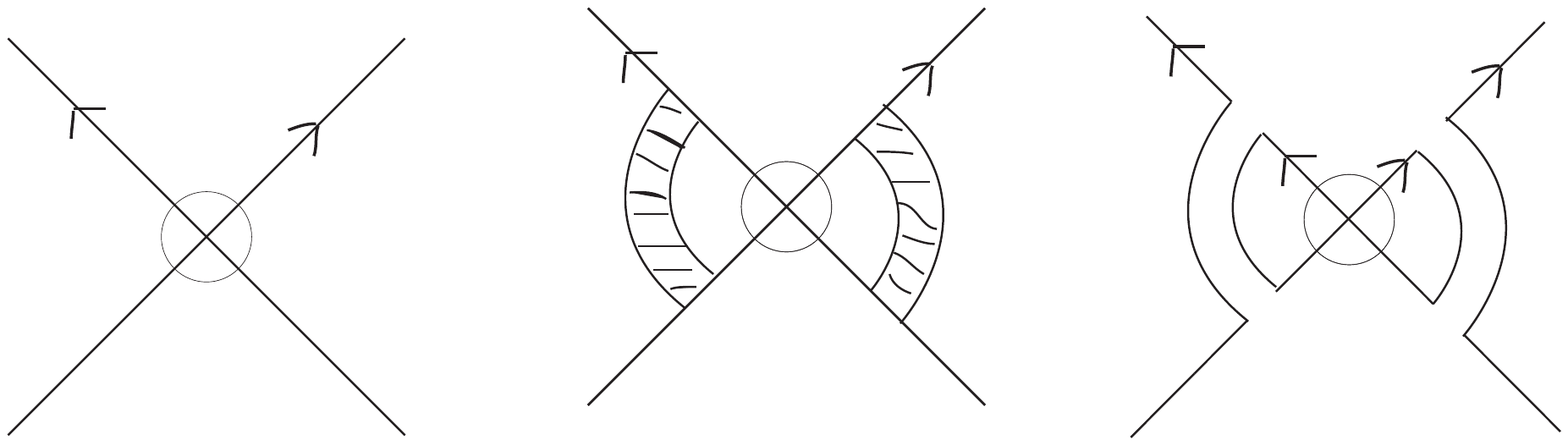}  
\vskip-70mm
\caption{{\bf Attaching bands.}\label{bands}}   
\end{figure}

\begin{figure}
\bigbreak
     \includegraphics[width=72mm]{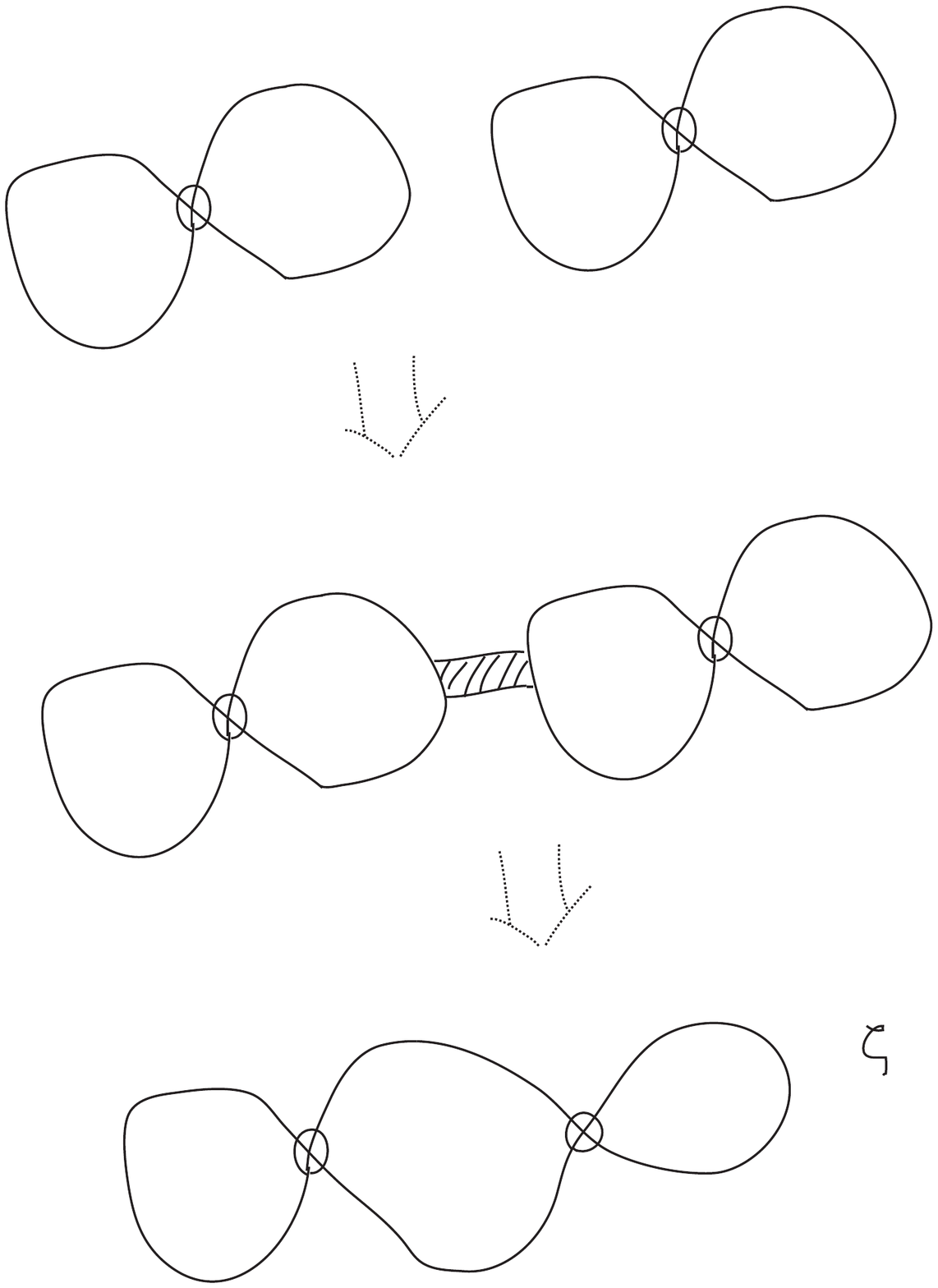}  
\smallbreak
\caption{{\bf  
A combination of two copies of the virtual figure $\infty$ knot diagram
}}\label{WestVirginia}   
\bigbreak  
\end{figure}

The `if' part:  
It suffices to prove that 
 $\alpha\amalg(-\beta)$ in $\R^2$, 
where $\amalg$ denote the disjoint union of the diagrams, 
is weakly fiberwise equivalent to the trivial 1-knot diagram. 
We can attach bands 
as drawn in Figure \ref{bands} 
so that the orientations of virtual knot diagrams and those of the bands are compatible.  
Thus $\alpha\amalg(-\beta)$ is 
weakly fiberwise equivalent to 
the disjoint union of 
nonnegative even  integer of  copies of the virtual figure $\infty$ knot and a classical link diagram. 
We can attach a band to two copies of the virtual figure $\infty$ knot diagram and 
 combine them as drawn in Figure \ref{WestVirginia}, 
so that the orientation of the band and those of the knot diagrams 
are compatible, 
and call the resultant diagram $\zeta$.   
Thus $\alpha\amalg(-\beta)$in $\R^2$ is 
weakly fiberwise equivalent to 
the disjoint union of a finite number of copies of $\zeta$. 
It is easy to prove that $\zeta$ is 
rotational welded equivalent to the trivial knot.  
Suppose that we obtain the $\mu$-component trivial 1-link diagram after that. 
Attach $\mu-1$ copies of 2-disc to   $(\mu-1)$ components of this trivial 1-link diagram.  
Hence $\alpha\amalg(-\beta)$ is weakly fiberwise equivalent to the trivial 1-knot diagram. 
Hence  $\alpha$ and $\beta$ are weakly fiberwise equivalent.  
This completes the proof of Theorem \ref{parityhozon}. \qed 

\begin{defn}\label{nono}  
We define 
the `{\it nonorientably weakly fiberwise equivalence}'.    
In Definition \ref{sigma} 
replace  `oriented surface' with  `non-orientable surface', 
and 
replace `weakly fiberwise equivalent' with 
`nonorientably weakly fiberwise equivalent'.    
\end{defn}

\begin{thm}\label{yuru}
Let $\alpha$ and $\beta$ be virtual 1-knot diagrams. 
 $\alpha$ and $\beta$ are nonorientably weakly fiberwise equivalent 
if and only if 
 $\alpha$ and $\beta$ are virtually parity equivalent. 
\end{thm}

\noindent{\bf Proof of Theorem \ref{yuru}.}
The `if' part: 
By Theorem \ref{parityhozon}, 
 $\alpha$ and $\beta$ are weakly fiberwise equivalent. 
It is trivial that 
if $\alpha$ and $\beta$ are weakly fiberwise equivalent,  
then $\alpha$ and $\beta$ are nonorientably weakly fiberwise equivalent.  
{\it Reason.} 
Take a generic oriented surface 
for $\alpha$ and $\beta$  
as in Definition \ref{sigma}. 
Take an immersed Klein bottle in $\R^2_b\x[0,1]$. 
Connect 
the generic oriented surface 
and 
the immersed Klein bottle 
by using an embedded 3-dimensional 1-handle in $\R^3$ 
such that 
the intersection of the 1-handle and 
the oriented surface (respectively, immersed Klein bottle) 
is  only the attaching part of the 1-handle. 
The resultant generic nonorientable surface implies that 
 $\alpha$ and $\beta$ are nonorientably weakly fiberwise equivalent. 
The proof of the `only if' part is 
the same as that of the `only if' part of Theorem \ref{parityhozon} 
if we replace the words `Definition \ref{sigma}' 
with `Definition \ref{nono}',   
and remove the sentence `(See Figure \ref{Texas}.)'.   
\qed

\bigbreak
Define {\it Whitney degree} of any virtual 1-knot diagram $\alpha$ to be 
Whitney degree which is defined by $\alpha$ 
when we regard $\alpha$ as an immersed oriented circle in $\R^2$.   
Two virtual 1-knot diagrams $\alpha$ and $\beta$ 
are said to be {\it Whitney parity equivalent}  
if the parity of Whitney degree of $\alpha$ is the same as that of $\beta$.   
Two virtual 1-knot diagrams $\alpha$ and $\beta$ 
are said to be {\it classically parity equivalent}  
if the parity of the classical crossing points of $\alpha$ is the same as that of $\beta$. 
The following holds. 
Let $\alpha$ and $\beta$ be virtual 1-knot diagrams which are rotational welded equivalent (respectively, fiberwise equivalent). (Note Theorem \ref{smooth}.) 
Then 
$\alpha$ and $\beta$ 
are classically parity equivalent 
if and only if 
$\alpha$ and $\beta$ 
are Whitney parity equivalent. \\ 

\noindent{\it Reason.} 
$\alpha$ and $\beta$ 
are Whitney parity equivalent 
if and only if 
 the number of the classical Reidemeister $I$ moves is even 
in a sequence of rotational welded moves which $\alpha$ is changed into $\beta$. 
Note that we cannot use  the virtual Reidemeister $I$ move by definition. 



\bigbreak

\begin{figure}
\bigbreak
     \includegraphics[width=120mm]{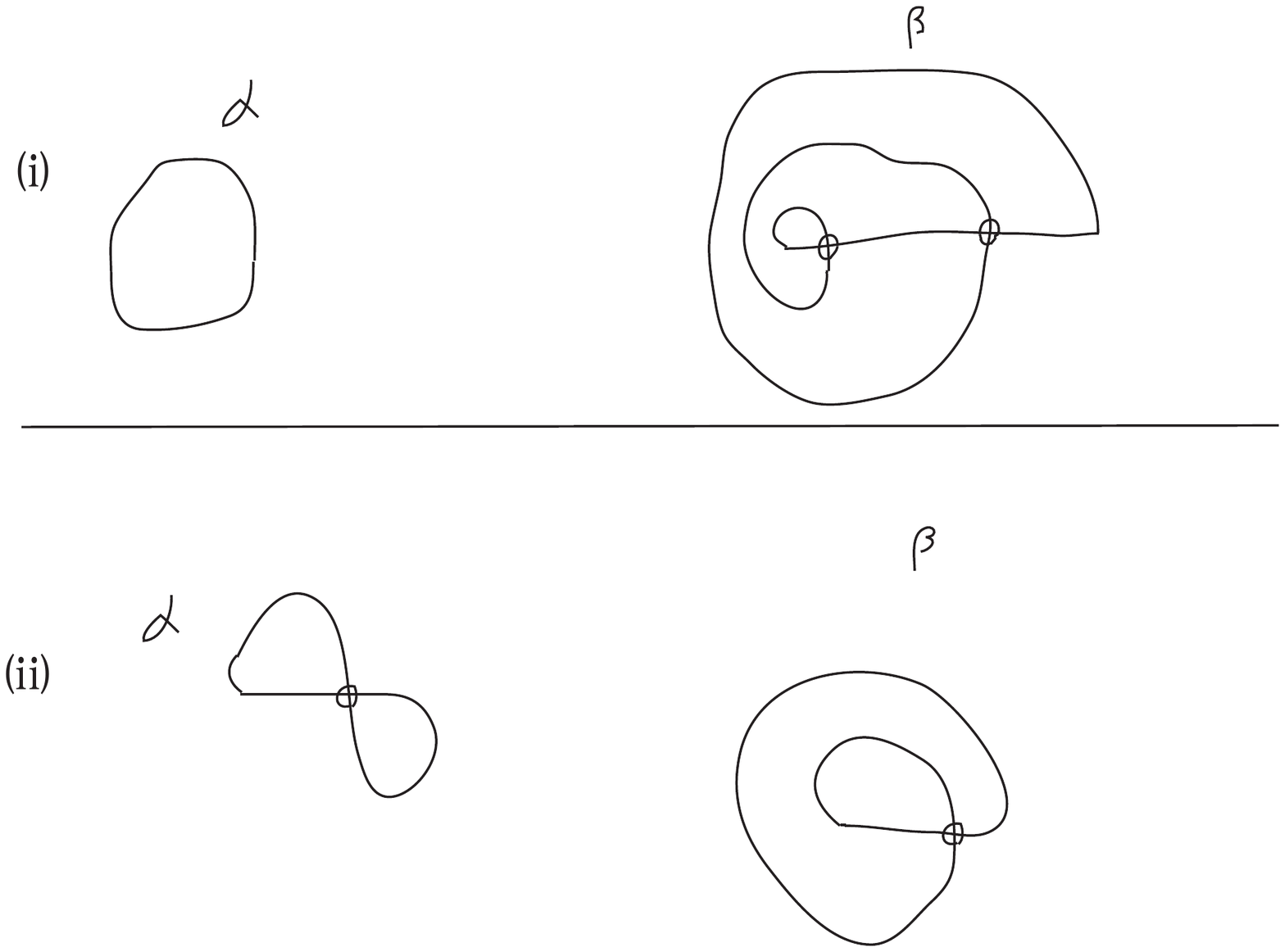}  
\smallbreak
\caption{{\bf  
Two pairs of virtual knot diagrams 
}}\label{Whitneydegree}   
\bigbreak  
\end{figure}

Some readers may ask the following question: 
Suppose that 
two virtual 1-knot diagrams $\alpha$ and $\beta$ 
do not have any classical crossing point 
and that Whitney degrees are different. 
Then is it valid that 
$\alpha$ and $\beta$ are not rotational welded equivalent?
The answer is negative. 
We show a counter example in Figure \ref{Whitneydegree}.(i) (respectively, \ref{Whitneydegree}.(ii)). 
The proof that each pair is rotational welded equivalent is left to the reader. 
\\

Two virtual 1-knot diagrams $\alpha$ and $\beta$ 
are said to be {\it mixed parity equivalent}  
if the parity of 
the sum of the classical and virtual crossing points of $\alpha$ is 
the same as 
that of $\beta$. 
The following holds. 
Let $\alpha$ and $\beta$ be virtual 1-knot diagrams which are welded equivalent. Then 
$\alpha$ and $\beta$ 
are mixed parity equivalent 
if and only if 
$\alpha$ and $\beta$ 
are Whitney parity equivalent.  \\

\noindent{\it Reason.} 
$\alpha$ and $\beta$ 
are Whitney parity equivalent 
if and only if 
the sum of the number of 
the classical and virtual Reidemeister $I$ moves 
is even 
in a sequence of welded moves which $\alpha$ is changed into $\beta$.


\bigbreak
\section{Virtual high dimensional knots}\label{vhigh}
\noindent 
See \cite{KauffmanOgasa, KauffmanOgasaII, KauffmanOgasaB, LevineOrr, Ogasa} 
for  codimension two  high dimensional knots. 
See \cite{Haefliger, Haefliger2, Levine} for high codimensional high dimensional knots. 
It is natural to  attempt to define virtual high dimensional knots and 
their one-dimensional-higher tubes. 
We could define $n$-dimensional virtual knots by using virtual $n$-knot diagrams in $\R^{n+1}$. 
We would make any virtual $n$-knot 
into a submanifold   
of
(a closed oriented $n$-manifold $M$) $\x[0,1]$ 
as we do in the virtual 1- and 2-dimensional cases.  
We want to make a bijection between 
the  set of such submanifolds and that of virtual $n$-knots. 
The 1-dimensional case is done (see Theorem \ref{vk} and Definition \ref{Jbase}). 
We should define a one-dimensional-higher tube   
as the spinning submanifold made from $K$ around $M$.  
Satoh's method makes no sense in the dimension greater than one. 
Rourke's way also makes non-sense by Theorem \ref{Rmuri}. 
Furthermore we must note 
that 
the $n$-dimensional case ($n\in\N-\{1,2\}$) of Theorems \ref{oh} and \ref{ohoh} 
does not necessarily hold 
in smooth category (respectively, PL category) 
because 
it is not trivial to produce an analogue of their proof 
by the following fact of \cite{Hudson}: 
There is an integer $p\geqq3$ and 
are two smooth (respectively, PL) $a$-dimensional submanifolds, $X$ and $Y$,  of $S^{a+p}$ which are diffeoomorphic (respectively, PL homeomorphic) each other 
but which are non-isotopic 
as smooth submanifolds (respectively, PL submanifolds). 
To complete these topics in this section is left to the readers as problems.

{\footnotesize \noindent{\bf Acknowledgment.} Kauffman's work was supported by the Laboratory of Topology and Dynamics, Novosibirsk State University
 (contract no.14.Y26.31.0025 with the Ministry of Education and Science of the Russian Federation).}

\bigbreak
\noindent
Louis H. Kauffman: 
Department of Mathematics, Statistics and Computer Science \\ 851 South Morgan Street    University of Illinois at Chicago
Chicago, Illinois 60607-7045, and\\ Department of Mechanics and Mathematics, Novosibirsk State University, Novosibirsk,  Russia\quad kauffman@uic.edu

\bigbreak\noindent
Eiji Ogasa:  Computer Science, Meijigakuin University, Yokohama, Kanagawa, 244-8539, Japan 
\quad pqr100pqr100@yahoo.co.jp  \quad
ogasa@mail1.meijigkakuin.ac.jp 

\bigbreak\noindent
Jonathan Schneider: 
Department of Mathematics, College of DuPage, 
425 Fawell Boulevard, Glen Ellyn,  Illinois, 60137, USA \quad
jschneider.math@gmail.com

\end{document}